\documentclass[11pt,reqno]{amsart}
\usepackage[T1]{fontenc}
\usepackage{lmodern}
\usepackage[a4paper,left=28mm,right=28mm,top=27mm,bottom=28mm]{geometry}
\usepackage{amsmath,amssymb,amsthm,mathtools,enumitem,microtype}
\usepackage{mathscinet}
\usepackage[colorlinks=true,linkcolor=blue,citecolor=blue,urlcolor=black]{hyperref}
\usepackage{bookmark}
\usepackage{tikz-cd}
\usepackage{booktabs,multicol,needspace}
\hypersetup{pdftitle={Symplectic structures on Lie algebras associated with acyclic quivers},pdfsubject={Combinatorial and recursive criteria for symplectic structures on path Lie algebras},pdfkeywords={acyclic quiver, path algebra, symplectic Lie algebra, forest, perfect matching}}
\setlist[enumerate]{label=\textup{(\roman*)},leftmargin=*,itemsep=2pt,topsep=4pt}
\newtheorem{theorem}{Theorem}[section]
\newtheorem{proposition}[theorem]{Proposition}
\newtheorem{lemma}[theorem]{Lemma}

\newtheorem{maintheorem}{Theorem}

\theoremstyle{definition}
\newtheorem{definition}[theorem]{Definition}
\newtheorem{example}[theorem]{Example}
\theoremstyle{remark}
\newtheorem{remark}[theorem]{Remark}
\newcommand{\kk}{\Bbbk}
\newcommand{\nq}{\mathfrak n_Q}
\newcommand{\sq}{\mathfrak s_Q}
\newcommand{\hq}{\mathfrak h_Q}
\newcommand{\pp}{\mathcal P(Q)}
\newcommand{\pr}[1]{\mathcal P_{#1}(Q)}
\newcommand{\Span}{\operatorname{span}}
\newcommand{\rank}{\operatorname{rank}}
\newcommand{\rad}{\operatorname{rad}}
\newcommand{\Pf}{\operatorname{Pf}}

\newcommand{\ad}{\operatorname{ad}}
\newcommand{\PM}{\operatorname{PM}}
\newcommand{\Fmax}{F_{\max}(Q)}
\newcommand{\Acal}{\mathcal A(Q)}
\newcommand{\Rcal}{\mathcal R(Q)}
\newcommand{\GR}{G_R(Q)}
\newcommand{\KQ}{\mathcal K_Q(\mathbf X)}
\newcommand{\FQ}{\kk(\mathbf X)}
\newcommand{\eps}{\varepsilon}
\newcommand{\id}{\operatorname{id}}
\title[Symplectic Lie groups and quivers]
{Symplectic Lie groups and nilmanifolds associated with acyclic quivers}

\hypersetup{
  pdftitle={Symplectic Lie groups and nilmanifolds associated with acyclic quivers}
}

\author{Yihao Zheng}
\thanks{%
\begin{tabular}[t]{@{}l@{}}
School of Mathematical Sciences, Fudan University, Shanghai 200433, China\\
Email: \href{mailto:yhzheng24@m.fudan.edu.cn}{yhzheng24@m.fudan.edu.cn}
\end{tabular}}

\date{}
\begin{document}

\begin{abstract}
We give combinatorial criteria for the existence of left-invariant symplectic forms on the positive and full path Lie groups of finite acyclic quivers. For the nilpotent positive path Lie groups, the criterion also determines whether their compact quotients admit symplectic forms and forces the nilpotency step to be at most three. For the completely solvable full path Lie groups, we give a direct combinatorial criterion and an equivalent recursive criterion. In the symplectic case, the recursive procedure is realized by successive Hamiltonian reductions ending at an abelian symplectic Lie group. We characterize the existence of Lagrangian normal subgroups for a given form in both families, with applications to Lagrangian torus fibrations of compact nilmanifolds. We also describe the cohomology rings of the compact symplectic nilmanifolds and determine all Lefschetz ranks.
\end{abstract}
\maketitle

\section{Introduction}

Lie groups with left-invariant symplectic structures and their
compact quotients are important objects of study in geometry.
The Lie algebraic description of left-invariant differential
forms facilitates explicit constructions and computations,
while compact quotients provide concrete examples for studying
global geometric and topological properties of symplectic
manifolds.
These manifolds have provided fundamental examples and
contributed to the study of generalized complex geometry, symplectic Hodge theory, de Rham cohomology, rational homotopy
theory, flat affine geometry, Lagrangian foliations, 
symplectic reduction, geometric quantization and representation theory; see, for example,
  \cite{AngellaKasuya}, \cite{BC}, \cite{CavalcantiGualtieri},
\cite{Hasegawa}, \cite{KirwinUribe}, \cite{Medina},
\cite{MedinaRevoy}, \cite{Nomizu}.

Compact nilmanifolds are particularly useful in this setting
because their topology can be studied through left-invariant
differential forms on the corresponding nilpotent Lie groups.
Compact nilmanifolds also provide a striking illustration of
the distinction between symplectic and K\"ahler geometry.
Thurston \cite{Thurston} constructs a compact symplectic
nilmanifold admitting no K\"ahler structure, while
Benson and Gordon \cite{BG} prove that a compact symplectic
nilmanifold satisfying the hard Lefschetz property must
be diffeomorphic to a torus.
These results motivate more precise questions about which
compact nilmanifolds admit symplectic forms, whether a given
form admits a Lagrangian fibration, and how multiplication
by its class acts on cohomology.

Finite acyclic quivers give families of Lie groups in which
these questions can be studied in arbitrary dimension.
Mizoguchi and Tamaru \cite{MT} construct nilpotent Lie algebras
from the non-trivial paths of such quivers and prove that the
corresponding simply connected Lie groups all admit
left-invariant Ricci solitons.
They also propose the study of other geometric structures
on these Lie groups.
Since the  families of Lie groups introduced above are defined by combinatorial data,
it is natural to ask which geometric properties are determined
by those data and which depend on the choice of symplectic form.
We determine which of these Lie groups admit
left-invariant symplectic forms and study the 
geometry of these groups and their compact quotients.
We also consider the completely solvable groups obtained by
including the trivial paths, and relate their symplectic
structures to explicit Hamiltonian reductions.

Let $Q$ be a finite acyclic quiver, with vertex set $Q_0$ and
arrow set $Q_1$.
Write $\mathbb RQ$ for its real path algebra and
$\mathbb RQ_{\ge1}$ for the ideal spanned by the non-trivial paths.
The positive path Lie group and the full path Lie group are
$N_Q=1+\mathbb RQ_{\ge1}$ and
$S_Q=(\mathbb RQ)^\times_0$, respectively, where the subscript
denotes the identity component.
Their Lie algebras are denoted by $\nq$ and $\sq$.
Both Lie groups are simply connected, $N_Q$ is nilpotent, and
$S_Q$ is completely solvable.
Moreover, $\Gamma_Q=1+\mathbb ZQ_{\ge1}$ is a cocompact
lattice in $N_Q$, whereas $S_Q$ has no lattice when $Q$
has an arrow.
These are the main objects of study in this paper.


Pouseele and Tirao \cite{PT} characterize symplectic existence
for compact nilmanifolds associated with graphs.
Mizoguchi \cite{Mizoguchi} shows that the two-step
positive path Lie algebras are graph Lie algebras and
determines exactly which of them admit symplectic structures.
Our first result goes beyond the two-step setting and
characterizes all finite acyclic quivers whose positive
path Lie groups admit left-invariant symplectic forms.


\begin{maintheorem}[Theorem~\ref{thm:nilclassification}
and Theorem~\ref{thm:geometricexistence}]\label{thm:intropositive}
The positive path Lie group $N_Q$ admits a left-invariant
symplectic form if and only if $m(Q)$ is even and one of
the following conditions holds.
\begin{enumerate}
\item $L(Q)\le2$ and
$d^-(v)d^+(v)\le d^-(v)+d^+(v)$
for every vertex $v$ of $Q$.
\item $L(Q)=3$, every active pair $(u,v)$ satisfies
$d^-(u)=d^+(v)=1$ and has even multiplicity $m_{uv}$,
and every connected component $H$ of $G_R(Q)$ satisfies
$|E(H)|\le|V(H)|$.
\end{enumerate}
For every cocompact lattice $\Gamma\subset N_Q$, these
conditions are also equivalent to the existence of a
symplectic form on $\Gamma\backslash N_Q$.
\end{maintheorem}

Thus symplecticity forces the nilpotency step to be at most
three, and when $L(Q)\ge4$, no compact quotient of $N_Q$
admits any symplectic form.
The same combinatorial conditions also characterize holomorphic
symplectic structures on the compact complex parallelizable
quotient
$\bigl(1+\mathbb Z[i]Q_{\ge1}\bigr)\backslash
\bigl(1+\mathbb CQ_{\ge1}\bigr)$.

For the full path Lie groups $S_Q$, symplectic existence is
governed by the underlying undirected multigraph $U(Q)$,
in which parallel arrows remain distinct edges.

\begin{maintheorem}[Theorem~\ref{thm:fullclassification} and
Theorem~\ref{thm:solvableexistence}]\label{thm:introfull}
The full path Lie group $S_Q$ admits a left-invariant symplectic
form if and only if $n(Q)+m(Q)$ is even and the following
conditions hold.
\begin{enumerate}
\item $U(Q)$ is a forest.
\item Every internal vertex has degree at most three in $Q$.
\item Every connected component of $U(Q)[I(Q)]$ contains at
most one vertex of degree three in $Q$, and any such vertex
is adjacent to a source or sink of $Q$.
\end{enumerate}
\end{maintheorem}

The criterion can be checked in $O(|Q_0|+|Q_1|)$ operations,
without enumerating the non-trivial paths.
In this family there is no upper bound on the path length or
the solvability step of a Lie group admitting a left-invariant
symplectic form.

The direct existence criterion for $S_Q$ has an equivalent
recursive form whose steps admit a geometric interpretation.
We use symplectic reduction, introduced by Marsden and
Weinstein \cite{MW} and studied for symplectic Lie groups
by Medina \cite{Medina}, to relate these groups to
abelian symplectic Lie groups.

We introduce terminal-path and source reductions  of $Q$. These operations
are defined without a symplectic hypothesis and yield the recursive  existence criterion below. For a given left-invariant symplectic form, the same
operations are realized by Hamiltonian reductions within
the class of full path Lie groups.
They also construct connected closed Lagrangian subgroups.
The existence of such subgroups follows more generally
from the work of Baues and Cort\'es \cite{BC} on completely
solvable symplectic Lie groups; here the construction
is determined by the successive quiver reductions.

\begin{maintheorem}[Theorem~\ref{thm:reduction}, 
Theorem~\ref{thm:recursive}, 
Theorem~\ref{thm:hamiltonianreduction}, and  
Theorem~\ref{thm:lagrangiansubgroup}]\label{thm:introreduction}
Starting from $Q$, apply terminal-path or source reductions
until neither is possible, and let $Q_\infty$ be the final
quiver.
\begin{enumerate}
\item The group $S_Q$ admits a left-invariant symplectic form
if and only if $Q_\infty$ has no arrows and an even number
of vertices.
When these conditions hold, reversing the reductions
constructs a left-invariant symplectic form on $S_Q$
from one on the final abelian group.
\item For every left-invariant symplectic form $\Omega$ on
$S_Q$, this sequence is realized by Hamiltonian reductions
by connected closed isotropic normal subgroups.
The successive reduced symplectic Lie groups are isomorphic
to the full path Lie groups of the reduced quivers equipped
with left-invariant symplectic forms, and the final group
is abelian.
A connected closed Lagrangian subgroup of $(S_Q,\Omega)$
can be constructed along the same sequence.
\end{enumerate}
\end{maintheorem}

Theorem~\ref{thm:introreduction} constructs Lagrangian
subgroups for every left-invariant symplectic form on $S_Q$.
We next determine when a Lagrangian subgroup can be chosen
to be normal and study the corresponding question for $N_Q$.
Baues and Cort\'es \cite{BC} develop a structure theory of
symplectic Lie groups based on isotropic normal subgroups
and symplectic reduction.
In particular, a Lagrangian normal subgroup gives a quotient
with a natural flat torsion-free connection, but such a
subgroup need not exist even for nilpotent symplectic Lie groups.
In dimension eight, A\"it Aissa and Mansouri \cite{AM}
classify the symplectic nilpotent Lie groups admitting a
Lagrangian normal subgroup.

We characterize the existence of Lagrangian normal subgroups
for a given left-invariant symplectic form in both families.
For full path Lie groups, we also determine the largest
possible dimension of a connected closed isotropic normal
subgroup.
For positive path Lie groups, the criterion leads to
Lagrangian torus fibrations of compact quotients.


\begin{maintheorem}[Theorem~\ref{thm:nilidealcriterion}
and Theorem~\ref{thm:fullisotropicrank}]
\label{thm:introlagrangian}
\begin{enumerate}
\item
Let $\Omega$ be a left-invariant symplectic form on $N_Q$,
with value $\omega$ at the identity.
The group $(N_Q,\Omega)$ has a connected closed Lagrangian
normal subgroup if and only if $G_A(Q)$ has a vertex cover
which is independent in $F_\omega$.

\item
Let $\Omega$ be a left-invariant symplectic form on $S_Q$,
with value $\omega$ at the identity, and write
$2r=\rank(\omega|_{Z(\sq)})$.
Then
\[
 \max_H\dim H=a(Q)+c(Q)-r,
\]
where $H$ ranges over the connected closed isotropic normal
subgroups of $(S_Q,\Omega)$.
Consequently, $(S_Q,\Omega)$ has a connected closed
Lagrangian normal subgroup if and only if
\[
 a(Q)+c(Q)-r=\tfrac12\bigl(n(Q)+m(Q)\bigr).
\]
\end{enumerate}
\end{maintheorem}



For connected $Q$, both this maximum and the existence
of a Lagrangian normal subgroup are independent of
the symplectic form.
For the upper-triangular groups with positive diagonal
entries, we construct examples showing that the difference
between half the group dimension and this maximum is
unbounded.

For a symplectic positive path Lie group $N_Q$, however,
a connected closed Lagrangian normal subgroup exists
for every left-invariant symplectic form if and only if
$G_A(Q)$ has matching number at most one.
The criterion in Theorem~\ref{thm:introlagrangian}\textup{(i)}
also governs Lagrangian torus fibrations of compact quotients.
A rational connected closed Lagrangian normal subgroup
of $N_Q$ is abelian and, by Mal\cprime cev's rational
subgroup theorem \cite{Malcev}, \cite{Raghunathan}, induces
a Lagrangian torus fibration of $\Gamma_Q\backslash N_Q$.
For a left-invariant symplectic form whose value at the
identity is rational in the path basis, this criterion
is equivalent to the existence of such a fibration arising
from a rational normal subgroup.
Every symplectic positive path Lie group $N_Q$ admits
a left-invariant symplectic form whose value at the identity
is rational in the path basis and for which such a
fibration exists.

For a given rational symplectic form, however, such a
fibration need not exist.
We construct a twenty-dimensional positive path Lie group
with two left-invariant symplectic forms, one admitting
a connected closed Lagrangian normal subgroup and the
other not.
Twenty is the smallest dimension in which this dependence
can occur among positive path Lie groups.
Both forms have rational values at the identity in the
path basis and descend to the same compact quotient,
on which such a fibration exists for the first induced
form but not for the second.




On a compact quotient of $N_Q$, multiplication by the
symplectic class provides another way to study the dependence
of the geometry on the form.
Mathieu \cite{Mathieu} proves that the hard Lefschetz property
is equivalent to the existence of a symplectically harmonic
representative in every de Rham cohomology class.
Ugarte and Villacampa \cite{UV} study the relation between
Lefschetz maps and symplectic cohomologies, including invariants
which vary with the symplectic form.
For non-toral compact nilmanifolds, the theorem of Benson and
Gordon \cite{BG} already excludes hard Lefschetz.
The ranks in individual degrees nevertheless contain finer
information, and one can ask which ranks are determined by
the manifold and which distinguish symplectic classes.

To study these maps in all degrees, we first describe
the cohomology ring.
Armstrong, Cairns and Jessup \cite{ACJ} compute the Betti
numbers of graph Lie algebras associated with star graphs,
and Aldi et al.~\cite{Aldi} study the cohomology of graph
Lie algebras more generally.
For the compact symplectic quotients considered here,
we determine all Betti numbers and describe the cup product.
We then give an explicit decomposition into graded Lefschetz
strings for multiplication by each symplectic class.
This determines the Lefschetz ranks in every degree.

\begin{maintheorem}[Theorem~\ref{thm:firstlefschetz}
and Theorem~\ref{thm:allstrings}]\label{thm:introlefschetz}
Let $M=\Gamma\backslash N_Q$ be a compact nilmanifold of
dimension $2d>0$, and let $\omega$ be a symplectic form on $M$.
Write $L_\omega$ for multiplication by $[\omega]$.

The construction in Theorem~\ref{thm:allstrings},
together with Appendix~\ref{app:square}, gives the starting 
degrees, lengths and multiplicities of all graded Lefschetz
strings of $H^*(M;\mathbb R)$.
It determines the ranks of all maps
\[
 L_\omega^u\colon H^k(M;\mathbb R)
 \longrightarrow H^{k+2u}(M;\mathbb R),
 \qquad k,u\ge0.
\]
In degree one,
\[
 \rank\bigl(
 L_\omega^{d-1}\colon H^1(M;\mathbb R)
 \longrightarrow H^{2d-1}(M;\mathbb R)
 \bigr)
 =|Q_1|-\nu(Q),
\]
where $\nu(Q)$ is the number of maximal paths in $Q$
of length at least two.
\end{maintheorem}

The first rank is therefore independent of the symplectic
class, and its kernel has dimension $\nu(Q)$.
Higher Lefschetz ranks, however, can distinguish symplectic
classes on the same compact nilmanifold, as we show by
an explicit example.


We also compare the two symplectic existence criteria,
give their explicit forms for $L(Q)\le2$, and characterize
simultaneous existence on $N_Q$ and $S_Q$.
We explain the different roles of parallel arrows and show
that, when $N_Q$ is nontrivial, every left-invariant symplectic
form on $S_Q$ restricts degenerately to $N_Q$.
We also give a sufficient condition in terms of splitting
graphs for an even-dimensional full path Lie group to admit
a left-invariant symplectic form, and exhibit the limitations
of dimension, matching and forest conditions.
In either family, whenever a left-invariant symplectic form
exists, we construct one whose matrix at the identity has
integer entries and determinant one in the path basis.

The paper is organized as follows.
In Section~\ref{sec:prelim}, we introduce the conventions
and preliminary results.
In Section~\ref{sec:nil}, we classify the quivers for which
$\nq$ admits a symplectic form.
In Section~\ref{sec:groups}, we study positive path Lie groups
and their compact quotients, including Lagrangian normal
subgroups and Lagrangian torus fibrations.
In Section~\ref{sec:full}, we establish direct and recursive
criteria for the existence of symplectic forms on $\sq$.
In Section~\ref{sec:solvablegroups}, we study Hamiltonian
reductions, Lagrangian subgroups and isotropic normal subgroups
of full path Lie groups, and compute the algebraic de Rham
cohomology of the complex unit group $(\mathbb CQ)^\times$.
In Section~\ref{sec:nilcohomlefschetz}, we describe the cohomology
and Lefschetz maps of the compact nilmanifolds.
In Section~\ref{sec:comparison}, we compare the two existence
criteria and derive further consequences.
In the appendix, we give the cohomology and Lefschetz
calculations for the residual component $K_{2,2}$.

\section{Preliminaries and general tools}\label{sec:prelim}

We introduce the two path Lie algebras and their associated
Lie groups, and collect the preliminaries for the
classification and geometric results that follow. Throughout, $\kk$ is a field of characteristic zero, and $Q$ is a finite acyclic quiver. Unless otherwise stated, all vector spaces and Lie algebras
are over $\kk$.

\subsection{Quivers, path Lie algebras and Lie groups}
We introduce the positive and full path Lie algebras
and their associated real Lie groups, and establish the
structural and lattice properties needed below.
For background on quivers and path algebras, see
\cite{ASS}, \cite{DW}, and \cite{Schiffler}; for lattices in Lie groups, see
\cite{CorwinGreenleaf} and \cite{Raghunathan}.

\begin{definition}\label{def:quiverpaths}
A quiver is a quadruple $Q=(Q_0,Q_1,s,t)$, where $Q_0$
is the vertex set, $Q_1$ is the arrow set, and
$s,t:Q_1\to Q_0$ are the source and target maps.
Distinct arrows with the same starting and terminal vertices
are parallel.
A non-trivial path is a sequence $p=a_1\cdots a_r$,
with $r\ge1$ and $t(a_i)=s(a_{i+1})$ for $1\le i<r$.
Its starting vertex is $s(p)=s(a_1)$, its terminal vertex
is $t(p)=t(a_r)$, and its length is $\ell(p)=r$.
For every vertex $v\in Q_0$, let $e_v$ be the trivial path
at $v$, with length zero and $s(e_v)=t(e_v)=v$.
A quiver is acyclic if $s(p)\ne t(p)$
for every non-trivial path $p$.
\end{definition}

\begin{definition}\label{def:downward}
A set $D\subseteq Q_0$ is downward closed if $u\in D$
whenever $v\in D$ and a path from $u$ to $v$ exists.
\end{definition}

Let $\pr r$ be the set of paths of length $r$, and let
$\pp=\bigcup_{r\ge1}\pr r$.
Write $L(Q)=\max\{\ell(p)\mid p\in\pp\}$, with $L(Q)=0$
if $Q$ has no arrows.
Let $n(Q)=|Q_0|$ and $m(Q)=|\pp|$.
The indegree $d^-(v)$ and outdegree $d^+(v)$ count,
respectively, arrows ending and starting at $v$,
with parallel arrows counted separately.
A vertex of indegree zero is a source,
and a vertex of outdegree zero is a sink.
A vertex is internal if both degrees are positive.
Let $I(Q)$ be the set of internal vertices.
Write $d(v)=d^-(v)+d^+(v)$ for the total degree of $v$ in $Q$,
and call $(d^-(v),d^+(v))$ its degree pair.

A non-trivial path is maximal if it starts at a source
and ends at a sink.
Let $\mathcal M(Q)$ be the set of maximal non-trivial paths.
Every non-trivial path extends to a maximal path.
A topological ordering of \(Q\) is a bijection \(\tau:Q_0\to\{1,\cdots,|Q_0|\}\) satisfying \(\tau(s(a))<\tau(t(a))\) for every arrow \(a\in Q_1\). Such an ordering exists for every finite acyclic quiver, since one can repeatedly remove a source and assign consecutive labels in this process.

\begin{definition}\label{def:pathlie}
The path algebra $\kk Q$ has  all paths as a basis,
with multiplication given by concatenation when
$t(p)=s(q)$ and zero otherwise.
The subspace $\kk Q_{\ge1}=\Span\pp$ is a two-sided ideal.

Let $\nq=\kk Q_{\ge1}$, $\sq=\kk Q$, and
$\hq=\Span\{e_v\mid v\in Q_0\}$, each equipped with the
Lie bracket $[x,y]=xy-yx$. $\nq$ is called the positive path Lie algebra of $Q$, and $\sq$ is called the full path Lie algebra of $Q$.
In particular, $[e_u,e_v]=0$ for all $u,v\in Q_0$,
so $\hq$ is abelian.
\end{definition}

For $p,q\in\pp$,
\begin{equation*}
 [p,q]=
 \begin{cases}
 pq,&t(p)=s(q),\\
 -qp,&t(q)=s(p),\\
 0,&\text{otherwise},
 \end{cases}
 \qquad
 [e_v,p]=(\delta_{v,s(p)}-\delta_{v,t(p)})p,
\end{equation*}
where $\delta$ is the Kronecker delta.
Let $\eps_v$ be dual to $e_v$ on $\hq$, and let
$\alpha_p=\eps_{s(p)}-\eps_{t(p)}$.
Then $[h,p]=\alpha_p(h)p$ for $h\in\hq$.
We call $\alpha_p$ the weight of $p$.

\begin{definition}
For a Lie algebra $\mathfrak g$, its center is
$Z(\mathfrak g)=\{z\in\mathfrak g\mid
[z,x]=0\text{ for all }x\in\mathfrak g\}$.
Let $C^1\mathfrak g=\mathfrak g$ and
$C^{r+1}\mathfrak g=[\mathfrak g,C^r\mathfrak g]$.
We call $\mathfrak g$ nilpotent of step $r$ if
$C^r\mathfrak g\ne0$ and $C^{r+1}\mathfrak g=0$.
Similarly, let $\mathfrak g^{(0)}=\mathfrak g$ and
$\mathfrak g^{(r+1)}=[\mathfrak g^{(r)},\mathfrak g^{(r)}]$.
We call $\mathfrak g$ solvable of step $r$ if
$\mathfrak g^{(r-1)}\ne0$ and $\mathfrak g^{(r)}=0$,
where $r\ge1$.
\end{definition}

\begin{proposition}\label{prop:structure}
If $L(Q)=\ell\ge1$, then $\nq$ is a Lie algebra nilpotent of step $\ell$ and $\sq$ is a Lie algebra solvable of step
$\lfloor\log_2\ell\rfloor+2$.
If $L(Q)=0$, then $\nq=0$ and $\sq$ is abelian.
Moreover, writing $\mathfrak n_Q^r=\Span\pr r$ for $r\ge1$, we have
$[\nq,\nq]=\bigoplus_{r\ge2}\mathfrak n_Q^r$,
$\sq=\hq\ltimes\nq$,
$[\sq,\sq]=\nq$, and
$Z(\nq)=\Span\mathcal M(Q)$.
\end{proposition}

\begin{proof}
It is easy to check that
$[\mathfrak n_Q^i,\mathfrak n_Q^j]=\mathfrak n_Q^{i+j}$
for $i,j\ge1$,
$[\nq,\nq]=\bigoplus_{r\ge2}\mathfrak n_Q^r$,
$\sq=\hq\ltimes\nq$, and $[\sq,\sq]=\nq$.

Suppose $L(Q)=\ell\ge1$.
Induction gives
$C^r\nq=\bigoplus_{j\ge r}\mathfrak n_Q^j$ for $r\ge1$.
Since $\mathfrak n_Q^\ell\ne0$ and
$\mathfrak n_Q^j=0$ for $j>\ell$,
we have $C^\ell\nq\ne0$ and $C^{\ell+1}\nq=0$.
Thus $\nq$ is nilpotent of step $\ell$.

Similarly, since $[\sq,\sq]=\nq$, induction gives
$\sq^{(r)}=\bigoplus_{j\ge2^{r-1}}\mathfrak n_Q^j$
for $r\ge1$.
Let $d=\lfloor\log_2\ell\rfloor+2$.
Then $2^{d-2}\le\ell<2^{d-1}$, so
$\sq^{(d-1)}\ne0$ and $\sq^{(d)}=0$.
Hence $\sq$ is solvable of step $d$.

Every maximal path belongs to $Z(\nq)$.
Conversely, if a central element $c$ has a nonzero
coefficient at a non-maximal path $p$, choose an arrow
$a$ extending $p$.
It is easy to check that $[c,a]\ne0$, a contradiction.
Therefore $Z(\nq)=\Span\mathcal M(Q)$.
\end{proof}

The opposite quiver $Q^{\mathrm{op}}$ is obtained by reversing
every arrow.
For each path $p$, let $p^{\mathrm{op}}$ be the reversed path.
The map $p\mapsto p^{\mathrm{op}}$ defines an anti-isomorphism
$\kk Q\to\kk Q^{\mathrm{op}}$,
so $p\mapsto-p^{\mathrm{op}}$ defines Lie algebra isomorphisms
$\nq\cong\mathfrak n_{Q^{\mathrm{op}}}$ and
$\sq\cong\mathfrak s_{Q^{\mathrm{op}}}$.
We use these isomorphisms to pass between statements
and their arrow-reversed versions.

For the remainder of this subsection, take $\kk=\mathbb R$
and consider the Lie groups associated with $\nq$ and $\sq$.

\begin{definition}
Let $1=\sum_{v\in Q_0}e_v$ be the identity element of
$\mathbb RQ$.
Let $N_Q=1+\mathbb RQ_{\ge1}$, with multiplication
$(1+x)(1+y)=1+x+y+xy$. Call it the positive path Lie group  of $Q$.
The identity element is $1$, and
$(1+x)^{-1}=\sum_{j=0}^{L(Q)}(-x)^j$.
The bijection $x\mapsto1+x$ from $\mathbb RQ_{\ge1}$ to $N_Q$
endows $N_Q$ with a smooth manifold structure for which
multiplication and inversion are smooth, making $N_Q$
a real Lie group.

Let $(\mathbb RQ)^\times$ be the group of invertible
elements of $\mathbb RQ$, with multiplication induced
by $\mathbb RQ$ and the topology
inherited from the real vector space $\mathbb RQ$.
 Let 
$S_Q=(\mathbb RQ)^\times_0$ be  the connected component
of $(\mathbb RQ)^\times$ containing $1$. We call $S_Q$ the full path Lie group of $Q$.
Explicitly, its elements are precisely 
$\sum_{v\in Q_0}t_ve_v+\sum_{p\in\pp}c_pp$,
where $t_v>0$ and $c_p\in\mathbb R$.
As an open subset of $\mathbb RQ$, the group $S_Q$
carries the inherited smooth structure and is a real
Lie group.
\end{definition}

\begin{definition}
For subgroups $H$ and $K$ of a  group $G$, let $[H,K]$
denote the subgroup generated by the commutators
$hkh^{-1}k^{-1}$, where $h\in H$ and $k\in K$.
 Now suppose that $G$ is a Lie group, let $C^1G=G$ and
$C^{r+1}G=[G,C^rG]$.
We call $G$ nilpotent of step $r$ if
$C^rG\ne\{1\}$ and $C^{r+1}G=\{1\}$.
Similarly, let $G^{(0)}=G$ and
$G^{(r+1)}=[G^{(r)},G^{(r)}]$.
We call $G$ solvable of step $r$ if
$G^{(r-1)}\ne\{1\}$ and $G^{(r)}=\{1\}$,
where $r\ge1$.
A real Lie group is completely solvable if every adjoint
operator of its Lie algebra has only real eigenvalues.
A discrete subgroup $\Gamma$ of a Lie group $G$ is a lattice
if $\Gamma\backslash G$ carries a finite nonzero
$G$-invariant Borel measure.
A lattice $\Gamma$ is cocompact if
$\Gamma\backslash G$ is compact.
\end{definition}

\begin{proposition}\label{prop:pathgroups}
The Lie group $N_Q$ is simply connected and nilpotent, with Lie
algebra $\nq$.
If $L(Q)=\ell\ge1$, then $N_Q$ is nilpotent of step $\ell$.
If $L(Q)=0$, then $N_Q=\{1\}$.
Moreover, $\Gamma_Q=1+\mathbb ZQ_{\ge1}$ is a cocompact
lattice in $N_Q$, where $\mathbb ZQ_{\ge1}$ denotes
the integer span of the non-trivial paths.
\end{proposition}

\begin{proof}
The map $x\mapsto1+x$ identifies $N_Q$ diffeomorphically
with the real vector space $\mathbb RQ_{\ge1}$,
so $N_Q$ is simply connected.
Under the induced identification, we have that
$T_1N_Q=\mathbb RQ_{\ge1}$, with the Lie bracket 
$[x,y]=xy-yx$.
Thus the Lie algebra of $N_Q$ is $\nq$.
By Proposition~\ref{prop:structure} and the correspondence
of the lower central series for connected Lie groups,
$N_Q$ is nilpotent of step $\ell$ when $L(Q)=\ell\ge1$.

The subgroup $\Gamma_Q$ is discrete.
Order the paths by nondecreasing length. Successively multiply
on the left by $1-\lfloor a_p\rfloor p$ to move each coefficient
$a_p$ into $[0,1)$; only coefficients of longer paths change.
Thus the compact set
$\{1+\sum_{p\in\pp}a_pp\mid0\le a_p\le1\}$ maps onto
$\Gamma_Q\backslash N_Q$.
Since $N_Q$ is unimodular, this compact quotient has finite
invariant measure, so $\Gamma_Q$ is a cocompact lattice.
\end{proof}

The existence of a cocompact lattice in $N_Q$ also follows
from Mal\cprime cev's lattice criterion \cite{Malcev}, \cite{Raghunathan}.

\begin{proposition}\label{prop:solvablegroups}
The Lie group $S_Q$ is simply connected and completely solvable,
with Lie algebra $\sq$.
If $L(Q)=\ell\ge1$, then $S_Q$ is
solvable of step $\lfloor\log_2\ell\rfloor+2$.
If $L(Q)=0$, then $S_Q$ is abelian.
If $Q$ has an arrow, then $S_Q$ has no lattice.
\end{proposition}

\begin{proof}
$S_Q=\{\sum_v t_ve_v+x\mid t_v>0,\ x\in\mathbb RQ_{\ge1}\}$,
which is diffeomorphic to $\mathbb R^{n(Q)+m(Q)}$. Hence  $S_Q$ is simply connected.
Under the induced identification $T_1S_Q=\mathbb RQ$,
the Lie bracket is $[x,y]=xy-yx$.
Thus the Lie algebra of $S_Q$ is $\sq$.

Order the paths by nondecreasing length, placing
the trivial paths first.
For $h\in\hq$ and $x\in\nq$, 
it is easy to check that the matrix of $\ad(h+x)$ is triangular, with
diagonal entries $0$ and $\alpha_p(h)$.
All these entries are real, so $S_Q$ is completely solvable.
By Proposition~\ref{prop:structure} and the correspondence
of the derived series for connected Lie groups,
$S_Q$ is solvable of step $\lfloor\log_2\ell\rfloor+2$
when $L(Q)=\ell\ge1$.

Suppose that $Q$ has an arrow.
Choose a topological ordering $\tau$ of $Q$ and let
$h_0=-\sum_{v\in Q_0}\tau(v)e_v$.
For every non-trivial path $p$, we have
$\alpha_p(h_0)=\tau(t(p))-\tau(s(p))>0$.
Since $\pp$ is nonempty, it follows that
$\operatorname{tr}(\ad h_0)=\sum_{p\in\pp}\alpha_p(h_0)>0$.
Thus $S_Q$ is not unimodular.
A Lie group admitting a lattice is unimodular
\cite{Raghunathan}, so $S_Q$ has no lattice.
\end{proof}

\subsection{Graphs and matchings}
We recall the graph and matching notions used to formulate
the symplecticity criteria.
Standard references  are \cite{Bollobas}, \cite{Diestel}, and \cite{LP}.

All undirected graphs and multigraphs considered below are finite and have no loops.
A multigraph may have parallel edges, which are counted separately.
A graph means a simple graph, with no parallel edges.
For an undirected multigraph $F$, write $V(F)$, $E(F)$,
and $c(F)$ for its vertex set, edge set, and number of
connected components, respectively. Two vertices are adjacent, or neighbors, if an edge joins them.
For $S\subseteq V(F)$, let $N(S)$ be the set of vertices
adjacent to at least one vertex of $S$.
The degree of a vertex in an undirected multigraph
is the number of edges incident to it,
with parallel edges counted separately. A vertex with degree zero is said to be isolated.
An undirected path has no repeated vertices.
Two vertices belong to the same connected component
if and only if they are joined by an undirected path.

\begin{definition}\label{def:forest}
A cycle of length $r\ge2$ consists of distinct vertices
$v_1,\cdots,v_r$ and distinct edges $e_1,\cdots,e_r$,
where $e_i$ joins $v_i$ to $v_{i+1}$ and $v_{r+1}=v_1$.
In particular, two parallel edges form a cycle of length two.
A forest is an undirected multigraph without cycles,
and a tree is a connected forest.
The cycle rank of an undirected multigraph $F$ is
$b_1(F)=|E(F)|-|V(F)|+c(F)$;
it is zero if and only if $F$ is a forest.
\end{definition}

In a forest, any two distinct vertices in the same component
are joined by a unique undirected path,
and every component containing an edge has a vertex
of degree one.

Let $F$ be a multigraph. For $S\subseteq V(F)$, the induced multigraph $F[S]$
has vertex set $S$ and contains all edges of $F$
joining vertices in $S$.

The underlying undirected multigraph $U(Q)$ is obtained
by forgetting the directions of the arrows of $Q$,
while retaining each arrow as a separate edge. Two vertices of $Q$ are adjacent if they are adjacent in $U(Q)$.
The connected components of $Q$ are those of $U(Q)$. A connected component is non-trivial if it is not an isolated vertex.
Write $c(Q)=c(U(Q))$, including isolated vertices.
The multigraph $U(Q)[I(Q)]$ has the internal vertices
of $Q$ as its vertex set and contains all edges
between them.

\begin{definition}\label{def:graphsubsets}
A vertex cover of a graph is a set containing at least one
vertex of every edge. A set is independent if no two of its
vertices are adjacent.
The complete graph $K_r$ has $r$ vertices and one edge joining
each pair of distinct vertices.
The complete bipartite graph $K_{r,s}$ has two parts of sizes
$r,s$, with all edges between the parts and none within either part.
The central vertex of $K_{1,\ell}$ is the vertex in the part
of size one; for $\ell=1$, either vertex may be chosen.
\end{definition}

\begin{definition}\label{def:matching}
A matching in a graph is a set of edges having no common
vertices.
It is perfect if it covers every vertex.
The matching number is the maximum size of a matching.

In a bipartite graph with vertex set $L\sqcup R$,
where every edge joins a vertex of $L$ to a vertex of $R$,
a matching saturates $L$ if it covers every vertex of $L$.
\end{definition}

Hall's theorem \cite{Hall} states that a matching saturating $L$
exists if and only if $|N(S)|\ge|S|$ for every $S\subseteq L$.

For a skew-symmetric matrix $A=(a_{ij})$ of order $2N$,
its Pfaffian is
\begin{equation*}
 \Pf(A)=\sum_{M\in\PM(\{1,\cdots,2N\})}\eps(M)
 \prod_{\{i,j\}\in M,\ i<j}a_{ij},
\end{equation*}
where $\PM(\{1,\cdots,2N\})$ is the set of perfect matchings
of the complete graph on $\{1,\cdots,2N\}$.
For $M=\{\{i_1,j_1\},\cdots,\{i_N,j_N\}\}$, ordered so that
$i_k<j_k$ and $i_1<\cdots<i_N$, let $\eps(M)$ be the sign
of the permutation $(i_1,j_1,\cdots,i_N,j_N)$. The identity $\det A=\Pf(A)^2$ (see \cite{Knuth}) shows that
$A$ is nondegenerate if and only if $\Pf(A)\ne0$.

\begin{definition}\label{def:supportgraph}
The support graph of a skew-symmetric matrix $(a_{ij})$
has one vertex for each index and an edge $\{i,j\}$
exactly when $a_{ij}\ne0$.
\end{definition}

If this graph has a unique perfect matching, the Pfaffian
has exactly one nonzero term, so the matrix is nonsingular.

\begin{definition}\label{def:graphlie}

For a finite simple graph $G$, choose an orientation of each edge. Define the graph Lie algebra $\mathfrak n_G$ to be the Lie algebra with basis
$\{x_v\mid v\in V(G)\}\cup\{z_e\mid e\in E(G)\}$ and bracket
$[x_u,x_v]=z_e$ for each edge $e$ oriented from $u$ to $v$,
with all remaining brackets zero except those determined
by skew-symmetry. Different choices of orientation give isomorphic Lie algebras.
\end{definition}

\subsection{Symplectic forms, reduction and cohomology}\label{subsec:symplectictools}\label{subsec:lefschetzprelim}
We recall the relation between symplectic forms on Lie
algebras and left-invariant symplectic forms on Lie groups,
together with reduction by isotropic ideals.
For the algebraic and geometric background used in this
subsection, see  \cite{BC}, \cite{Huybrechts}, and \cite{Weibel}.

For a finite-dimensional vector space $V$, let $V^*$ be
its dual and $\Lambda^rV^*$ the space of alternating
$r$-linear forms.
The exterior product is denoted by $\wedge$.

\begin{definition}
\label{def:radical}
The radical of an alternating two-form $\omega$ on $V$ is
$\rad\omega=\{x\in V\mid\omega(x,V)=0\}$.
The form is nondegenerate if its radical is zero.
A subspace $W\subseteq V$ is isotropic with respect to
$\omega$ if $\omega(W,W)=0$.
\end{definition}

We use the Chevalley--Eilenberg differential \cite{CE}
with coefficients in $\kk$, where $x\cdot a=0$
for all $x\in\mathfrak g$ and $a\in\kk$.
Our sign convention is
\begin{align*}
 (d\lambda)(x,y)&=-\lambda([x,y]),\\
 (d\omega)(x,y,z)&=-\omega([x,y],z)+\omega([x,z],y)
 -\omega([y,z],x).
\end{align*}
The Jacobi identity gives $d^2=0$.

\begin{definition}
\label{def:symplectic}
An alternating two-form $\omega$ on a Lie algebra
$\mathfrak g$ is closed if $d\omega=0$ and exact if
$\omega=d\lambda$ for some $\lambda\in\mathfrak g^*$. A symplectic structure on $\mathfrak g$ is a closed
nondegenerate alternating two-form.

Write
$Z^2(\mathfrak g)
=\ker(d:\Lambda^2\mathfrak g^*\to\Lambda^3\mathfrak g^*)$
for the space of closed two-forms and
$B^2(\mathfrak g)=d(\mathfrak g^*)$
for the space of exact two-forms.
The second cohomology space is
$H^2(\mathfrak g)=Z^2(\mathfrak g)/B^2(\mathfrak g)$.

\end{definition}

By the definition of $d\omega$,  closedness is equivalent to
\begin{equation}\label{eq:cyclic}
 \omega([x,y],z)+\omega([y,z],x)+\omega([z,x],y)=0.
\end{equation}
A Lie algebra admitting a symplectic structure has even
dimension.
In dimension $2N$, an alternating two-form $\omega$ is
nondegenerate if and only if $\omega^N\ne0$.

Pouseele and Tirao \cite{PT} give a necessary and sufficient
condition for a graph Lie algebra to admit a symplectic structure.

\begin{theorem}\label{thm:PT}
The graph Lie algebra $\mathfrak n_G$ is symplectic if and only if
$|V(G)|+|E(G)|$ is even and $|E(H)|\le|V(H)|$
for every connected component $H$ of $G$.
\end{theorem}

\begin{lemma}\label{lem:centerorth}
For every closed two-form $\omega$ on a Lie algebra $\mathfrak g$, one has $\omega(Z(\mathfrak g),[\mathfrak g,\mathfrak g])=0$. If $\omega$ is nondegenerate, the map $Z(\mathfrak g)\to(\mathfrak g/[\mathfrak g,\mathfrak g])^*$ given by $z\mapsto(\bar x\mapsto\omega(z,x))$ is well-defined and injective.
\end{lemma}
\begin{proof}
For central $z$, \eqref{eq:cyclic} gives $\omega([x,y],z)=0$. An element in the kernel of the displayed map is therefore orthogonal to the whole algebra and must be zero.
\end{proof}

The next lemma identifies the radical vectors that cannot be removed
by adding a form pulled back from a quotient.

For a linear map $\pi:V\to H$ and a form $\beta$ on $H$, its pullback is $(\pi^*\beta)(x,y)=\beta(\pi x,\pi y)$.
\begin{lemma}\label{lem:completion}
Let $\pi:V\to H$ be surjective, with kernel $N$,
and let $\eta\in\Lambda^2V^*$. There exists $\beta\in\Lambda^2H^*$ such that $\eta+\pi^*\beta$ is nondegenerate if and only if $\dim V$ is even and $N\cap\rad\eta=0$.
\end{lemma}
\begin{proof}
If $\eta+\pi^*\beta$ is nondegenerate, then
$N\cap\rad\eta=\{0\}$, since every vector in this intersection
belongs to $\rad(\eta+\pi^*\beta)$. Conversely, let $R_\eta=\rad\eta$ and write $\dim V=2N_0$, $\rank\eta=2r$. Then $\pi|_{R_\eta}$ is injective and $\dim R_\eta=2(N_0-r)$. Choose $\beta_0$ on $H$ whose restriction to $\pi(R_\eta)$ is nondegenerate. If $V=E\oplus R_\eta$, then $\eta|_E$ is nondegenerate.

In $(\eta+t\pi^*\beta_0)^{N_0}$, all coefficients of degrees below $N_0-r$ vanish, and the coefficient of $t^{N_0-r}$ is $\binom{N_0}{r}\eta^r\wedge(\pi^*\beta_0)^{N_0-r}$. 
This coefficient is nonzero, since $\eta|_E$ and
$(\pi^*\beta_0)|_{R_\eta}$ are both nondegenerate.
Hence $(\eta+t\pi^*\beta_0)^{N_0}$ is not identically zero
as a polynomial in $t$. So, we can choose
$t_0\in\kk$ such that $(\eta+t_0\pi^*\beta_0)^{N_0}\ne0$.
Thus $\eta+\pi^*(t_0\beta_0)$ is nondegenerate, and we take
$\beta=t_0\beta_0$.
\end{proof}

\begin{definition}\label{def:symplecticgroups}
A symplectic form on a real manifold is a closed nondegenerate
two-form. A form $\Omega$ on a Lie group $G$ is left-invariant
if $L_g^*\Omega=\Omega$ for every left translation $L_g$.
A symplectic Lie group is a Lie group with a left-invariant
symplectic form \cite{Medina}.

On a complex manifold, a holomorphic symplectic form is a
closed holomorphic two-form that is nondegenerate over $\mathbb C$.
A complex manifold is complex parallelizable if its holomorphic
tangent bundle is trivial, or equivalently if it has a global
frame of holomorphic one-forms.

A compact nilmanifold is a quotient $\Gamma\backslash N$,
where $N$ is a simply connected nilpotent real Lie group
and $\Gamma$ is a cocompact lattice.

On a smooth complex algebraic variety, an algebraic symplectic
form is a closed nondegenerate regular two-form.
\end{definition}

Every alternating two-form on a real Lie algebra extends
uniquely by left translation to a left-invariant two-form
on any Lie group with that Lie algebra.
This correspondence preserves closedness and nondegeneracy,
and hence gives a bijection between symplectic forms on the
Lie algebra and left-invariant symplectic forms on the Lie group.
The same statement holds for complex alternating two-forms
on a complex Lie algebra and left-invariant holomorphic
two-forms on a complex Lie group.

\begin{definition}\label{def:lagrangian}
For a subspace $W$ of a symplectic vector space $(V,\omega)$,
let $W^{\perp_\omega}=\{x\in V\mid\omega(x,W)=0\}$.
We call $W$ coisotropic if $W^{\perp_\omega}\subseteq W$ and
Lagrangian if $W=W^{\perp_\omega}$.
Thus an isotropic subspace is Lagrangian exactly when
$2\dim W=\dim V$.
An isotropic or Lagrangian Lie subalgebra, ideal, or subgroup
has the corresponding property with respect to the restricted
form. For a subgroup with a left-invariant form, this condition
can be checked on its Lie algebra at the identity.

Following \cite{BC}, the symplectic rank is
$\sigma(\mathfrak g,\omega)=\max_I\dim I$, where $I$ runs
through the isotropic ideals.
\end{definition}

In a simply connected completely solvable real Lie group,
Lie subalgebras correspond to connected closed subgroups,
and ideals correspond to connected closed normal subgroups
\cite{Raghunathan}. This applies in particular to $N_Q$ and $S_Q$.

We use the following standard facts about symplectic reduction
\cite{BC}.

\begin{lemma}\label{lem:isotropicreduction}
Let $(\mathfrak g,\omega)$ be a symplectic Lie algebra and let
$I$ be an isotropic ideal. Then $I$ is abelian,
$I^{\perp_\omega}$ is a subalgebra, and $\omega$ induces a
symplectic form $\bar\omega$ on $I^{\perp_\omega}/I$.
If $\bar L$ is a Lagrangian subalgebra of this quotient,
its inverse image in $I^{\perp_\omega}$ is a Lagrangian
subalgebra of $\mathfrak g$.
\end{lemma}
\begin{proof}
For $x,y\in I$ and $z\in\mathfrak g$, the last two terms of
\eqref{eq:cyclic} vanish by the ideal condition and isotropy.
Hence $\omega([x,y],z)=0$, so $[x,y]=0$.
For $x,y\in I^{\perp_\omega}$ and $z\in I$, the same identity
shows that $\omega([x,y],z)=0$.
Thus $I^{\perp_\omega}$ is a subalgebra.
The radical of the restricted form is
$I^{\perp_\omega}\cap(I^{\perp_\omega})^{\perp_\omega}=I$,
which proves the quotient assertion.
The inverse image of $\bar L$ is isotropic and has dimension
$\dim I+\tfrac12(\dim\mathfrak g-2\dim I)=\tfrac12\dim\mathfrak g$.
\end{proof}

We call $(I^{\perp_\omega}/I,\bar\omega)$ the symplectic
reduction by $I$. 
\begin{definition}\label{def:momentum}
For an action of a Lie group $G$ preserving a symplectic form
$\Omega$ on $M$, write
$x^\#_p=\left.\frac{d}{dt}\right|_{t=0}\exp(tx)\cdot p$
for the fundamental vector field of $x\in\mathfrak g$.
A momentum map is a map $\mu:M\to\mathfrak g^*$ satisfying
$d\langle\mu,x\rangle=\iota_{x^\#}\Omega$ for every $x$.
We call the action Hamiltonian if it admits such a map.
\end{definition}

\begin{definition}\label{def:rational}
For $\kk=\mathbb R$, a subspace of $\nq$ or $\sq$ is rational
if it has a basis with rational coordinates in the path basis.
A two-form is rational if its matrix in the path basis
has rational entries.
A connected subgroup of $N_Q$ is rational if its Lie algebra
is rational.
\end{definition}

\begin{remark}\label{rem:canonicalconnection}
A symplectic form $\omega$ on a real Lie algebra defines
a torsion-free flat left-invariant connection by
$\omega(\nabla_xy,z)=-\omega(y,[x,z])$.
This connection is complete exactly when  the corresponding Lie group is
unimodular \cite{Medina}; it need not preserve $\omega$.
For symplectic $N_Q$ this connection is complete.
For symplectic $S_Q$ it is incomplete whenever $Q$ has an arrow,
by Proposition~\ref{prop:solvablegroups}.

\end{remark}

For the cohomological calculations in Section~\ref{sec:nilcohomlefschetz},
we record the Lefschetz notation and two rules for graded strings.

\begin{definition}\label{def:lefschetz}
For a closed two-form $\omega$ on a manifold $M$, the Lefschetz
operator $L_\omega$ is multiplication by $[\omega]$ on cohomology.
Thus $L_\omega^u:H^k(M)\to H^{k+2u}(M)$ is defined by
$\eta\mapsto\eta[\omega]^u$.
A symplectic manifold of dimension $2d$ has the hard Lefschetz
property if $L_\omega^{d-k}$ is an isomorphism for every
$0\le k\le d$. We use the same notation on Lie algebra cohomology.
\end{definition}

\begin{definition}\label{def:gradedstrings}
For a nilpotent operator $L$ of degree two, a graded string
$(a,\ell)$ is a basis $v,Lv,\cdots,L^{\ell-1}v$ with
$\deg v=a$ and $L^\ell v=0$.
Every finite-dimensional graded space with such an operator
is a direct sum of strings, counted with multiplicity.
\end{definition}

\begin{lemma}\label{lem:stringrules}
For the operator $L\otimes1+1\otimes L$, the tensor product
$(a,\ell)\otimes(b,h)$ is
$\bigoplus_{0\le j<\min(\ell,h)}(a+b+2j,\ell+h-1-2j)$.
On the exterior algebra of a space of dimension $e$, multiplication
by an alternating form of rank $2r$ gives strings
$(v+w,r-v+1)$ with multiplicity
$\bigl(\binom{2r}{v}-\binom{2r}{v-2}\bigr)\binom{e-2r}{w}$,
where $0\le v\le r$ and $0\le w\le e-2r$,
with $\binom ab=0$ for $b<0$ or $b>a$.
\end{lemma}
\begin{proof}
The first formula follows from the tensor product rule
for $\mathfrak{sl}_2$ representations \cite{FH},
and the second from linear Lefschetz decomposition
\cite{Huybrechts}.
\end{proof}

\section{Symplectic structures on \texorpdfstring{$\nq$}{nQ}}\label{sec:nil}

We classify the quivers for which $\nq$ admits a symplectic form.
Closedness first gives an obstruction at path length four.
When $L(Q)=3$, closedness and nondegeneracy impose degree
and parity conditions on the active families.
Combining these conditions with the symplectic criterion
for the residual graph Lie algebra gives the classification
in Theorem~\ref{thm:nilclassification}.

\begin{lemma}\label{lem:maxsupportnil}
Let $p=a_1\cdots a_r$ be maximal, with successive vertices $v_0,\cdots,v_r$, and write $I(p)=\{v_1,\cdots,v_{r-1}\}$ for the set of  its internal vertices. For a closed two-form $\omega$ on $\nq$, $\omega(p,q)=0$ for every path $q$ of length at least two. If $r\ge2$ and $a$ is an arrow, then $\omega(p,a)\ne0$ implies $I(p)\subseteq\{s(a),t(a)\}$.
\end{lemma}
\begin{proof}
The first assertion follows from Proposition~\ref{prop:structure} and Lemma~\ref{lem:centerorth}. If $v_i\notin\{s(a),t(a)\}$, split $p=uw$ at $v_i$. The products $wa$ and $au$ are zero because $v_r$ is a sink and $v_0$ is a source. The products $aw$ and $ua$ are zero because $t(a)\ne v_i$ and $s(a)\ne v_i$. Thus $[w,a]=[a,u]=0$, while $[u,w]=p$. Equation~\eqref{eq:cyclic} gives $\omega(p,a)=0$.
\end{proof}

\begin{definition}\label{def:hallgraph}
Let $\mathcal H_Q$ be the bipartite graph with left
vertex set $\mathcal M(Q)$ and right vertex set $Q_1$.
A maximal path $p$ of length at least two is joined to $a$
when $I(p)\subseteq\{s(a),t(a)\}$.
A maximal arrow is joined to every arrow other than itself.
\end{definition}

\begin{proposition}\label{prop:nilHall}
If $\nq$ is symplectic, then $m(Q)$ is even and $\mathcal H_Q$ has a matching saturating $\mathcal M(Q)$.
\end{proposition}

\begin{proof}
Let $\omega$ be a symplectic form on $\nq$.
Then $m(Q)=\dim\nq$ is even.
Moreover, $Z(\nq)=\Span\mathcal M(Q)$, and the classes
of the arrows form a basis of $\nq/[\nq,\nq]$.
By Lemma~\ref{lem:centerorth}, $\omega$ induces an injective
map from $Z(\nq)$ to $(\nq/[\nq,\nq])^*$.
Hence the matrix $A=(\omega(p,a))$, with rows indexed by
$\mathcal M(Q)$ and columns indexed by $Q_1$, has full row rank.
By Lemma~\ref{lem:maxsupportnil} and the alternation of $\omega$,
$A_{p,a}\ne0$ only if $p$ and $a$ are adjacent in $\mathcal H_Q$.

Let $r=|\mathcal M(Q)|$.
Since $A$ has full row rank, some $r$ columns give an
$r\times r$ submatrix $A'$ with $\det A'\ne0$.
Write its row indices as $p_1,\cdots,p_r$ and its column
indices as $a_1,\cdots,a_r$.
The determinant expansion gives a permutation $\sigma\in S_r$
such that $\prod_{i=1}^r A'_{i,\sigma(i)}\ne0$.
Thus $p_i$ and $a_{\sigma(i)}$ are adjacent in $\mathcal H_Q$
for every $i$.
Since the arrows $a_{\sigma(i)}$ are distinct, these edges
form a matching saturating $\mathcal M(Q)$.
\end{proof}

The support restriction also produces radical vectors common to
all closed two-forms, giving the following bound on path length.

\begin{theorem}\label{thm:long}
Every maximal path of length at least four belongs to the radical of every closed two-form on $\nq$. Consequently, $L(Q)\ge4$ implies that $\nq$ is not symplectic.
\end{theorem}
\begin{proof}
Such a maximal path has at least three distinct internal vertices. Lemma~\ref{lem:maxsupportnil} makes it orthogonal to every arrow and every longer path; hence it is a nonzero radical vector. 
\end{proof}

\begin{definition}\label{def:activearrows}
Let $L(Q)\le3$.
An arrow $b$ from $u$ to $v$ is active if it is the middle arrow
of a path of length three, equivalently if $d^-(u)>0$ and $d^+(v)>0$.
Let $\Acal$ be the set of active arrows, and let
$\mathcal A_{uv}$ be the set of active arrows from $u$ to $v$.
Write $m_{uv}=|\mathcal A_{uv}|$.
A pair $(u,v)$ is active if $m_{uv}>0$;
we then call $\mathcal A_{uv}$ an active family.
Let $\Rcal=Q_1\setminus\Acal$ be the set of residual arrows.
The residual graph $\GR$ has vertex set $\Rcal$,
where distinct arrows $a,b$ are adjacent exactly when
$ab$ or $ba$ is a path of length two.
\end{definition}

If $L(Q)=0$, then $\nq=0$.
If $L(Q)=1$, then $\nq$ is abelian and is symplectic
exactly when $|Q_1|$ is even.
If $L(Q)=2$, the existence of a symplectic structure on $\nq$
is characterized in \cite{Mizoguchi}.

For $L(Q)\le2$, there are no active arrows, so $\Rcal=Q_1$
and $\nq\cong\mathfrak n_{\GR}$.
Thus Theorem~\ref{thm:PT} also gives the following criterion.

\begin{theorem}\label{thm:niltwo}
Let $L(Q)\le2$.
Then $\nq$ is symplectic if and only if
$|Q_1|+\sum_{v\in Q_0}d^-(v)d^+(v)$ is even and
$d^-(v)d^+(v)\le d^-(v)+d^+(v)$ for every $v\in Q_0$.
At an internal vertex, the latter condition is equivalent
to $(d^-(v),d^+(v))$ being $(1,r)$, $(r,1)$ or $(2,2)$,
where $r\ge1$.
\end{theorem}

Now assume $L(Q)=3$.
The next two lemmas give necessary degree and parity conditions
on the active families for $\nq$ to admit a symplectic form.

\begin{lemma}\label{lem:threeboundary}
If $\nq$ is symplectic, then every active pair $(u,v)$ satisfies $d^-(u)=d^+(v)=1$.
\end{lemma}
\begin{proof}
Let $\omega$ be a symplectic form on $\nq$. Every path of length three is maximal. Let $p$ be such a path with internal vertices $u,v$, and let $a$ be an arrow. By Lemma~\ref{lem:maxsupportnil}, $\omega(p,a)\ne0$ only if $\{u,v\}=\{s(a),t(a)\}$. By acyclicity, any such arrow $a$ must be one of the $m_{uv}$ arrows from $u$ to $v$. There are $d^-(u)m_{uv}d^+(v)$ such paths of length three. 

Let $Z_{uv}$ be the span of these paths of length three.
Then $Z_{uv}\subseteq Z(\nq)$.
By Lemma~\ref{lem:centerorth} and the preceding restriction,
$\omega$ induces an injective map
$Z_{uv}\to(\Span\mathcal A_{uv})^*$.
Hence $d^-(u)m_{uv}d^+(v)\le m_{uv}$. Both positive degrees are one.

\end{proof}

Under this degree condition, let $a$ be the unique arrow entering $u$ and $c$ the unique arrow leaving $v$. Write $\mathcal A_{uv}=\{b_1,\cdots,b_m\}$, and set $X_i=ab_i$, $Y_i=b_ic$, and $Z_i=ab_ic$. The brackets used below are
 $[a,b_i]=X_i$, $[b_i,c]=Y_i$, 
 $[a,Y_i]=Z_i$, and $[X_i,c]=Z_i$.

\begin{lemma}\label{lem:threeeven}
If $\nq$ is symplectic, every active multiplicity $m_{uv}$ is even.
\end{lemma}

\begin{proof}
Let $\omega$ be symplectic and set $D_{ij}=\omega(Z_i,b_j)$
for one active family $\mathcal A_{uv}=\{b_1,\cdots,b_m\}$.
By Lemma~\ref{lem:maxsupportnil}, each $Z_i$ is orthogonal
to all arrows outside $\mathcal A_{uv}$.
Since each $Z_i$ is central, Lemma~\ref{lem:centerorth}
also makes it orthogonal to $[\nq,\nq]$.
Nondegeneracy therefore makes the induced map
$\Span\{Z_1,\cdots,Z_m\}\to\Span\{b_1,\cdots,b_m\}^*$
injective. Both spaces have dimension $m$, so $D$ is nonsingular.
Closedness on $(a,Y_i,b_j)$ and $(X_j,c,b_i)$ gives
$\omega(X_j,Y_i)=D_{ij}=-D_{ji}$.
Thus $D$ is skew-symmetric and has even order.
\end{proof}

We separate the paths containing active arrows from the residual
paths; the quotient will be a graph Lie algebra.

Assume now that $d^-(u)=d^+(v)=1$ for every active pair $(u,v)$.
For an active arrow $b$, let $a_b,c_b$ be its unique
predecessor and successor, and let
$K=\Span\{b,a_bb,bc_b,a_bbc_b\mid b\in\Acal\}$.
A path cannot contain two active arrows.
Indeed, extending before the first and after the second
would give a path of length at least four.
These paths are therefore distinct, and $\dim K=4|\Acal|$. 

\begin{lemma}\label{lem:corequotient}
The space $K$ is an ideal of $\nq$, and $\nq/K\cong\mathfrak n_{\GR}$. Every path of length two made of residual arrows is maximal.
If $H_R$ is the span of the residual arrows and these
paths of length two, then $\nq=K\oplus H_R$ and
$\dim\nq=4|\Acal|+|V(\GR)|+|E(\GR)|$.

\end{lemma}
\begin{proof}

Let $x,y\in\Rcal$ with $xy$ a path.
If $axy$ were a path for some arrow $a$, then $x$ would
be active; similarly, a path $xyc$ would make $y$ active.
Thus $xy$ is maximal.
Every path of length three has an active middle arrow.
The degree conditions imply that every path of length two
containing an active arrow $b$ is $a_bb$ or $bc_b$,
and every path of length three is $a_bbc_b$.

The bracket of any path in the path basis of $K$
with an arrow belongs to $K$.
Since the arrows generate $\nq$ as a Lie algebra,
the Jacobi identity gives $[\nq,K]\subseteq K$.
The path basis of $K$, together with the residual arrows
and the paths of length two formed from them,
forms a basis of $\nq$.
Hence $\nq=K\oplus H_R$, and the classes of the latter paths
form a basis of $\nq/K$.
Their brackets are precisely those defining
$\mathfrak n_{\GR}$, so $\nq/K\cong\mathfrak n_{\GR}$.
Counting these paths gives the dimension formula.

\end{proof}

\begin{lemma}\label{lem:residualnecessary}
If $\nq$ is symplectic, then $|E(H)|\le|V(H)|$ for every connected component $H$ of $\GR$, and $|V(\GR)|+|E(\GR)|$ is even.
\end{lemma}
\begin{proof}
Let $\omega$ be a symplectic form on $\nq$.
Let $Z_3=\mathfrak n_Q^3$, and let $Z_R$ be the span of the residual paths of length two. They are central. 
Consider the matrix $(\omega(p,a))$, whose rows are indexed
first by the paths of length three and then by the paths of length two
consisting of residual arrows, and whose columns are
indexed first by the active arrows and then by the residual
arrows.
It has the block form
$\left(\begin{smallmatrix}D&0\\E&F\end{smallmatrix}\right)$,
where $D$ is the direct sum of the nonsingular matrices from Lemma~\ref{lem:threeeven}. The whole matrix has full row rank by Lemma~\ref{lem:centerorth}. Eliminating $E$ shows that $F$ has full row rank.

If the edge $\{x,y\}$ corresponding to the path $xy$
lies in one connected component of $\GR$
and the residual arrow $z$ lies in another,
then $[y,z]=[z,x]=0$.
Equation~\eqref{eq:cyclic} gives $\omega(xy,z)=0$. 
Thus $F$ is block diagonal with respect to the connected
components of $\GR$.
For each component $H$, the corresponding block has
$|E(H)|$ rows and $|V(H)|$ columns.
Since $F$ has full row rank, each block has full row rank,
so $|E(H)|\le|V(H)|$.
Finally, $\dim\nq$ is even, and the dimension formula
in Lemma~\ref{lem:corequotient}
implies that $|V(\GR)|+|E(\GR)|$ is even.
\end{proof}

For sufficiency, we construct a closed two-form on $\nq$
whose restriction to $K$ is nondegenerate.
Adding the pullback of a symplectic form on $\nq/K$
will then give a symplectic form on $\nq$.

For each active pair $(u,v)$, write
$\mathcal A_{uv}=\{b_1,\cdots,b_m\}$ and assume that $m$ is even.
Choose a nonsingular skew-symmetric matrix
$J=(J_{ij})\in M_m(\kk)$.
Using the paths $X_i,Y_i,Z_i$ defined above,
define an alternating bilinear form $\omega_{uv}$ on $\nq$ by
$\omega_{uv}(Z_i,b_j)=\omega_{uv}(X_j,Y_i)=J_{ij}$
for $1\le i,j\le m$.
Extend by alternation and set all other values on pairs of paths to zero. Let $\omega_K=\sum_{(u,v)\text{ active}}\omega_{uv}$.

\begin{lemma}\label{lem:coreformclosed}
The form $\omega_K$ is closed, its restriction to $K$ is nondegenerate, and it vanishes whenever one entry lies in $H_R$.
\end{lemma}
\begin{proof}

For each active pair $(u,v)$, the chosen matrix $J$
is nonsingular, so the restriction of $\omega_{uv}$
to $\Span\{b_i,X_i,Y_i,Z_i\mid 1\le i\le m\}$ is nondegenerate.
These subspaces form a direct sum equal to $K$, and
$\omega_K$ vanishes between distinct summands.
Hence $\omega_K|_{K\times K}$ is nondegenerate.
By definition, $\omega_K(H_R,\nq)=0$.

For closedness, let $p,q,r$ be non-trivial paths. A nonzero term $\omega_{uv}([p,q],r)$ must have its bracket equal, up to sign, to some $X_j,Y_i$, or $Z_i$; a bracket cannot be an arrow. The paths of length two have unique splittings, and $Z_i$ has the two splittings $a\cdot Y_i$ and $X_i\cdot c$. Thus, up to permutation, it suffices to verify \eqref{eq:cyclic} for
$(a,b_j,Y_i)$ and $(X_j,c,b_i)$.
The corresponding sums are $J_{ij}-J_{ij}=0$ and
$J_{ji}+J_{ij}=0$, respectively.
The same verification applies when different active families
share their predecessor or successor.
Thus every $\omega_{uv}$ is closed.

\end{proof}

\begin{theorem}\label{thm:nilthree}
Let $L(Q)=3$. The algebra $\nq$ is symplectic if and only if every active pair satisfies $d^-(u)=d^+(v)=1$ and has even multiplicity $m_{uv}$, every connected component $H$ of $\GR$ satisfies $|E(H)|\le|V(H)|$, and $|V(\GR)|+|E(\GR)|$ is even.
\end{theorem}
\begin{proof}
Assume that $\nq$ is symplectic. The assertions follow from Lemma~\ref{lem:threeboundary}, Lemma~\ref{lem:threeeven}, and Lemma~\ref{lem:residualnecessary}.  Conversely, take $\omega_K$ as in
Lemma~\ref{lem:coreformclosed}. The graph conditions in the statement and Theorem~\ref{thm:PT} give a symplectic form $\bar\omega_R$ on $\nq/K\cong\mathfrak n_{\GR}$. If $\pi:\nq\to\nq/K$ is the quotient map, then $\omega_K+\pi^*\bar\omega_R$ is closed. The form $\omega_K+\pi^*\bar\omega_R$ is nondegenerate
on both $K$ and $H_R$ and vanishes on $K\times H_R$.
Since $\nq=K\oplus H_R$, it is nondegenerate on $\nq$
and hence symplectic.
\end{proof}
Combining the length obstruction with the criteria for
$L(Q)\le2$ and $L(Q)=3$ gives the complete classification.

\begin{theorem}\label{thm:nilclassification}
The Lie algebra $\nq$ is symplectic if and only if $m(Q)$ is even
and one of the following holds:
\begin{enumerate}
\item $L(Q)\le2$ and $d^-(v)d^+(v)\le d^-(v)+d^+(v)$
for every vertex $v$;
\item $L(Q)=3$, every active pair $(u,v)$ satisfies
$d^-(u)=d^+(v)=1$, every $m_{uv}$ is even, and every
connected component $H$ of $\GR$ satisfies $|E(H)|\le|V(H)|$.
\end{enumerate}
\end{theorem}
\begin{proof}
Theorem~\ref{thm:long} excludes $L(Q)\ge4$, and
Theorem~\ref{thm:niltwo} gives the assertion for $L(Q)\le2$.
For $L(Q)=3$, apply Theorem~\ref{thm:nilthree}.
By the dimension formula in Lemma~\ref{lem:corequotient},
the condition that $|V(\GR)|+|E(\GR)|$ is even is
equivalent to $m(Q)$ being even.
\end{proof}
Moreover, Theorem~\ref{thm:long} and Theorem~\ref{thm:nilthree}
imply that if $Q$ has no parallel arrows and $L(Q)\ge3$,
then $\nq$ is not symplectic.

\section{Geometry of positive path Lie groups and compact quotients}\label{sec:groups}

The classification of $\nq$ gives criteria for the existence
of left-invariant symplectic forms on $N_Q$ and symplectic
forms on its compact quotients.
We then study Lagrangian normal subgroups for a fixed
left-invariant symplectic form.
The resulting criterion distinguishes the quivers for which
such a subgroup exists for every form and, for rational forms,
determines when it induces a Lagrangian torus fibration.
In this section the Lie groups and manifolds are real unless
stated otherwise, and the corresponding Lie algebras are taken
over $\mathbb R$.

\subsection{\texorpdfstring{Symplectic forms on $N_Q$ and its compact quotients}{Symplectic forms on NQ and its compact quotients}}\label{subsec:nilforms}

By Proposition~\ref{prop:pathgroups}, $N_Q$ admits the
cocompact lattice $\Gamma_Q=1+\mathbb ZQ_{\ge1}$.
For every cocompact lattice $\Gamma\subset N_Q$, the quotient
$\Gamma\backslash N_Q$ is a compact nilmanifold.
The following theorem characterizes the existence of
left-invariant symplectic forms on $N_Q$ and symplectic forms
on its compact quotients, with no invariance assumption
on the latter.

\begin{theorem}\label{thm:geometricexistence}
The Lie group $N_Q$ admits a left-invariant symplectic form exactly
when the conditions of Theorem~\ref{thm:nilclassification} hold.
For every cocompact lattice $\Gamma$ in $N_Q$, the same
conditions are equivalent to the existence of a symplectic
form on $\Gamma\backslash N_Q$.
\end{theorem}
\begin{proof}
A left-invariant form is determined by its value at the identity,
and closedness is the Chevalley--Eilenberg equation.
Thus the first assertion follows from Section~\ref{sec:nil}.
A left-invariant symplectic form descends to every left quotient
by a lattice.

Conversely, suppose that $M=\Gamma\backslash N_Q$ has dimension
$2d$ and admits a symplectic form $\eta$.
Nomizu's theorem \cite{Nomizu} identifies
$H^*(\nq;\mathbb R)$ with $H^*(M;\mathbb R)$ as graded algebras.
Choose a left-invariant closed two-form $\omega$ representing
$[\eta]$. Since $[\omega]^d=[\eta]^d\ne0$, we have
$\omega^d\ne0$. Left invariance makes $\omega$ nondegenerate
everywhere. Hence $\nq$ is symplectic.
\end{proof}

The argument using invariant representatives is classical;
see \cite{BG} and \cite{Nomizu}.
In particular, if $L(Q)\ge4$, no compact quotient of $N_Q$
admits a symplectic form.
Moreover, for $M=\Gamma\backslash N_Q$ with $\dim M=2d>0$,
a class $u\in H^2(M;\mathbb R)$
has a symplectic representative exactly when $u^d\ne0$.

For the complex quotient below, compactness forces the
coefficients of every holomorphic two-form in a global
holomorphic coframe to be constant.

\begin{proposition}\label{prop:complexquotients}
Let $N_Q(\mathbb C)=1+\mathbb CQ_{\ge1}$ and
$\Gamma_Q(\mathbb C)=1+\mathbb Z[i]Q_{\ge1}$.
Then the quotient
$X_Q=\Gamma_Q(\mathbb C)\backslash N_Q(\mathbb C)$
is a compact complex parallelizable manifold. Moreover,
the following conditions are equivalent:
\begin{enumerate}
\item $Q$ satisfies the conditions of
Theorem~\ref{thm:nilclassification}.
\item $N_Q(\mathbb C)$ admits a left-invariant
holomorphic symplectic form.
\item $X_Q$ admits a holomorphic symplectic form.
\end{enumerate}
\end{proposition}

\begin{proof}
The subgroup $\Gamma_Q(\mathbb C)$ is discrete.
Applying the coefficient reduction in
Proposition~\ref{prop:pathgroups} to the real and imaginary
parts of the path coefficients proves compactness of $X_Q$.
Choose a basis of left-invariant holomorphic one-forms
on $N_Q(\mathbb C)$.
Their descent gives a global holomorphic frame on $X_Q$,
so $X_Q$ is complex parallelizable.
Denote this frame by $\theta_1,\cdots,\theta_{m(Q)}$.

Every holomorphic two-form on $X_Q$ is a linear combination
of the forms $\theta_i\wedge\theta_j$ with holomorphic
coefficients, which are constant by compactness and
connectedness.
Consequently, pullback identifies holomorphic symplectic
forms on $X_Q$ with left-invariant holomorphic symplectic
forms on $N_Q(\mathbb C)$.
The equivalence follows from Theorem~\ref{thm:nilclassification}
applied to $\nq$ over $\mathbb C$.

\end{proof}

\begin{remark}\label{rem:complexkahler}
If $L(Q)\ge2$, the lattice $\Gamma_Q(\mathbb C)$ is nonabelian,
as $1+a$ and $1+b$ do not commute whenever $ab\ne0$.
Thus $X_Q$ is not a complex torus and is not K\"ahler
by Wang's theorem \cite{Wang}.
\end{remark}

\subsection{Lagrangian normal subgroups and torus fibrations}\label{subsec:nilnormal}\label{subsec:niltori}
We study connected closed Lagrangian normal subgroups
of $N_Q$ through the corresponding Lagrangian ideals
of $\nq$.
We give an existence criterion for a given left-invariant
symplectic form and determine when it holds for every
such form.
For rational forms, the criterion also determines when
$\Gamma_Q\backslash N_Q$ admits a Lagrangian torus
fibration induced by a rational normal subgroup.
From this subsection to the end of Section~\ref{sec:groups},
assume that $\nq$ admits a symplectic form, so $L(Q)\le3$.

For $L(Q)=3$, use $K$ and $H_R$ from
Lemma~\ref{lem:corequotient}.
For $L(Q)\le2$, set $K=0$, $H_R=\nq$ and
$\Acal=\varnothing$; use the graph $\GR$ defined in
Section~\ref{sec:nil} for $L(Q)\le2$.
In either case $H_R$ is a subalgebra nilpotent of step at most two and
$\nq=K\rtimes H_R$. Write $\omega_K(J)$ for the form constructed before
Lemma~\ref{lem:coreformclosed} from a family
$J=(J_{uv})_{(u,v)\text{ active}}$ of nonsingular
skew-symmetric matrices
$J_{uv}\in M_{m_{uv}}(\mathbb R)$.

To study Lagrangian ideals for a given symplectic form,
we first use an automorphism to make $K$ and $H_R$ orthogonal.

\begin{lemma}\label{lem:corenormal}
The ideal $K$ is abelian. For every symplectic form $\omega$,
its restriction to $K$ is nondegenerate.
There is a Lie algebra automorphism preserving $K$ which
transforms $\omega$ into $\omega_K(J)+\omega_R$,
where $\omega_K(J)$ is the form of
Lemma~\ref{lem:coreformclosed} and $\omega_R$ is symplectic
on $H_R$. The transformed form vanishes on $K\times H_R$.
 If $\omega$ has rational coefficients in the path basis,
the automorphism can be chosen to have a rational matrix
in that basis.
\end{lemma}

\begin{proof}
If $K=0$, the assertions are immediate. Suppose $L(Q)=3$.
A product of two paths in $K$ would contain two active arrows,
which is impossible by $L(Q)=3$. Hence $K$ is abelian.

Recall from Section~\ref{sec:nil} that, for each active pair
$(u,v)$, $a$ is the unique arrow entering $u$,
$c$ is the unique arrow leaving $v$,
$\mathcal A_{uv}=\{b_1,\cdots,b_m\}$,
and $X_i=ab_i$, $Y_i=b_ic$, $Z_i=ab_ic$.
Choose $b_i\in\mathcal A_{uv}$ and let $b_j$ range over
all active arrows, with $X_j,Y_j,Z_j$ defined
using the active family containing $b_j$.
Closedness on $(a,b_i,X_j)$ and $(b_i,c,Y_j)$,
together with $L(Q)=3$, gives
$\omega(X_i,X_j)=\omega(Y_i,Y_j)=0$.
Closedness on $(a,b_i,Y_j)$ gives
$\omega(X_i,Y_j)=\omega([a,Y_j],b_i)$.
The bracket $[a,Y_j]$ is either $Z_j$ or zero.
By Lemma~\ref{lem:maxsupportnil},
$\omega(X_i,Y_j)=0$ whenever $b_i,b_j$
belong to different active families.
Within one family, Lemma~\ref{lem:threeeven} gives
$\omega(Z_i,b_j)=\omega(X_j,Y_i)=J_{ij}$,
where $J$ is nonsingular and skew-symmetric.
The $Z_i$ are orthogonal to all paths of length at least two
and to all arrows outside their active family.
For $k\in K$, the equations
$\omega(k,Z_i)=0$, $\omega(k,X_i)=\omega(k,Y_i)=0$,
and $\omega(k,b_i)=0$ imply successively that its
coefficients at $b_i$, at $X_i,Y_i$, and at $Z_i$ vanish,
since each $J$ is nonsingular.
Here the indices range over all active families.
Hence $\omega(k,K)=0$ implies $k=0$, so $\omega|_K$
is nondegenerate.

For an abelian ideal, closedness shows that its orthogonal
complement is a subalgebra.
Here $K^{\perp_\omega}$ is a complement to $K$ and projects
isomorphically onto $H_R$.
It is the graph of a linear map $f:H_R\to K$.
Since $K$ is abelian, the subalgebra condition gives
$f([x,y])=[x,f(y)]-[y,f(x)]$ for $x,y\in H_R$.
Therefore $F(k+x)=k+x+f(x)$ defines an automorphism,
and $F^*\omega$ vanishes on $K\times H_R$.

Write $\omega_0=\omega_K(J)|_K$ and define $S$ on $K$ by
$\omega(x,y)=\omega_0(x,Sy)$ for $x,y\in K$.
The two forms agree on pairs of paths whose lengths sum
to four, and both vanish when the sum exceeds four.
Hence $S=\id_K+T$, where
$T$ strictly increases path length, and $T^3=0$.
Both forms are invariant under the $H_R$-action, so
$T([h,k])=[h,T(k)]$ for $h\in H_R$ and $k\in K$.
Also, $\omega_0(Sx,y)=\omega_0(x,Sy)$ for $x,y\in K$.
Let $A=\id_K-\tfrac12T+\tfrac38T^2$.
Since $T^3=0$, we have $A^2S=\id_K$, and hence
$\omega(Ax,Ay)=\omega_0(x,y)$ for $x,y\in K$.
The operator $A$ commutes with the $H_R$-action,
so it extends to an automorphism $\widetilde A$
by the identity on $H_R$. Both $F$ and $\widetilde A$ preserve $K$, so their
composition also preserves $K$.
The restriction $\omega_R=(F^*\omega)|_{H_R}$ is closed
and nondegenerate.
Extending $\omega_R$ by zero whenever one argument lies
in $K$, we obtain
$(F\circ\widetilde A)^*\omega=\omega_K(J)+\omega_R$.

If $\omega$ has rational coefficients in the path basis,
then the linear systems and polynomial changes of basis above
are defined over $\mathbb Q$, so the automorphism is rational.
\end{proof}

If the matrix of $\omega$ in $K\oplus H_R$ is
$\left(\begin{smallmatrix}A&B\\-B^T&C\end{smallmatrix}\right)$,
then the induced form $\omega_R$ on $H_R$ has matrix
$C+B^TA^{-1}B$.
Thus the inverse form graph in Definition~\ref{def:inversegraph}
can be computed without finding the normalizing automorphism.
\begin{lemma}\label{lem:residualstars}
Every connected component of $\GR$ is an isolated vertex,
a complete bipartite graph $K_{1,\ell}$ with $\ell\ge1$, or a complete bipartite graph $K_{2,2}$.
Each non-trivial component contains at most one arrow that
is a predecessor or successor of an active family.
If such an arrow occurs, the component is a $K_{1,\ell}$ and the arrow can
be chosen as its central vertex.
\end{lemma}
\begin{proof}
Delete the active arrows from $Q$.
The remaining quiver has length at most two, since the middle
arrow of a  path of length three would be active.
No arrow of this quiver joins two internal vertices.
Thus every non-trivial component of its concatenation graph
is $K_{p,q}$, arising from the incoming and outgoing arrows
at one internal vertex.
Theorem~\ref{thm:PT} and Lemma~\ref{lem:residualnecessary}
give $pq\le p+q$; hence $(p-1)(q-1)\le1$.
This gives precisely the stated shapes.

A predecessor of an active family starts at a source and
is the only arrow entering the starting vertex of that family.
If this arrow has a neighbor in $\GR$, its connected
component is a $K_{1,\ell}$ with the arrow as its central vertex.
 No other arrow of this $K_{1,\ell}$ can be a
predecessor, by the degree condition, or a successor, since
that would yield a path containing two active arrows.
The reversed argument applies to a successor.
\end{proof}

\begin{definition}\label{def:activedata}
We call the predecessors and successors of active families
boundary arrows.
The set of predecessors and the set of successors are disjoint,
since an arrow belonging to both would give a path
of length at least four.
The active-family graph $G_A(Q)$ is the bipartite graph with these two sets
as its parts, with one edge joining the predecessor
and successor of each active family.

In each connected component of $\GR$ isomorphic to
$K_{1,\ell}$, choose a central vertex, choosing the boundary
arrow if the component contains one.
This is possible by Lemma~\ref{lem:residualstars}.
Together with all isolated vertices of $\GR$,
denote the chosen vertices by $a_1,\cdots,a_b$.
We call $a_1,\cdots,a_b$ the chosen residual arrows.
Every boundary arrow occurs among them.
If $a_i$ belongs to a component $K_{1,\ell}$,
let $\ell_i=\ell$; if $a_i$ is isolated, let $\ell_i=0$.
Let $c_\square$ be the number of components isomorphic
to $K_{2,2}$.

For each active family with predecessor $a_i$ and
successor $a_j$, let $r_{ij}=r_{ji}$ be the number of its active arrows.
Let $r_{ii}=0$, and let $r_{ij}=0$ if no active family
has $a_i,a_j$ as its predecessor and successor
in either order.
Sums over $\{i,j\}$ are taken over unordered pairs
of distinct indices, so each active family is counted once.
\end{definition}

Counting the non-trivial paths gives
$m(Q)=b+2\sum_{i=1}^b\ell_i
+4\sum_{\{i,j\}}r_{ij}+8c_\square$.
Since $m(Q)=\dim\nq$ is even, $b$ is even.

Let $\Omega$ be a left-invariant symplectic form on $N_Q$,
and write $\omega$ for its value at the identity.
By the correspondence recalled in Section~\ref{sec:prelim},
Lagrangian ideals of $(\nq,\omega)$ correspond to connected
closed Lagrangian normal subgroups of $(N_Q,\Omega)$.
We use this to decide when $(N_Q,\Omega)$ has a Lagrangian
normal subgroup and when such a subgroup induces a torus
fibration on a compact quotient.

The graph $G_A(Q)$ records the active families.
The next graph records the constraint imposed by $\omega$ on
the choice of a vertex cover.

\begin{definition}\label{def:inversegraph}
The projection $\nq=K\oplus H_R\to H_R$ restricts
to an isomorphism from $K^{\perp_\omega}$ onto $H_R$.
Use this isomorphism to transfer
$\omega|_{K^{\perp_\omega}}$ to a form $\omega_R$ on $H_R$.
For a boundary arrow $a$, let $a^*$ be the element
of the basis of $H_R^*$ dual to the path basis
that corresponds to $a$.
The inverse form $\omega_R^{-1}$ is the alternating form
on $H_R^*$ whose matrix in the dual basis is the inverse
of the matrix of $\omega_R$.
The inverse form graph $F_\omega$ has the boundary arrows
as vertices, with $a,c$ adjacent exactly when
$\omega_R^{-1}(a^*,c^*)\ne0$.
\end{definition}

\begin{theorem}\label{thm:nilidealcriterion}
The symplectic Lie group $(N_Q,\Omega)$ has a connected closed
Lagrangian normal subgroup if and only if $G_A(Q)$ has a
vertex cover which is independent in $F_\omega$.
If $\omega$ is rational and $G_A(Q)$ has a vertex cover
which is independent in $F_\omega$, then $(N_Q,\Omega)$
has a rational connected closed Lagrangian normal subgroup.
\end{theorem}
\begin{proof}
It suffices to prove the corresponding assertion for
Lagrangian ideals of $(\nq,\omega)$.
By Lemma~\ref{lem:corenormal}, we may assume that
$\omega=\omega_K+\omega_R$.
Let $L$ be a Lagrangian ideal and let $x\in L$.
On the symplectic space $K$, the map $\rho(x)=\ad(x)|_K$
is skew-adjoint, by closedness and $[K,K]=0$.
Its image is contained in $L\cap K$ and is isotropic with respect to $\omega_K$.
Thus $\rho(x)^2=0$, since
$\omega_K(\rho(x)y,\rho(x)z)=-\omega_K(y,\rho(x)^2z)$.

Choose an active family $\mathcal A_{uv}$ with predecessor $a$
and successor $c$, and let $b_i\in\mathcal A_{uv}$.
For $x\in L$, denote by $x_a,x_c$ its coefficients at
$a,c$ in the path basis.
Then $\rho(x)^2b_i=[x,[x,b_i]]=-2x_ax_cab_ic$.
Since $\rho(x)^2=0$ and $ab_ic\ne0$, we obtain
$x_ax_c=0$ for every $x\in L$.

The product of two linear maps cannot vanish identically
on a vector space over an infinite field unless one of them
is zero. 
Consequently, for each edge $\{a,c\}$ of $G_A(Q)$,
either $x_a=0$ for every $x\in L$ or
$x_c=0$ for every $x\in L$.
Let $T$ be the set of all boundary arrows $a$ such that
$x_a=0$ for every $x\in L$, and let
$H_T=\bigcap_{a\in T}\ker a^*\subseteq H_R$.
Since $L\subseteq K\oplus H_T$ and $L=L^{\perp_\omega}$,
taking orthogonal complements gives
$H_T^{\perp_{\omega_R}}\subseteq L\subseteq K\oplus H_T$.
Hence $H_T^{\perp_{\omega_R}}\subseteq H_T$.
Under the isomorphism defined by $\omega_R$, this says
$\omega_R^{-1}(a^*,c^*)=0$ for all $a,c\in T$.
Thus $T$ is independent in $F_\omega$.

Conversely, let $T$ have these two properties.
Then $H_T$ is coisotropic and contains
$C_R=[H_R,H_R]$. Since $C_R\subseteq Z(H_R)$, Lemma~\ref{lem:centerorth} implies that $C_R$
is isotropic with respect to $\omega_R$.
Also $C_R\subseteq H_T$, so
$C_R+H_T^{\perp_{\omega_R}}$ is isotropic with respect to $\omega_R$.
Let $\pi:H_T\to H_T/H_T^{\perp_{\omega_R}}$
be the natural projection.
Extend $\pi(C_R)$ to a Lagrangian subspace $\overline P$
with respect to the symplectic form induced by $\omega_R$.
Then $P=\pi^{-1}(\overline P)$ is a Lagrangian subspace
of $(H_R,\omega_R)$ satisfying $C_R\subseteq P\subseteq H_T$.
Since $[H_R,P]\subseteq [H_R,H_R]=C_R\subseteq P$,
the subspace $P$ is an ideal of $H_R$.

For each active family, if its predecessor is not in $T$,
choose $U_{ac}=\Span\{ab_i,ab_ic\}_i$.
If its successor is not in $T$, choose
$U_{ac}=\Span\{b_ic,ab_ic\}_i$.
The cover condition prevents both cases from occurring.
If both vertices belong to $T$, choose the first space.
Let $U$ be the sum of the subspaces $U_{ac}$ over all
active families.
Then $U$ is Lagrangian in $(K,\omega_K)$ and satisfies
$[H_R,U]\subseteq U$ and $[P,K]\subseteq U$.
Since $K$ is abelian and $P$ is an ideal of $H_R$,
the subspace $U\oplus P$ is an ideal of $\nq$.
It is isotropic with respect to $\omega_K+\omega_R$
and has dimension $\frac12\dim\nq$, so it is Lagrangian.
If $\Phi$ is the automorphism from Lemma~\ref{lem:corenormal}
satisfying $\Phi^*\omega=\omega_K+\omega_R$,
then $\Phi(U\oplus P)$ is a Lagrangian ideal
with respect to the original form $\omega$.
When $\omega$ is rational, all the spaces and extensions
can be chosen over $\mathbb Q$.
\end{proof}

\begin{remark}\label{rem:twosat}
Let $t_a$ be a Boolean variable indicating
whether the boundary arrow $a$ belongs to $T$.
The graph condition requires $t_a\lor t_c$ for every
edge $\{a,c\}$ of $G_A(Q)$ and $\neg t_a\lor\neg t_c$
for every edge $\{a,c\}$ of $F_\omega$.
Thus, after computing $\omega_R^{-1}$,  the existence
of a vertex cover of $G_A(Q)$ which is independent
in $F_\omega$ reduces to a $2$-SAT problem;
see \cite{APT} for details.
The existence of a Lagrangian normal subgroup  can depend on
the choice of the left-invariant symplectic form.
\end{remark}

\begin{theorem}\label{thm:allnilideals}
The Lie group $N_Q$ admits a connected closed Lagrangian normal
subgroup for each left-invariant symplectic form if and
only if $G_A(Q)$ has matching number at most one.
\end{theorem}
\begin{proof}

A bipartite graph of matching number at most one has all
its edges incident to one vertex, or has no edges.
Every singleton vertex set is independent in $F_\omega$, so
Theorem~\ref{thm:nilidealcriterion} applies.

Suppose that $G_A(Q)$ contains two disjoint edges.
For each component of $\GR$ isomorphic to $K_{1,\ell}$,
let $a$ be the chosen central vertex and let
$u_1,\cdots,u_\ell$ be the other vertices.
For each $j$, let $z_j$ be the nonzero product
of the arrows $a$ and $u_j$, taken in the order
in which they concatenate.
On $\Span\{u_1,\cdots,u_\ell,z_1,\cdots,z_\ell\}$,
choose the nondegenerate alternating form
$\sum_{j=1}^{\ell}u_j^*\wedge z_j^*$.
For each $K_{2,2}$ component, Theorem~\ref{thm:PT}
gives a symplectic form on its graph Lie algebra.
The form can be chosen over $\mathbb Q$, since rational
closed forms are dense in the real space of closed forms
and nondegeneracy is an open condition.

Let $V=\Span\{a_1,\cdots,a_b\}$, where
$a_1,\cdots,a_b$ are the previously chosen residual arrows.
On $V$, choose the alternating form whose matrix $B$
in the  basis $(a_1,\cdots,a_b)$ has all entries
above the diagonal equal to one.
Since $b$ is even, a direct calculation shows that $B$
is invertible and $(B^{-1})_{ij}=(-1)^{i+j}$ for $i<j$.
Hence every entry of $B^{-1}$ outside the diagonal is nonzero.
For each active pair $(u,v)$, choose a nonsingular
skew-symmetric matrix $J_{uv}$ with rational entries,
and define $\omega_K$ by the construction preceding
Lemma~\ref{lem:coreformclosed}.
Define $\omega$ on $\nq$ as the direct sum of
$\omega_K|_K$ and the forms chosen above.
These choices give a rational symplectic form
$\omega=\omega_K+\omega_R$ for which $F_\omega$ is complete.
Every vertex cover of $G_A(Q)$ contains at least two vertices,
so none is independent in $F_\omega$.
By Theorem~\ref{thm:nilidealcriterion}, $N_Q$ has no
connected closed Lagrangian normal subgroup with respect
to the left-invariant symplectic form determined by $\omega$.

\end{proof}

Even when the condition for all forms fails, a suitable rational
symplectic form always admits a rational Lagrangian ideal.
Rationality then allows the corresponding subgroup to descend
to a torus fibration of the standard compact quotient.

\begin{proposition}\label{prop:rationalnilideal}
Since $\nq$ admits a symplectic form, there is a rational
symplectic form on $\nq$ with a rational Lagrangian ideal.
Equivalently, $N_Q$ admits a left-invariant symplectic form
with a rational connected closed Lagrangian normal subgroup.
\end{proposition}
\begin{proof}
Let $V=\Span\{a_1,\cdots,a_b\}$.
Let $T$ be the part of the bipartition
of $G_A(Q)$ with the smaller number of vertices,
choosing either part if the two parts have the same size.
Since every boundary arrow belongs to
$\{a_1,\cdots,a_b\}$, we have
$|T|\le\frac12|V(G_A(Q))|\le b/2$.
There is a rational nondegenerate alternating form $\Pi$ on
$V^*$ such that $\Pi(a^*,c^*)=0$ for all $a,c\in T$.
Indeed, extend these coordinate vectors to a rational basis
$f_1,\cdots,f_{2s}$ of $V^*$, with
$f_1,\cdots,f_{|T|}$ corresponding to the arrows in $T$.
Define $\Pi$ by $\Pi(f_i,f_{s+j})=\delta_{ij}$ for
$1\le i,j\le s$, with all other values zero except those
determined by skew-symmetry.
Use $\Pi^{-1}$ on $V$ and keep the other rational forms
chosen in the proof of Theorem~\ref{thm:allnilideals}.
Then $T$ is a vertex cover of $G_A(Q)$ independent in
$F_\omega$, so the rational assertion of
Theorem~\ref{thm:nilidealcriterion} applies.
\end{proof}

\begin{theorem}\label{thm:torusfibration}

Let $\dim\nq=2d$, let $\omega$ be a rational symplectic
form on $\nq$, and let $M_Q=\Gamma_Q\backslash N_Q$.
$M_Q$, equipped with the symplectic form induced by $\omega$,
has a Lagrangian torus fibration induced by a rational
normal subgroup if and only if $G_A(Q)$ has a vertex cover
which is independent in $F_\omega$.
Since $\nq$ admits a symplectic form, it admits a rational
symplectic form $\omega$ such that $M_Q$, equipped with the
symplectic form induced by $\omega$, has a Lagrangian torus
fibration induced by a rational normal subgroup of $N_Q$.
\end{theorem}

\begin{proof}
Suppose that $G_A(Q)$ has a vertex cover which is independent
in $F_\omega$.
By Theorem~\ref{thm:nilidealcriterion}, there exists a rational
Lagrangian ideal $L$ of $(\nq,\omega)$.
By Lemma~\ref{lem:isotropicreduction}, $L$ is abelian.
The connected subgroup $H=\exp L$ is closed and normal.
Mal\cprime cev's rational subgroup theorem \cite{Malcev}, \cite{Raghunathan}
shows that $\Gamma_Q\cap H$ is a lattice in $H$ and that
the image $\bar\Gamma_Q$ of $\Gamma_Q$ in $N_Q/H$ is a lattice.
The map
$\Gamma_Q\backslash N_Q\to\bar\Gamma_Q\backslash(N_Q/H)$
has $d$-dimensional torus fibers.
Since $\omega|_{L\times L}=0$, left invariance implies that
the induced symplectic form vanishes on each fiber.
Hence the fibers are Lagrangian.

Conversely, the Lie algebra of a rational normal subgroup
inducing such a fibration is a rational Lagrangian ideal.
Theorem~\ref{thm:nilidealcriterion} therefore implies that
$G_A(Q)$ has a vertex cover which is independent in $F_\omega$.

The existence assertion follows from
Proposition~\ref{prop:rationalnilideal} and the preceding construction.
\end{proof}

\begin{remark}
This criterion concerns Lagrangian torus fibrations induced
by rational normal subgroups. Its failure does not exclude
other Lagrangian fibrations.
\end{remark}

The following example shows that the existence of a Lagrangian
normal subgroup can change with the symplectic form on a fixed
$N_Q$, and gives the smallest dimension in the family of
positive path Lie algebras where failure occurs.
The two forms constructed below are rational, so the distinction
also holds for Lagrangian torus fibrations induced by rational
normal subgroups on the same compact quotient.

\begin{example}\label{ex:twenty}
Among the Lie algebras $\nq$ associated with finite acyclic
quivers, the minimum dimension admitting a symplectic form
with no Lagrangian ideal is $20$.

This minimum is attained by the following quiver $Q$.
\begin{center}
\begin{tikzcd}[column sep=large, row sep=large]
& u_1
  \arrow[r, shift left=1, "b_{11}"]
  \arrow[r, shift right=1, "b_{12}"']
& v_1 \arrow[dr, "c_1"] & \\
s \arrow[ur, "a_1"] \arrow[dr, "a_2"']
& & & t \\
& u_2
  \arrow[r, shift left=1, "b_{21}"]
  \arrow[r, shift right=1, "b_{22}"']
& v_2 \arrow[ur, "c_2"'] &
\end{tikzcd}
\end{center}
Here $\dim\nq=20$.
For $i=1,2$, let
$K_i=\Span\{b_{ij},a_ib_{ij},b_{ij}c_i,a_ib_{ij}c_i
\mid j=1,2\}$.
Let $\omega_K$ be the form constructed before
Lemma~\ref{lem:coreformclosed}, taking
$J=\begin{psmallmatrix}0&1\\-1&0\end{psmallmatrix}$
for each of the two active families.
 Let $\theta_1,\cdots,\theta_4$ be the basis
of $H_R^*$ dual to $(a_1,c_1,a_2,c_2)$.
Define
$\omega_0=\omega_K+\theta_1\wedge\theta_2+
\theta_3\wedge\theta_4$ and
$\omega_1=\omega_K+\sum_{i<j}\theta_i\wedge\theta_j$.
Both are symplectic forms with integer coefficients
in the path basis.

For $i=0,1$, let $\Omega_i$ be the left-invariant
extension of $\omega_i$ to $N_Q$.
Let $\overline\Omega_i$ be the form on
$M_Q=\Gamma_Q\backslash N_Q$ induced by $\Omega_i$.
Then $(N_Q,\Omega_0)$ has a rational connected closed
Lagrangian normal subgroup, whereas $(N_Q,\Omega_1)$
has no connected closed Lagrangian normal subgroup.
Consequently, $(M_Q,\overline\Omega_0)$ admits a
Lagrangian torus fibration induced by a rational normal
subgroup, with ten-dimensional fibers, whereas
$(M_Q,\overline\Omega_1)$ admits no such fibration.

\end{example}

\begin{proof}
Let $\omega$ be a symplectic form on $\nq$ such that
$(\nq,\omega)$ has no Lagrangian ideal.
By Theorem~\ref{thm:allnilideals}, $G_A(Q)$ contains
two disjoint edges.
These correspond to two active families with four
distinct boundary arrows.
Each active family has positive even multiplicity,
so the two families contain at least four active arrows.
Each active arrow $b$, with predecessor $a$ and
successor $c$, contributes four distinct  paths
$b,ab,bc,abc$.
Thus there are at least sixteen distinct paths, in addition
to the four boundary arrows, and $\dim\nq\ge20$.

For the displayed quiver, a direct calculation gives $\dim\nq=20$ and shows that
$\omega_0,\omega_1$ are symplectic forms with integer
coefficients.
Moreover, $G_A(Q)$ and $F_{\omega_0}$ both consist of
the two edges $\{a_1,c_1\}$ and $\{a_2,c_2\}$, whereas
$F_{\omega_1}$ is complete.
The set $\{a_1,a_2\}$ is a vertex cover of $G_A(Q)$
and is independent in $F_{\omega_0}$.
Every vertex cover of $G_A(Q)$ has at least two vertices,
whereas every independent set in $F_{\omega_1}$ has at most one.
Theorem~\ref{thm:nilidealcriterion} therefore gives
the assertions about Lagrangian normal subgroups.
Both left-invariant forms descend to
$M_Q=\Gamma_Q\backslash N_Q$, and the assertions
about Lagrangian torus fibrations follow from
Theorem~\ref{thm:torusfibration}.

\end{proof}

\section{Symplectic structures on
\texorpdfstring{$\sq$}{sQ}}\label{sec:full}

We characterize symplecticity of $\sq$ by the structure of the
underlying forest and its internal vertices.
The description of closed two-forms reduces the problem to full
row rank of a matrix determined by path concatenation.
We then prove the structural criterion by combining a
forest-matching obstruction with a terminal-path reduction.
The reduction also gives a recursive criterion for symplecticity and a
construction of symplectic forms; its geometric realization
is given in Section~\ref{sec:solvablegroups}.

We extend $\beta\in\Lambda^2\hq^*$ and $\lambda\in\nq^*$
to $\sq=\hq\oplus\nq$ by setting
$\beta(h_1+x_1,h_2+x_2)=\beta(h_1,h_2)$ and
$\lambda(h+x)=\lambda(x)$ for $h,h_1,h_2\in\hq$
and $x,x_1,x_2\in\nq$.
We use the same symbols for these extensions.

Every path weight $\alpha_p$, for $p\in\pp$, is nonzero, and $\alpha_p+\alpha_q\ne0$ for $p,q\in\pp$. Indeed, for a topological ordering $\tau$, evaluation at $h_0=-\sum_{v\in Q_0}\tau(v)e_v$ gives $\alpha_p(h_0)=\tau(t(p))-\tau(s(p))$, a positive integer. Also, $\alpha_{pq}=\alpha_p+\alpha_q$ whenever $pq$ is a path.

\begin{theorem}\label{thm:closedforms}
Every closed two-form $\omega$ on $\sq$ has a unique expression $\omega=\beta+d\lambda$, with $\beta\in\Lambda^2\hq^*$ and $\lambda\in\nq^*$. Thus
\begin{equation*}
 Z^2(\sq)=\Lambda^2\hq^*\oplus d(\nq^*),\qquad
 B^2(\sq)=d(\nq^*),\qquad H^2(\sq)\cong\Lambda^2\hq^*.
\end{equation*}
Their dimensions are respectively $\binom{n(Q)}2+m(Q)$, $m(Q)$ and $\binom{n(Q)}2$.
\end{theorem}
\begin{proof}
For a closed $\omega$, the equality
$(d\omega)(h_1,h_2,p)=0$ gives  $\alpha_p(h_1)\omega(p,h_2)=\alpha_p(h_2)\omega(p,h_1)$, where $h_1,h_2\in\hq$ and $p\in\pp$. Hence $h\mapsto\omega(p,h)$ is a unique scalar multiple of $\alpha_p$. Define $\lambda(p)$ by $\omega(p,h)=\lambda(p)\alpha_p(h)$. This defines a linear map on the path basis, and $(d\lambda)(p,h)=\lambda(p)\alpha_p(h)$.

The closed form $\zeta=\omega-d\lambda$ vanishes on $\nq\times\hq$. Applying \eqref{eq:cyclic} to $(h,p,q)$, for $p,q\in\pp$, gives $(\alpha_p+\alpha_q)(h)\zeta(p,q)=0$. Therefore $\zeta(p,q)=0$, and $\zeta$ is the extension of its restriction $\beta$ to $\hq$.

Conversely, every such $\beta$ is closed because $[\sq,\sq]\subseteq\nq$, and $d\lambda$ is closed by $d^2=0$. If $\beta+d\lambda=0$, evaluation on $(p,h)$ gives $\lambda(p)\alpha_p(h)=0$, so $\lambda=0$ and then $\beta=0$. This also proves that $d:\nq^*\to\Lambda^2\sq^*$ is injective. Linear maps on $\hq$ have zero differential, since $[\sq,\sq]=\nq$. All exact forms therefore come from $\nq^*$. The quotient and dimension assertions follow.
\end{proof}

The decomposition is also a consequence of the weight zero
cohomology argument in \cite[Theorem~2.1]{CG}.
Its extension to all degrees appears in
Proposition~\ref{prop:fullcohomology}.
The remaining problem is to determine when $\beta+d\lambda$
can be nondegenerate.

\begin{proposition}\label{prop:noexact}
If $n(Q)>0$, no exact two-form on $\sq$ is symplectic. If $\nq\ne0$, the restriction of every closed two-form on $\sq$ to $\nq$ is degenerate.
\end{proposition}
\begin{proof}
The identity $1=\sum_v e_v$ is a nonzero central vector of $\sq$ and lies in the radical of every exact form. On $\nq$, the restriction of $\beta+d\lambda$ is exact, and every maximal non-trivial path is a nonzero central vector in its radical.
\end{proof}

For $\lambda\in\nq^*$, set $B_\lambda(x,y)=\lambda([x,y])$, and denote its radical by $\sq^\lambda$. $B_\lambda$ and $d\lambda$ have the same radical. Applying Lemma~\ref{lem:completion} to the natural projection
$\sq\to\hq$ and using Theorem~\ref{thm:closedforms},
we see that $\sq$ is symplectic if and only if
$\dim\sq$ is even and $\sq^\lambda\cap\nq=0$
for some $\lambda\in\nq^*$.

Write $n=n(Q)$ and $m=m(Q)$. Choose orderings of $Q_0$ and $\pp$, and associate an independent indeterminate $X_p$ to every $p\in\pp$. Let $\kk[\mathbf X]$ be their polynomial ring and $\FQ$ its fraction field. Define the  matrix $\KQ$ with rows indexed by $\pp$ and columns indexed by $Q_0\sqcup\pp$ by
\begin{align*}
 \mathcal K_Q(p,v)&=(\delta_{v,t(p)}-\delta_{v,s(p)})X_p,\\
 \mathcal K_Q(p,q)&=\begin{cases}
 X_{pq},&pq\text{ is a path},\\
 -X_{qp},&qp\text{ is a path},\\
 0,&\text{otherwise}.
 \end{cases}
\end{align*}
These are the values $\lambda_{\mathbf X}([p,z])$ of the linear map $\lambda_{\mathbf X}(p)=X_p$. Let $\gamma(Q)=\rank_{\FQ}\KQ$. This is the generic rank, or the largest size of a square minor that is not the zero polynomial. Full row rank means $\gamma(Q)=m$.

\begin{theorem}\label{thm:fullrank}
The algebra $\sq$ is symplectic if and only if $n+m$ is even and $\gamma(Q)=m$.
\end{theorem}

\begin{proof}
For each $p\in\pp$, choose $t_p\in\kk$ and define
$\lambda_{\mathbf t}:\nq\to\kk$ by
$\lambda_{\mathbf t}(p)=t_p$ on the path basis.
Let $\mathcal K_Q(\mathbf t)$ be obtained from
$\mathcal K_Q(\mathbf X)$ by substituting $X_p=t_p$
for every $p\in\pp$.
The transpose $\mathcal K_Q(\mathbf t)^T$ represents the map
$y\mapsto B_{\lambda_{\mathbf t}}(y,\cdot)$ from $\nq$ to $\sq^*$,
whose kernel is $\nq\cap\sq^{\lambda_{\mathbf t}}$.
The preceding radical criterion therefore gives symplecticity
exactly when $n+m$ is even and some such matrix has row rank $m$.

If some $\mathcal K_Q(\mathbf t)$ has row rank $m$,
choose a nonzero $m\times m$ minor of this matrix.
The corresponding minor of $\mathcal K_Q(\mathbf X)$
is then a nonzero polynomial.
Conversely, a nonzero polynomial
minor has a nonzero value for some choice of the $t_p$.
Thus a matrix of row rank $m$ exists if and only if
$\gamma(Q)=m$.

\end{proof}

For $\lambda\in\nq^*$, write $\mathcal K_Q(\lambda)$
for the matrix obtained from $\mathcal K_Q(\mathbf X)$
by substituting $X_p=\lambda(p)$ for every $p\in\pp$.

Repeated path variables allow different determinant terms to
cancel. We first record the signed matching criterion, before
passing to the structural characterization of full row rank.

Define a bipartite graph $\mathcal B_Q$ with left vertex set $\pp$ and right vertex set $Q_0\sqcup\pp$. An edge $p$--$z$ is present exactly when $\mathcal K_Q(p,z)\ne0$,
and it has that signed monomial as label.

Choose a set $C$ of $m$ right vertices, with the inherited order. A matching $M$ from all left vertices onto $C$ determines a permutation. Multiply the sign of this permutation by the signs of the corresponding edge labels to obtain $\sigma(M)\in\{1,-1\}$, and write $\mu(M)$ for the product of the corresponding unsigned labels. 
For a monomial \(\mu\), set \(c_{C,\mu}=\sum_{\mu(M)=\mu}\sigma(M)\), where the sum runs over matchings \(M\) saturating the left vertex set \(\mathcal P(Q)\) and having \(C\) as their set of matched right vertices.
Then $\det\mathcal K_Q[\pp,C]=\sum_\mu c_{C,\mu}\mu$. In particular, Theorem~\ref{thm:fullrank} gives the following statement.
\begin{proposition}\label{prop:signedmatching}
The algebra $\sq$ is symplectic if and only if $n+m$ is even and $c_{C,\mu}\ne0$ for some $C,\mu$.

If $\sq$ is symplectic, then $\mathcal B_Q$ has a matching
saturating $\pp$; equivalently, $|N(S)|\ge |S|$ for every
$S\subseteq\pp$, where $N(S)$ is the set of right vertices
adjacent to at least one vertex in $S$.

Suppose that $n+m$ is even and that there exist
a subset $C\subseteq Q_0\sqcup\pp$ with $|C|=m$ and a matching
$M$ from $\pp$ onto $C$ such that $\mu(M')\ne\mu(M)$ for every
matching $M'\ne M$ from $\pp$ onto $C$. Then $\sq$ is symplectic.
\end{proposition}
\begin{proof}
The first assertion follows from the determinant expansion and Theorem~\ref{thm:fullrank}. A nonzero determinant has a nonzero term, proving the second assertion. A unique monomial has coefficient $1$ or $-1$, proving the third assertion.
\end{proof}

The vertex columns of $\KQ$ are determined by the path weights.
The following lemma expresses linear independence of these
weights as a forest condition.

\begin{definition}\label{def:pathforest}
For $T\subseteq\pp$, let $F_T$ be the undirected multigraph
on $Q_0$ with one edge $\{s(p),t(p)\}$ for each $p\in T$.
Distinct paths with the same starting and terminal vertices
give parallel edges. Let $\Fmax=F_{\mathcal M(Q)}$.
\end{definition}
\begin{lemma}\label{lem:incidence}
The family $(\alpha_p)_{p\in T}$ is linearly independent if and only if $F_T$ is a forest.
\end{lemma}

\begin{proof}

Suppose first that $F_T$ contains a cycle. Write its successive
vertices as $v_1,\cdots,v_r,v_{r+1}=v_1$, and let
$p_1,\cdots,p_r$ be the distinct paths corresponding to its
successive edges. 
For each $i$, 
take $\delta_i=1$ if $p_i$ starts at $v_i$, and
$\delta_i=-1$ otherwise, so 
$\delta_i\alpha_{p_i}=\eps_{v_i}-\eps_{v_{i+1}}$. Then
$\sum_{i=1}^r\delta_i\alpha_{p_i}
=\sum_{i=1}^r(\eps_{v_i}-\eps_{v_{i+1}})=0$. Hence linear independence
implies that $F_T$ is a forest.

Conversely, suppose that $F_T$ is a forest, and consider a
linear relation $\sum_{p\in T}c_p\alpha_p=0$, with
$c_p\in\kk$. We prove by induction on $|T|$ that every
$c_p$ is zero. The assertion is trivial when $T$ is empty.
If $T$ is nonempty, the forest $F_T$ has a vertex $v$
incident to exactly one edge.
Let $p_0\in T$ correspond to the unique edge incident to $v$.
Evaluating the linear relation at $e_v$ gives
$\sum_{p\in T}c_p\alpha_p(e_v)=0$.
For $p\ne p_0$, the vertex $v$ is neither the starting nor
the terminal vertex of $p$, so $\alpha_p(e_v)=0$.
Since $\alpha_{p_0}(e_v)\in\{1,-1\}$, it follows that
$c_{p_0}=0$.

The remaining relation is
$\sum_{p\in T\setminus\{p_0\}}c_p\alpha_p=0$.
Deleting the edge corresponding to $p_0$ leaves a forest,
so the induction hypothesis gives $c_p=0$ for every
$p\in T\setminus\{p_0\}$. This proves linear independence.
\end{proof}

\begin{definition}\label{def:forestmatching}
Let $D(Q)$ be the graph on $\pp$ in which distinct paths $p,q$
are adjacent if and only if $pq$ or $qp$ is a path.
A forest-matching configuration is a pair $(M,T)$, where
$M$ is a matching in $D(Q)$, $T\subseteq\pp\setminus V(M)$,
and $F_T$ is a forest.
Here $V(M)$ is the set of paths that occur as vertices of
the edges in $M$. Define
\begin{equation*}
 \rho(Q)=\max_{(M,T)}\bigl(2|M|+|T|\bigr).
\end{equation*}
The configuration is a cover if $V(M)\sqcup T=\pp$.
Thus a cover exists exactly when $\rho(Q)=m(Q)$.
\end{definition}

\begin{theorem}\label{thm:rho}
For every finite acyclic quiver, $\gamma(Q)\le\rho(Q)$. Consequently, full row rank of $\KQ$ implies that $Q$ has a forest-matching cover.
\end{theorem}

\begin{proof}
We work over $\kk(\mathbf X)$, extending all vector
spaces and forms from $\kk$ to $\kk(\mathbf X)$.
Let $m=m(Q)$ and $\eta=d\lambda_{\mathbf X}$.
Let $D=\eta|_{\nq\times\nq}$, and define
$\psi:\nq\to\hq^*$ by $\psi(x)(h)=\eta(x,h)$.
Consider the linear map $L:\nq\to\sq^*$ given by
$L(x)=\eta(x,\cdot)$. Since
$\eta(x,y)=-\lambda_{\mathbf X}([x,y])$, its rank is
$\gamma(Q)$.

Write $R=\rad D$. An element $x\in\nq$ belongs to $\ker L$
if and only if $D(x,y)=0$ for every $y\in\nq$ and
$\psi(x)(h)=0$ for every $h\in\hq$. Thus
$\ker L=R\cap\ker\psi=\ker(\psi|_R)$.
Write $\rank D=2r$ and $\rank(\psi|_R)=d$.
Then $\dim R=m-2r$ and
$\dim\ker L=m-2r-d$, so $\gamma(Q)=2r+d$.

Extend $D$ to $\sq$ by setting it equal to zero whenever
one of its two entries belongs to $\hq$.
Define an alternating two-form $\eta_{\mathrm{mix}}$ on $\sq$ by
$\eta_{\mathrm{mix}}(x,h)=\psi(x)(h)$ for $x\in\nq$, $h\in\hq$,
and by requiring its restrictions to $\nq\times\nq$ and
$\hq\times\hq$ to be zero.
We first show that $D^r\wedge\eta_{\mathrm{mix}}^d\ne0$.

Choose a complement $E$ to $R$ in $\nq$.
The restriction $D|_E$ is nondegenerate.
Hence $\dim E=2r$, and for a basis
$u_1,\cdots,u_{2r}$ of $E$ we have
$D^r(u_1,\cdots,u_{2r})\ne0$.
Choose $z_1,\cdots,z_d\in R$ such that
$\psi(z_1),\cdots,\psi(z_d)$ are linearly independent.
The matrix of their values on the vertex basis of $\hq$
has row rank $d$. We therefore choose $d$ elements of this
vertex basis, $h_1,\cdots,h_d$ such that the matrix
$A=(\psi(z_i)(h_j))_{i,j=1}^d$ has nonzero determinant.
Evaluate $D^r\wedge\eta_{\mathrm{mix}}^d$ on
$u_1,\cdots,u_{2r},z_1,\cdots,z_d,h_1,\cdots,h_d$.
Since $D$ vanishes whenever one entry lies in $R$ or $\hq$,
every nonzero term in this evaluation assigns all the
$u_i$ to $D^r$. The remaining entries are assigned to
$\eta_{\mathrm{mix}}^d$. Since $\eta_{\mathrm{mix}}$ vanishes on
$\nq\times\nq$ and $\hq\times\hq$, these entries must be
paired as $z_i,h_j$. Consequently,
\[
 (D^r\wedge\eta_{\mathrm{mix}}^d)
 (u_1,\cdots,u_{2r},z_1,\cdots,z_d,h_1,\cdots,h_d)
 =
 \pm d!\,D^r(u_1,\cdots,u_{2r})\det A\ne0.
\]

For each $p\in\pp$, let $p^*$ be the dual basis
 of $p$ in $\nq^*$, extended by zero on $\hq$.
Also extend each $\alpha_p$ by zero on $\nq$.
The nonzero terms in the basis expansion of $D$ have
the form $c_{p,q}p^*\wedge q^*$, where
$c_{p,q}\in\kk(\mathbf X)$ is nonzero and either $pq$
or $qp$ is a path. Moreover,
$\eta(p,h)=X_p\alpha_p(h)$ gives
$\eta_{\mathrm{mix}}=\sum_{p\in\pp}X_p p^*\wedge\alpha_p$.

Expand $D^r\wedge\eta_{\mathrm{mix}}^d$ using these expressions.
Since it is nonzero, at least one term in the expansion
is nonzero. Such a term is a nonzero coefficient times
\[
 p_1^*\wedge q_1^*\wedge\cdots\wedge p_r^*\wedge q_r^*
 \wedge t_1^*\wedge\cdots\wedge t_d^*
 \wedge\alpha_{t_1}\wedge\cdots\wedge\alpha_{t_d},
\]
where either $p_iq_i$ or $q_ip_i$ is a path for each $i$.
The paths $p_1,q_1,\cdots,p_r,q_r,t_1,\cdots,t_d$
must all be distinct, since a repeated dual basis
element would make this exterior product zero.
Also, $\alpha_{t_1},\cdots,\alpha_{t_d}$ must be linearly
independent. Their coordinates lie in $\kk$, so their
linear independence over $\kk(\mathbf X)$ is equivalent
to linear independence over $\kk$.

Set $M=\{\{p_i,q_i\}\mid 1\le i\le r\}$ and
$T=\{t_1,\cdots,t_d\}$.
Then $M$ is a matching in $D(Q)$ and
$T\subseteq\pp\setminus V(M)$.
Lemma~\ref{lem:incidence} shows that $F_T$ is a forest.
Thus $(M,T)$ is a forest-matching configuration with
$2|M|+|T|=2r+d=\gamma(Q)$, proving
$\gamma(Q)\le\rho(Q)$.

If $\gamma(Q)=m(Q)$, then $\rho(Q)=m(Q)$. This configuration uses all
non-trivial paths and is therefore a cover.
\end{proof}

We now determine the restrictions that a forest-matching cover
imposes on the arrows and vertices of $Q$.

\begin{lemma}\label{lem:maxforest}
If $Q$ has a forest-matching cover or $\gamma(Q)=m(Q)$,
then $\Fmax$ is a forest and 
$(\alpha_p)_{p\in\mathcal M(Q)}$ is linearly independent.
Moreover, for every symplectic form $\omega=\beta+d\lambda$
on $\sq$, we have $\lambda(p)\ne0$ for all
$p\in\mathcal M(Q)$.
\end{lemma}
\begin{proof}
A maximal path has no neighbor in $D(Q)$, so it belongs
to $T$ for every forest-matching cover $(M,T)$. Thus $\mathcal M(Q)\subseteq T$, hence $\Fmax$ is a forest. If $\gamma(Q)=m(Q)$, Theorem~\ref{thm:rho} gives a
forest-matching cover, so $\Fmax$ is also a forest. Lemma~\ref{lem:incidence} then shows that
$(\alpha_p)_{p\in\mathcal M(Q)}$ is linearly independent.

For a maximal path $p$, the form $\omega=\beta+d\lambda$
satisfies $\omega(p,\nq)=0$ and
$\omega(p,h)=\lambda(p)\alpha_p(h)$ for $h\in\hq$.
Thus $\lambda(p)=0$ would imply $p \in \rad\omega$,
contrary to nondegeneracy.

\end{proof}

The next identity transfers the forest condition from maximal
paths to the underlying multigraph.
To eliminate an internal vertex $v$ of current indegree $r$
and outdegree $s$, remove $v$ and replace each path of length
two through $v$ by an arrow with the same starting and terminal vertices.
\begin{proposition}\label{prop:elimination}
Eliminate all original internal vertices in this way.
At step $i$, let $r_i,s_i$ be the current degrees
of the vertex removed. Then
$b_1(\Fmax)=b_1(U(Q))+\sum_i(r_i-1)(s_i-1)$.
In particular, if $\Fmax$ is a forest, then $U(Q)$ is a forest.
\end{proposition}
\begin{proof}
A current arrow records an original path whose internal vertices have been removed. Concatenating at the next removed vertex preserves this correspondence. Acyclicity persists, and every remaining original internal vertex still lies between an original source and an original sink, so its current indegree and outdegree are positive.

In the underlying multigraph, removing a vertex of degrees $r,s$ changes the edge count by $rs-r-s$ and the vertex count by $-1$. The number of connected components is unchanged. Hence the cycle rank increases by $(r-1)(s-1)$. At the end, the arrows correspond precisely to the original maximal paths. Adding the removed vertices back as isolated vertices gives $\Fmax$ without changing cycle rank. Summing these changes gives the stated equality.  All its terms on the right are nonnegative. So, if $b_1(\Fmax)=0$, then $b_1(U(Q))=0$, which proves the last assertion.
\end{proof}

\begin{proposition}\label{prop:underlyingforest}
If $Q$ has a forest-matching cover or $\gamma(Q)=m(Q)$, then $U(Q)$ is a forest. Thus $Q$ has no parallel arrows, has at most one path between any two vertices, and $m(Q)\le\sum_C\binom{|C|}2\le\binom{n(Q)}2$, where $C$ runs through the connected components of $U(Q)$. If two distinct paths have the same starting and terminal vertices, then $\sq$ is not symplectic. 
\end{proposition}
\begin{proof}
The first assertion follows from Lemma~\ref{lem:maxforest} and Proposition~\ref{prop:elimination}. Since $U(Q)$ is a forest, $Q$ has no parallel arrows and has at most one path between any two vertices. A forest has a unique undirected path between any two vertices in a connected component, so the connected component $C$ 
contains at most $\binom{|C|}{2}$ non-trivial paths.
Summing over the connected components gives
$m(Q)\le\sum_C\binom{|C|}{2}\le\binom{n(Q)}{2}$. 

Two distinct paths with the same starting and terminal vertices extend to two distinct maximal paths with the same starting and terminal vertices, giving parallel edges in $\Fmax$. So, $\Fmax$ is not a forest. By Theorem~\ref{thm:fullrank} and Lemma~\ref{lem:maxforest}, $\sq$ is not symplectic. 
\end{proof}

Once $U(Q)$ is a forest, separating the arrows at sources and
sinks decomposes both the rank problem and the cover problem
into the following blocks.

\begin{definition}\label{def:blocks}
Suppose that $U(Q)$ is a forest.
For $a,b\in Q_1$, define $a\sim b$ if there exist arrows
$a_0,\cdots,a_\ell$, with $\ell\ge0$, such that
$a_0=a$, $a_\ell=b$, and each consecutive pair
$a_{i-1},a_i$ has a common vertex which is an internal vertex in $Q$.
This is an equivalence relation on $Q_1$.

For each equivalence class $E$, let $B$ be the subquiver
with arrow set $E$ and vertex set
$\{s(a)\mid a\in E\}\cup\{t(a)\mid a\in E\}$,
with starting and terminal vertices inherited from $Q$.
We call $B$ a block of $Q$.

An arrow directly from a source to a sink is a one-arrow block.
Isolated vertices are not blocks.
Distinct blocks may share an original source or sink.
\end{definition}

\begin{lemma}\label{lem:blocks}
The underlying undirected multigraph of each block is a tree.
Its sources and sinks have degree one in the block,
and all other vertices are internal.
Every non-trivial path lies in one block, and non-trivial  paths from different blocks cannot concatenate. Moreover, $\gamma(Q)=\sum_B\gamma(B)$ and
$\rho(Q)=\sum_B\rho(B)$.
In particular, a cover exists on $Q$ if and only if it exists on every block.
\end{lemma}
\begin{proof}
All arrows incident to an internal vertex belong to the same equivalence class. If two arrows at a source (respectively, sink) belonged to one equivalence class, 
the definition of the equivalence relation, together with the fact that the source (respectively, sink) is not internal, would give an undirected path between their terminal vertices (respectively, starting vertices) avoiding that source (respectively, sink). Together with the two arrows, this would form an undirected cycle.
Thus sources and sinks have degree one in each block.

All arrows of the same path
belong to the same equivalence class. Hence every non-trivial
path lies in a unique block. If two non-trivial paths
concatenate, their concatenation lies in one block,
so both paths lie in that block. Thus non-trivial paths from different
blocks cannot concatenate.

The arrow weights of $Q$ are independent by Lemma~\ref{lem:incidence}. For each block $B$, let
$W_B=\Span_{\kk}\{\alpha_a\mid a\in B_1\}\subseteq\hq^*$.
 Thus the spaces $W_B$ form
a direct sum.
If $p=a_1\cdots a_\ell$ is a path in $B$, then
$\alpha_p=\sum_{i=1}^{\ell}\alpha_{a_i}$, so
$\alpha_p\in W_B$. Let $H_B\subseteq\mathcal P(B)$
for each block, and let $H=\bigcup_BH_B$. Since the spaces $W_B$ form
a direct sum,
 $(\alpha_p)_{p\in H}$ is linearly independent
if and only if each 
$(\alpha_p)_{p\in H_B}$ is linearly independent.
By Lemma~\ref{lem:incidence}, it follows that $F_H$
is a forest if and only if every $F_{H_B}$ is a forest.

Let $(M,T)$ be a forest-matching configuration on $Q$.
Every edge of $M$ joins paths in the same block.
For each block $B$, let $M_B=M\cap E(D(B))$
and let $T_B=T\cap\mathcal P(B)$.
Then $(M_B,T_B)$ is a forest-matching configuration on $B$,
so $2|M|+|T|=\sum_B(2|M_B|+|T_B|)\le\sum_B\rho(B)$.
Conversely, choose a forest-matching configuration
$(M_B,T_B)$ on each block $B$ such that
$2|M_B|+|T_B|=\rho(B)$.
Then $M=\bigcup_BM_B$ is a matching in $D(Q)$,
and $T=\bigcup_BT_B$ is disjoint from $V(M)$.
The preceding independence criterion shows that $F_T$ is a forest.
Thus $(M,T)$ is a forest-matching configuration on $Q$
with $2|M|+|T|=\sum_B\rho(B)$.
Therefore $\rho(Q)=\sum_B\rho(B)$.

For the rank assertion, work over $\kk(\mathbf X)$ and let
$\eta=d\lambda_{\mathbf X}$.
Write $x=\sum_Bx_B$ with $x_B\in\mathfrak n_B$.
Since paths in different blocks cannot concatenate,
$\eta(x,\nq)=0$ if and only if
$\eta(x_B,\mathfrak n_B)=0$ for every $B$.
If $x_B=\sum_{p\in\mathcal P(B)}c_pp$, then
$\eta(x_B,\cdot)|_{\hq}=\sum_{p\in\mathcal P(B)}c_pX_p\alpha_p$,
which belongs to the $\kk(\mathbf X)$-span of $W_B$.
These spans form a direct sum.
For each $B$, every form in the $\kk(\mathbf X)$-span of $W_B$
vanishes on $e_v$ for $v\notin B_0$.
Hence $\eta(x,\hq)=0$ if and only if
$\eta(x_B,\mathfrak h_B)=0$ for every $B$.

It follows that the kernel of
$x\mapsto\eta(x,\cdot)$ from $\nq$ to $\sq^*$
is the direct sum of the kernels of the maps
$x_B\mapsto\eta(x_B,\cdot)|_{\mathfrak s_B}$.
Each block map has rank $\gamma(B)$, since extending
its coefficient field to $\kk(\mathbf X)$ does not change
its matrix rank. 
Since $\nq=\bigoplus_B\mathfrak n_B$ and the kernels
decompose in the same way, we obtain
$\gamma(Q)=\sum_B\gamma(B)$.
The assertion about covers follows from $m(Q)=\sum_Bm(B)$.

\end{proof}

We now restrict to a block $B$ with a cover $(M,T)$.
Every source and sink of $B$ is incident to exactly one
arrow of $B$. For vertices $x,y$ of $B$, write
$x\rightsquigarrow y$ for the unique directed path
in $B$ from $x$ to $y$, when such a path exists.
Every maximal path of $B$ belongs to $T$.

\begin{lemma}\label{lem:twoextremes}
Every source of $B$ reaches at most two sinks, and every sink is reached from at most two sources.
\end{lemma}
\begin{proof}
Suppose a source $s$ reaches three distinct sinks $t_1,t_2,t_3$. Let $a$ be its unique outgoing arrow, from $s$ to $u$. If $u$ had another incoming arrow, its incoming branch would contain a different source $s'$. Both sources would reach all three sinks, and the six
corresponding maximal paths would create a cycle in $F_T$. Thus $d^-(u)=1$, and the only non-trivial path ending at $u$ is $a$.

Each $u\rightsquigarrow t_i$, denoted by $p_i$, has $a$ as its only neighbor in $D(B)$.
Since $M$ is a matching, at most one of the edges
$\{a,p_1\}$, $\{a,p_2\}$, and $\{a,p_3\}$ belongs to $M$. Hence at least two of these paths belong to $T$.
Without loss of generality, assume that $p_1,p_2\in T$.
The maximal paths $s\rightsquigarrow t_i$ where $i=1,2$ also belong to $T$. The two paths and $p_1,p_2$ create a cycle in $F_T$, a contradiction. Reversing arrows proves the other assertion.
\end{proof}

A branch vertex means a vertex of total degree at least three. A vertex of total degree three is said to be a trivalent vertex. An internal trivalent vertex is a split if its degree pair is $(1,2)$ and a merge if it is $(2,1)$.
\begin{lemma}\label{lem:onebranch}
Every internal vertex of $B$ has total degree at most three, and $B$ has at most one branch vertex.
\end{lemma}
\begin{proof}
Outdegree at least three gives three distinct sinks, all reachable from a source, contrary to Lemma~\ref{lem:twoextremes}. Indegree at least three is excluded dually. A vertex with degree pair $(2,2)$ gives two sources and two sinks, whose four maximal paths form a cycle in $F_T$. Thus every internal vertex has total degree at most three.

Suppose there are at least two branch vertices. Choose distinct branch vertices \(x,y\) such that the unique undirected path between them contains no other branch vertex. The intermediate vertices have degree pair $(1,1)$. Thus, after interchanging $x$ and $y$ if necessary,
this path is directed from $x$ to $y$. A split followed by a split gives three sinks reachable from a source. A merge followed by a merge is excluded dually. A merge followed by a split gives a cycle in $F_T$. It remains to consider a split followed by a merge.

Choose a source $s$ reaching the split $x$ and a sink $t$
reachable from the merge $y$.
There is another sink $t'\ne t$ reachable from $x$
and another source $s'\ne s$ reaching $y$.
The maximal paths $s\rightsquigarrow t$,
$s\rightsquigarrow t'$, and $s'\rightsquigarrow t$
belong to $T$.
Let $\Gamma$ be the connected subgraph of $F_T$
consisting of their three edges.

Let $a$ be the first arrow from $s$, with terminal vertex $u$.
If $u$ had another incoming arrow, a different source
would reach both $t,t'$, producing a cycle in $F_T$.
Hence the only non-trivial path ending at $u$ is $a$.

Each of $u\rightsquigarrow t$ and
$u\rightsquigarrow t'$ has $a$ as its only neighbor
in $D(B)$.
Since $M$ is a matching, at most one of
$\{a,u\rightsquigarrow t\}$ and
$\{a,u\rightsquigarrow t'\}$ belongs to $M$.
Consequently, at least one of the two paths belongs to $T$.
They cannot both belong to $T$, since their edges,
together with $\{s,t\}$ and $\{s,t'\}$,
would give a cycle in $F_T$.
Hence exactly one of $u\rightsquigarrow t$ and
$u\rightsquigarrow t'$ belongs to $T$,
and the other forms an edge of $M$ with $a$.
Add the edge corresponding to the path in $T$
to $\Gamma$.
Then $\Gamma$ is connected and contains $u$,
and $a$ belongs to an edge of $M$.

Dually, let $c$ be the last arrow to $t$, with starting vertex $v$.
The only non-trivial path starting at $v$ is $c$.
Exactly one of $s\rightsquigarrow v$ and
$s'\rightsquigarrow v$ belongs to $T$,
and the other forms an edge of $M$ with $c$.
Add the edge corresponding to the path in $T$
to $\Gamma$.
Then $\Gamma$ is connected and contains both $u$ and $v$,
and $c$ belongs to an edge of $M$.

The path $u\rightsquigarrow v$ contains
$x\rightsquigarrow y$, so it has positive length.
Its only neighbors in $D(B)$ are $a$ and $c$.
Both already belong to edges of $M$ with paths
different from $u\rightsquigarrow v$.
Therefore $u\rightsquigarrow v$ belongs to $T$.
Its edge $\{u,v\}$ is not among the edges used
to construct $\Gamma$, but $\Gamma$ already connects
$u$ and $v$, which creates a cycle
in $F_T$. This is a contradiction.

\end{proof}

It follows that $B$ either consists of a single directed
path or has precisely one trivalent vertex $w$.
In the latter case, $B$ is the union of three non-trivial
directed paths, any two of which have $w$ as their
only common vertex.
After reversing all arrows if necessary, one of these
paths ends at $w$, and the other two start at $w$.
We denote this quiver by $Y(a;b,c)$, where $a$ is the
length of the path ending at $w$, and $b,c$ are the
lengths of the two paths starting at $w$.
In particular, $a,b,c\ge1$.

\begin{lemma}\label{lem:shortarm}
If $Y(a;b,c)$ has a forest-matching cover, then $\min\{a,b,c\}=1$.
\end{lemma}
\begin{proof}
Suppose $a,b,c\ge2$. Let $s$ be the source, $u$ the first vertex after $s$, and $g$ the second vertex after $s$. Let $t_1,t_2$ be the sinks and $v_1,v_2$ the first vertices before $t_1,t_2$, respectively. Then $g,v_1,v_2$ are distinct, and all paths below have positive length. Let $a_0$ be the arrow from $s$ to $u$, and let $c_i$ be the arrow from $v_i$ to $t_i$ for $i=1,2$.

The maximal paths $s\rightsquigarrow t_1$ and
$s\rightsquigarrow t_2$ belong to $T$.
Let $\Gamma$ be the connected subgraph of $F_T$
consisting of their two edges.
Each of $u\rightsquigarrow t_1$ and
$u\rightsquigarrow t_2$ has $a_0$ as its only neighbor
in $D(B)$.
By the same argument as in the proof of
Lemma~\ref{lem:onebranch}, exactly one of these paths
belongs to $T$, and the other forms an edge of $M$
with $a_0$.
Add the edge corresponding to the path in $T$
to $\Gamma$.
Then $\Gamma$ is connected and contains $u$,
and $a_0$ belongs to an edge of $M$.

For each $i=1,2$, the only neighbor of
$s\rightsquigarrow v_i$ in $D(B)$ is $c_i$,
whereas the only neighbors of $u\rightsquigarrow v_i$
are $a_0$ and $c_i$.
Since $a_0$ already belongs to an edge of $M$
with a different path, 
at most one of
$s\rightsquigarrow v_i$ and $u\rightsquigarrow v_i$
can belong to $V(M)$.
They cannot both belong to $T$, since $\Gamma$
already connects $s$ and $u$.
Hence exactly one belongs to $T$,
and the other forms an edge of $M$ with $c_i$.
Add the edge corresponding to the path in $T$
to $\Gamma$ for each $i$.
Then $\Gamma$ is connected and contains
$s,u,t_1,t_2,v_1,v_2$, and both $c_1,c_2$
belong to edges of $M$.

Consider the four paths $g\rightsquigarrow v_1$,
$g\rightsquigarrow t_1$, $g\rightsquigarrow v_2$,
and $g\rightsquigarrow t_2$.
The only non-trivial paths ending at $g$ are
$u\rightsquigarrow g$ and $s\rightsquigarrow g$.
The only non-trivial path starting at $v_i$ is $c_i$,
and no non-trivial path starts at $t_i$.
Since $c_1,c_2$ already belong to edges of $M$
with different paths, any edge of $M$ incident to
one of the four paths must have
$u\rightsquigarrow g$ or $s\rightsquigarrow g$
as its other vertex.
Thus at most two of the four paths belong to $V(M)$,
so at least two belong to $T$.
Their edges join $g$, which is not a vertex of $\Gamma$,
to two distinct vertices of $\Gamma$.
Since $\Gamma$ is connected, adding these two edges
creates a cycle in $F_T$, a contradiction.

\end{proof}

\begin{proposition}\label{prop:necessaryshape}
A quiver admitting a forest-matching cover satisfies condition (iii)
of Theorem~\ref{thm:fullclassification}.
\end{proposition}
\begin{proof}
Proposition~\ref{prop:underlyingforest} shows that $U(Q)$ is a forest. Restrict a cover to each block by Lemma~\ref{lem:blocks}, then apply Lemma~\ref{lem:onebranch} and Lemma~\ref{lem:shortarm}. 
The condition $\min\{a,b,c\}=1$ in Lemma~\ref{lem:shortarm} holds exactly when
the trivalent vertex is adjacent to a source or sink
of $Q$.
The vertices of every connected component of $U(Q)[I(Q)]$
belong to one block, and their degrees in that block
are the same as in $Q$.
One-arrow blocks have no internal vertex.
\end{proof}

To prove sufficiency, we introduce a terminal-path reduction
that preserves both symplecticity and the difference
$m(Q)-\gamma(Q)$. 
It constructs forms on the permitted blocks and yields a
recursive criterion for arbitrary quivers.
Section~\ref{sec:solvablegroups} identifies this operation with
symplectic reduction on the corresponding Lie groups.
\begin{definition}\label{def:terminal}
A reducible terminal-path is a path $P=a_1\cdots a_\ell$,
with $\ell\ge1$, where $a_i$ is an arrow from $v_{i-1}$
to $v_i$ for $1\le i\le\ell$, such that
$s=v_0$ is a source, $t=v_\ell$ is a sink,
and $d^-(v_i)=1$ for $1\le i\le\ell$.
Define $R_P(Q)$ as follows:
\begin{enumerate}
\item delete $t$ and its unique incoming arrow $a_\ell$;
\item if $\ell\ge2$, delete $a_1$ as well;
\item for each $1\le i<\ell$ and each arrow $b$ from $v_i$ to $z$ other than $a_{i+1}$, retain $b$ and add a new arrow $b'$ from $s$ to $z$.
\end{enumerate}
All other vertices and arrows remain. A source reduction of $Q$ is defined as follows.
First reverse all arrows of $Q$ to obtain $Q^{\mathrm{op}}$.
Choose a reducible terminal-path $g$ in $Q^{\mathrm{op}}$
and construct $R_g(Q^{\mathrm{op}})$.
Then reverse all arrows of $R_g(Q^{\mathrm{op}})$.
 
\end{definition}
Every new arrow represents the original path $a_1\cdots a_i b$,
so $R_P(Q)$ is acyclic.
Constructions representing paths by edges or arrows also occur in
\cite{Mizoguchi} and \cite{MT}.
The deletion rule here is chosen to preserve the  difference $m(Q)-\gamma(Q)$  and symplecticity, as proved in Theorem~\ref{thm:reduction}.

Choose such a path $P$ and write $R=R_P(Q)$. Let
\begin{equation*}
 q_i=a_{i+1}\cdots a_\ell\quad(0\le i<\ell),\qquad
 r_i=a_1\cdots a_i\quad(1\le i<\ell).
\end{equation*}
These are distinct non-trivial paths, with $q_0=P$ and $r_iq_i=q_0$. Let $W=\mathcal P(Q)\setminus
\bigl(\{q_i\mid 0\le i<\ell\}\cup\{r_i\mid 1\le i<\ell\}\bigr)$.

\begin{lemma}\label{lem:reductionembedding}
Replacing each new arrow $b'$ by $a_1\cdots a_i b$ defines a bijection $\mathcal P(R)\to W$ preserving concatenation. There is an injective homomorphism $\iota:\kk R\to\kk Q$ whose image has basis
$
 W\ \sqcup\ \{e_s+e_t\}\ \sqcup\ \{e_v\mid v\ne s,t\}$.
It sends the trivial path at $s$ to $e_s+e_t$ and the other remaining trivial paths to themselves. Moreover,
$m(Q)=m(R)+2\ell-1$, and
 $\dim\sq=\dim\mathfrak s_R+2\ell$.

\end{lemma}
\begin{proof}
The only non-trivial paths ending at $t$ are the $q_i$, because all indegrees along $P$ are one. Every path ending at an internal vertex of $P$ is a segment of $P$. A path leaving $P$ along an arrow not belonging to $P$
cannot return to $P$, since every arrow ending at
a vertex of $P$ is one of $a_1,\cdots,a_\ell$.

A path $p\in W$ starting at $s$ and using $a_1$ cannot stop at an internal vertex of $P$ or at $t$. It first leaves $P$ at an arrow $b$ which does not belong to $P$. Replacing that initial segment by $b'$ gives its unique corresponding path in $R$. Every other path in $W$ consists entirely of arrows
that remain in $R$, so it is a path in $R$
and corresponds to itself. This describes the inverse bijection.

The replacement preserves starting and terminal vertices, none of the paths in $W$ ends at $t$, and no non-trivial path starts at the sink $t$. Hence it preserves both nonzero and zero products. The product of two paths in $W$ is either zero
or another path in $W$.
The idempotent $e_s+e_t$ acts on paths in $W$ as
the trivial path at $s$ in $R$,
and the other vertex actions are unchanged. This shows that $\iota$ is an injective homomorphism. The path count follows from the bijection $\mathcal P(R)\to W$. One vertex is deleted, giving the dimension formula.
\end{proof}

\begin{lemma}\label{lem:normalization}
Let $\lambda\in\nq^*$ satisfy $\lambda(q_0)\ne0$. There exist an inner automorphism $\varphi$ of $\kk Q$
with $\varphi(\nq)=\nq$ and a nonzero $c\in\kk$
such that the linear map $\widetilde\lambda\in\nq^*$,
defined by $\widetilde\lambda(x)=c\lambda(\varphi(x))$
for $x\in\nq$, satisfies
$\widetilde\lambda(q_0)=1$ and
$\widetilde\lambda(q_i)=0$ for $1\le i<\ell$. The matrices $\mathcal K_Q(\lambda)$ and
$\mathcal K_Q(\widetilde\lambda)$ have the same rank.
For every closed two-form $\omega=\beta+d\lambda$,
the two-form $\widetilde\omega=c\beta+d\widetilde\lambda$
satisfies
$\widetilde\omega(x,y)=c\omega(\varphi(x),\varphi(y))$
for all $x,y\in\sq$.
Consequently, $\widetilde\omega$ is closed, and $\widetilde\omega$ is
nondegenerate if and only if $\omega$ is nondegenerate.

\end{lemma}
\begin{proof}
For $1\le i<\ell$, let
$c_i=\lambda(q_i)/\lambda(q_0)$, and set
$u=\sum_{i=1}^{\ell-1}c_ir_i$.
 $r_ir_j=0$ for all $1\le i,j<\ell$.
Hence $u^2=0$, and $g=1+u$ has inverse $1-u$.
Define the inner automorphism $\varphi$ of $\kk Q$
by $\varphi(x)=g^{-1}xg$.
Since $\kk Q_{\ge1}$ is a two-sided ideal,
$\varphi(\nq)=\nq$.

It is easy to check that $\varphi(q_0)=q_0$ and
$\varphi(q_i)=q_i-c_iq_0$ for $1\le i<\ell$.
Set $c=\lambda(q_0)^{-1}$ and define
$\widetilde\lambda(x)=c\lambda(\varphi(x))$
for $x\in\nq$.
Then $\widetilde\lambda(q_0)=1$, and
$\widetilde\lambda(q_i)
=c(\lambda(q_i)-c_i\lambda(q_0))=0$
for $1\le i<\ell$.

For $x\in\nq$ and $y\in\sq$, we have
$B_{\widetilde\lambda}(x,y)
=c\lambda(\varphi([x,y]))
=cB_\lambda(\varphi(x),\varphi(y))$.
The restrictions of $\varphi$ to $\nq$ and $\sq$
are invertible, and $c\ne0$.
Therefore $\mathcal K_Q(\widetilde\lambda)$ and
$\mathcal K_Q(\lambda)$ have the same rank.

Finally, let $\omega=\beta+d\lambda$ be a closed two-form.
Since $u\in\nq$, we have $\varphi(x)-x\in\nq$
for every $x\in\sq$.
The form $\beta$ vanishes whenever one of its two
entries belongs to $\nq$, so
$\beta(\varphi(x),\varphi(y))=\beta(x,y)$.
Moreover, since $[\sq,\sq]\subseteq\nq$,
we have
$(d\widetilde\lambda)(x,y)
=-c\lambda(\varphi([x,y]))
=c(d\lambda)(\varphi(x),\varphi(y))$.
Consequently,
$(c\beta+d\widetilde\lambda)(x,y)
=c\omega(\varphi(x),\varphi(y))$
for all $x,y\in\sq$.
Since $\varphi$ is also a Lie algebra automorphism,
this two-form is closed.
Since $\varphi$ is invertible and $c\ne0$,
it is nondegenerate if and only if $\omega$
is nondegenerate.
\end{proof}

A linear map $\lambda\in\nq^*$ satisfying $\lambda(q_0)=1$
and $\lambda(q_i)=0$ for $1\le i<\ell$ is said to be normalized.

Define $\lambda_R\in\mathfrak n_R^*$ by
$\lambda_R(p)=\lambda(\iota(p))$
for every $p\in\mathcal P(R)$. Let $b_1,\cdots,b_N$ be the path basis of $\mathfrak s_R$,
where $N=\dim\mathfrak s_R$.

\begin{lemma}\label{lem:reductionblocks}
For a normalized $\lambda$, the two vector spaces $
\nq\cap\rad B_\lambda$ and $
\mathfrak n_R\cap\rad B_{\lambda_R}$ are isomorphic. 
For every closed form $\omega=\beta+d\lambda$, its matrix with respect to the basis
$(q_0,\cdots,q_{\ell-1},e_t,r_1,\cdots,r_{\ell-1},
\iota(b_1),\cdots,\iota(b_N))$ has the form
\begin{equation*}
 \begin{pmatrix}
 0&P_\ell&0\\
 -P_\ell^T&0&\Theta\\
 0&-\Theta^T&\Omega_R
 \end{pmatrix},
 \qquad P_\ell=\operatorname{diag}(-1,1,\cdots,1),
\end{equation*}
where $\Omega_R$ is the matrix of $\iota^*\omega$
in the basis $(b_1,\cdots,b_N)$  and $\Theta$ is the $\ell\times N$ matrix defined by
$\Theta_{1j}=\omega(e_t,\iota(b_j))$ for $1\le j\le N$
and $\Theta_{i+1,j}=\omega(r_i,\iota(b_j))$
for $1\le i<\ell$ and $1\le j\le N$. 
The displayed block matrix has rank $2\ell+\rank\Omega_R$
and determinant $\det\Omega_R$.
\end{lemma}
\begin{proof}
For $i\ge1$, the only nonzero value of $B_\lambda(q_i,p)$ on a non-trivial path $p$ is $B_\lambda(q_i,r_i)=-1$. Also, $B_\lambda(q_0,e_t)=1$ and $B_\lambda(q_i,e_t)=0$ for $i\ge1$. For every $0\le i<\ell$ and $z\in\mathfrak s_R$,
 $B_\lambda(q_i,\iota(z))=0$.
Indeed, $q_0$ commutes with every path in $W$,
with $e_s+e_t$, and with $e_v$ for $v\ne s,t$.
For $1\le i<\ell$, since
$r_i\notin W$, $B_\lambda(q_i,p)=0$ for every $p\in W$.
Moreover, $\lambda(q_i)=0$ implies
$B_\lambda(q_i,e_v)=0$ for every vertex $v$,
and hence also $B_\lambda(q_i,e_s+e_t)=0$.

Write $x=\sum_{i=0}^{\ell-1}f_iq_i
+\sum_{i=1}^{\ell-1}g_ir_i+w$,
where $f_i,g_i\in\kk$ and $w\in\Span W$.
Suppose first that $x\in\nq\cap\rad B_\lambda$.
The equations $B_\lambda(x,e_t)=0$ and $B_\lambda(x,q_i)=0$ give $f_0=0$ and all $g_i=0$.
So, $x=\sum_{i=1}^{\ell-1}f_iq_i+w$.
Let $z=\iota^{-1}(w)\in\mathfrak n_R$.
For every $y\in\mathfrak s_R$, we have
$B_\lambda(x,\iota(y))=B_\lambda(w,\iota(y))$,
since $B_\lambda(q_i,\iota(y))=0$.
Since $\iota$ preserves brackets and
$\lambda_R=\lambda\circ\iota$ on $\mathfrak n_R$,
$B_\lambda(w,\iota(y))=B_{\lambda_R}(z,y)$.
Hence $B_\lambda(x,\iota(y))=0$ for every
$y\in\mathfrak s_R$ is equivalent to
$z\in\mathfrak n_R\cap\rad B_{\lambda_R}$.
For $1\le i<\ell$, the equation
$B_\lambda(x,r_i)=0$ becomes
$-f_i+B_\lambda(w,r_i)=0$.
Therefore $f_i=B_\lambda(w,r_i)$,
so $x$ is uniquely determined by $w$.

Conversely, let
$z\in\mathfrak n_R\cap\rad B_{\lambda_R}$,
let $w=\iota(z)$, and define
$x=w+\sum_{i=1}^{\ell-1}B_\lambda(w,r_i)q_i$.
The same calculations show that $B_\lambda(x,e_t)=0$,
$B_\lambda(x,q_i)=0$ for $0\le i<\ell$,
$B_\lambda(x,r_i)=0$ for $1\le i<\ell$,
and $B_\lambda(x,\iota(y))=0$
for every $y\in\mathfrak s_R$.
These elements span $\sq$, so
$x\in\nq\cap\rad B_\lambda$.
Thus the map $x\mapsto\iota^{-1}(w)$ is an
isomorphism from $\nq\cap\rad B_\lambda$ to
$\mathfrak n_R\cap\rad B_{\lambda_R}$.

We now compute the matrix of $\omega=\beta+d\lambda$. It is easy to compute
that the matrix of $\omega$ has the stated block form,
with $P_\ell=\operatorname{diag}(-1,1,\cdots,1)$.
Let $N=\dim\mathfrak s_R$ and set
\[
S=
\begin{pmatrix}
I_\ell&0&(P_\ell^T)^{-1}\Theta\\
0&I_\ell&0\\
0&0&I_N
\end{pmatrix}.
\]
It is easy to check that
\[
S^T\Omega S=
\begin{pmatrix}
0&P_\ell&0\\
-P_\ell^T&0&0\\
0&0&\Omega_R
\end{pmatrix},
\]
where $\Omega$ is the matrix of $\omega$.
Since $\det S=1$ and the upper left
$2\ell\times2\ell$ block has rank $2\ell$
and determinant $\det(P_\ell)^2=1$,
we obtain $\rank\Omega=2\ell+\rank\Omega_R$
and $\det\Omega=\det\Omega_R$.
\end{proof}

The block matrix separates a nondegenerate part of dimension
$2\ell$ and leaves the rank obstruction on the reduced quiver.

\begin{theorem}\label{thm:reduction}
For a reducible terminal-path $P$ of length $\ell$ and $R=R_P(Q)$,
$ \gamma(Q)=2\ell-1+\gamma(R)$, $
 m(Q)-\gamma(Q)=m(R)-\gamma(R)$.
Moreover, $\sq$ is symplectic if and only if $\mathfrak s_R$ is symplectic. The same assertions hold for a source reduction.
\end{theorem}
\begin{proof}

Let $f(\mathbf X)$ be the determinant of a
$\gamma(Q)\times\gamma(Q)$ minor of
$\mathcal K_Q(\mathbf X)$ that is not the zero polynomial.
Then $X_{q_0}f(\mathbf X)$ is also nonzero.
 And there exist values
$t_p\in\kk$ such that
$t_{q_0}f((t_p)_{p\in\mathcal P(Q)})\ne0$.
Define $\lambda_0\in\nq^*$ by $\lambda_0(p)=t_p$
for every $p\in\mathcal P(Q)$.
Then $\lambda_0(q_0)\ne0$ and
$\rank\mathcal K_Q(\lambda_0)=\gamma(Q)$.
By Lemma~\ref{lem:normalization}, there exists
$\lambda\in\nq^*$ satisfying $\lambda(q_0)=1$,
$\lambda(q_i)=0$ for $1\le i<\ell$, and
$\rank\mathcal K_Q(\lambda)
=\rank\mathcal K_Q(\lambda_0)=\gamma(Q)$.

For every $\lambda\in\nq^*$ satisfying
$\lambda(q_0)=1$ and $\lambda(q_i)=0$
for $1\le i<\ell$, define
$\lambda_R(p)=\lambda(\iota(p))$
for $p\in\mathcal P(R)$.
By Lemma~\ref{lem:reductionblocks},
$\nq\cap\rad B_\lambda$ and
$\mathfrak n_R\cap\rad B_{\lambda_R}$
have the same dimension.
Since $m(Q)=m(R)+2\ell-1$, 
$\rank\mathcal K_Q(\lambda)
=2\ell-1+\rank\mathcal K_R(\lambda_R)$.
Applying this equality to the $\lambda$ chosen above
yields $\gamma(Q)\le2\ell-1+\gamma(R)$.

Conversely, choose $\lambda_R\in\mathfrak n_R^*$
such that $\rank\mathcal K_R(\lambda_R)=\gamma(R)$.
Define $\lambda$ on $W$ by
$\lambda(\iota(p))=\lambda_R(p)$
for every $p\in\mathcal P(R)$.
Set $\lambda(q_0)=1$ and $\lambda(q_i)=0$
for $1\le i<\ell$, and choose the values
$\lambda(r_i)$ arbitrarily.
These assignments define a linear map $\lambda\in\nq^*$.
The same rank equality gives
$\rank\mathcal K_Q(\lambda)=2\ell-1+\gamma(R)$,
so $\gamma(Q)\ge2\ell-1+\gamma(R)$.
Therefore $\gamma(Q)=2\ell-1+\gamma(R)$.
Together with $m(Q)=m(R)+2\ell-1$, this proves
$m(Q)-\gamma(Q)=m(R)-\gamma(R)$.
Since $\dim\sq=\dim\mathfrak s_R+2\ell$,
Theorem~\ref{thm:fullrank} shows that $\sq$ is symplectic
if and only if $\mathfrak s_R$ is symplectic.

We also give an explicit construction of the form. Recall that $t=v_\ell$ is the terminal vertex of $P$
and is deleted in the construction of $R$.
Let $\omega_R=\beta_R+d\lambda_R$ be a symplectic form
on $\mathfrak s_R$.
Define $\lambda$ on $W$ by
$\lambda(\iota(p))=\lambda_R(p)$
for every $p\in\mathcal P(R)$, and set
$\lambda(q_0)=1$ and
$\lambda(q_i)=\lambda(r_i)=0$ for $1\le i<\ell$.
For each vertex $v$ of $R$, let $e_v^R$ denote
the corresponding trivial path.
Define $\beta\in\Lambda^2\hq^*$ by
$\beta(e_v,e_w)=\beta_R(e_v^R,e_w^R)$ for $v,w\ne t$
and $\beta(e_t,h)=0$ for every $h\in\hq$.
Since $\iota(e_s^R)=e_s+e_t$ and
$\iota(e_v^R)=e_v$ for $v\ne s,t$,
it is easy to check that the closed form
$\omega=\beta+d\lambda$ satisfies
$\omega(\iota(x),\iota(y))=\omega_R(x,y)$
for all $x,y\in\mathfrak s_R$.
The determinant formula in
Lemma~\ref{lem:reductionblocks} therefore shows
that $\omega$ is nondegenerate.

Conversely, let $\omega=\beta+d\lambda$ be symplectic.
Since $q_0$ is maximal, Lemma~\ref{lem:maxforest}
gives $\lambda(q_0)\ne0$.
Let $\widetilde\omega$ be the symplectic form obtained
by Lemma~\ref{lem:normalization}.
Lemma~\ref{lem:reductionblocks} shows that
$\iota^*\widetilde\omega$ is symplectic on $\mathfrak s_R$.

The proof of the assertions for a source reduction is similar.

\end{proof}

The lifting construction also preserves the values $0$ and $1$.
Thus a map $\lambda_R$ taking these values on paths and
satisfying $\rank\mathcal K_R(\lambda_R)=m(R)$ lifts to a map
$\lambda$ with the same properties on $Q$.
\begin{proposition}\label{prop:permittedblocks}

Let $Q$ consist of a single directed path or,
up to reversal of all arrows, be $Y(a;b,c)$
with $\min\{a,b,c\}=1$.
Then there exists a linear map $\lambda\in\nq^*$
such that $\lambda(p)\in\{0,1\}$ for every
$p\in\mathcal P(Q)$ and
$\rank\mathcal K_Q(\lambda)=m(Q)$.
In particular, $\gamma(Q)=m(Q)$.
\end{proposition}
\begin{proof}
If $Q$ consists of a single directed path $P$, construct $R_P(Q)$. If $\ell(P)=1$, $R_P(Q)$ is an isolated vertex. If $\ell(P) > 1$, $R_P(Q)$ consists of the isolated starting
vertex of $P$ and the subquiver obtained from $Q$
by deleting the starting and terminal vertices of $P$. Induction on the number of vertices, using the lifting statement after Theorem~\ref{thm:reduction}, proves the assertion.

For the trivalent blocks, we argue by induction on $a+b+c$.
By interchanging $b$ and $c$ if necessary,
we may assume that $b\le c$. 
Suppose first that $b=1$.
Let $P$ be the concatenation of the path of length $a$
ending at the trivalent vertex and the path of length one starting there.
Then $P$ is a reducible terminal-path.
Set $R=R_P(Q)$.
If $c=1$, then $R$ is the union of two directed paths
of lengths one and $a$, whose only common vertex
is their terminal vertex.
These two paths are precisely the blocks of $R$.
By the preceding case and the proof of Lemma~\ref{lem:blocks},
there exists $\lambda_R\in\mathfrak n_R^*$
such that $\lambda_R(p)\in\{0,1\}$ for every
$p\in\mathcal P(R)$ and
$\rank\mathcal K_R(\lambda_R)=m(R)$.
The extension described after
Theorem~\ref{thm:reduction} gives
$\lambda\in\nq^*$ with the same properties for $Q$.
If $c\ge2$, reversing all arrows of $R$ gives
$Y(c-1;a,1)$.
Since $(c-1)+a+1<a+1+c$, the induction hypothesis
applies.
Reversing the arrows back and using the extension
described after Theorem~\ref{thm:reduction}
proves the assertion for $Q$.

It remains to consider $a=1$, $b,c\ge2$. Let $P$ be the concatenation of the arrow ending at
the trivalent vertex and the path of length $b$
starting there.
Then $P$ is a reducible terminal-path.
Set $R=R_P(Q)$.
The trivalent vertex of $Q$ becomes a source in $R$,
and $R$ has two blocks.
One consists of a single directed path of length $b-1$;
reversing all arrows of the other gives $Y(c-1;1,1)$.
Both blocks have fewer arrows than $Q$.
The preceding directed path case, the induction hypothesis,
and the proof of Lemma~\ref{lem:blocks}
therefore give the assertion for $R$.
The extension described after
Theorem~\ref{thm:reduction} proves it for $Q$,
preserving the values in $\{0,1\}$ on paths.
Finally, reversing all arrows and assigning the same
values to the corresponding paths changes the
 matrix $\mathcal K_Q(\lambda)$ only by multiplication by $-1$ and
permutations of rows and columns, so its rank is unchanged.

\end{proof}

The necessary block restrictions and the preceding construction
now give the complete characterization of full row rank.

\begin{theorem}\label{thm:fullclassification}
For every finite acyclic quiver $Q$, the following are equivalent:
\begin{enumerate}
\item $\gamma(Q)=m(Q)$;
\item $\rho(Q)=m(Q)$, that is, $Q$ has a forest-matching cover;
\item $U(Q)$ is a forest, every internal vertex has
degree at most three in $Q$, and every connected component
of $U(Q)[I(Q)]$ contains at most one vertex of degree
three in $Q$; any such vertex is adjacent to a source
or sink of $Q$;
\item $U(Q)$ is a forest and every block consists of a single directed path or, up to reversal of all arrows, is $Y(a;b,c)$ with $\min\{a,b,c\}=1$.
\end{enumerate}
If these conditions hold, there exists a linear map $\lambda\in\nq^*$
such that $\lambda(p)\in\{0,1\}$ for every
$p\in\mathcal P(Q)$ and
$\rank\mathcal K_Q(\lambda)=m(Q)$.
 The Lie algebra $\sq$ is symplectic exactly when these conditions hold and $n(Q)+m(Q)$ is even.
\end{theorem}
\begin{proof}
Theorem~\ref{thm:rho} gives (i)$\Rightarrow$(ii), and Proposition~\ref{prop:necessaryshape} gives (ii)$\Rightarrow$(iii). 
For each block $B$ containing an internal vertex,
the set of its internal vertices is the vertex set of
a connected component of $U(Q)[I(Q)]$.
Each of these internal vertices has the same degree
in $B$ as in $Q$.

Condition (iii) is therefore equivalent to (iv).
Indeed, a block with no branch vertex consists of
a non-trivial directed path.
A block with one branch vertex is the union of three
non-trivial directed paths.
At least one of these paths has length one exactly when
the branch vertex is adjacent to a source or sink.
Proposition~\ref{prop:permittedblocks} and
Lemma~\ref{lem:blocks} give (iv)$\Rightarrow$(i)
and also give a linear map $\lambda\in\nq^*$
such that $\lambda(p)\in\{0,1\}$ for every
$p\in\mathcal P(Q)$ and
$\rank\mathcal K_Q(\lambda)=m(Q)$. The symplecticity assertion follows from Theorem~\ref{thm:fullrank}.
\end{proof}

The theorem replaces the generic rank condition by a criterion
on the arrows and vertices of $Q$.
A forest-matching cover is sufficient because it forces the
permitted block shapes, on which the reduction constructs a
map with full row rank.
Together with even dimension, these conditions completely
decide symplecticity of $\sq$.

The reduction gives a second form of the classification, in
which the obstruction is detected on the final quiver.
Starting from $Q$, apply terminal-path or source reductions
whenever possible. Each step removes one vertex, so the
procedure terminates.
Call a quiver final if neither reduction is possible.

\begin{lemma}\label{lem:irreducible}
If a final quiver $H$ still has an arrow, then $\gamma(H)<m(H)$ and every closed two-form on $\mathfrak s_H$ is degenerate.
\end{lemma}
\begin{proof}

Suppose that $v$ has degree one in $F_{\max}(H)$.
Since every maximal path starts at a source and ends at a sink,
$v$ is a source or sink of $H$.
If $v$ is a sink, write $t=v$ and let $P$ be the unique
maximal path ending at $t$.
Every vertex of $P$ other than its starting vertex has indegree one.
Thus $P$ is a reducible terminal-path, contrary to the
assumption that $H$ is final.
If $v$ is a source, the same argument with all arrows reversed
gives a source reduction, again a contradiction.
Therefore all nonisolated vertices of $F_{\max}(H)$ have degree at least two.

Since $H$ has an arrow, it has a maximal non-trivial path,
so $F_{\max}(H)$ has an edge.
Every nonisolated vertex of $F_{\max}(H)$ has degree at least two.
A finite forest containing an edge has a vertex of degree one,
so $F_{\max}(H)$ contains a cycle.
By Lemma~\ref{lem:incidence}, there exist
$c_p\in\kk$, indexed by $p\in\mathcal M(H)$ and not all zero,
such that $\sum_{p\in\mathcal M(H)}c_p\alpha_p=0$.

Let $R_p$ denote the row of $\mathcal K_H$ indexed by $p \in \mathcal M(H)$.
Since $p$ is maximal, every entry of $R_p$ in a column
indexed by a non-trivial path is zero, while its entry in the column indexed by $v\in H_0$
is $-X_p\alpha_p(e_v)$.
Consequently,
$\sum_{p\in\mathcal M(H)}(c_p/X_p)R_p=0$.
This is a non-trivial linear relation over $\kk(\mathbf X)$.
Thus $\mathcal K_H$ does not have full row rank,
and $\gamma(H)<m(H)$.

Let $\omega=\beta+d\lambda$ be any closed two-form on
$\mathfrak s_H$.
For every maximal non-trivial path $p$ in $H$,
we have $\omega(p,q)=0$ for every non-trivial path $q$ in $H$,
and $\omega(p,h)=\lambda(p)\alpha_p(h)$ for
$h\in\mathfrak h_H$.
If $\lambda(p)=0$ for some such $p$, then
$\omega(p,y)=0$ for every $y\in\mathfrak s_H$.
Since $p\ne0$, the form $\omega$ is degenerate.

Suppose now that $\lambda(p)\ne0$ for every
$p\in\mathcal M(H)$.
Using the coefficients $c_p$ chosen above, let
$x=\sum_{p\in\mathcal M(H)}(c_p/\lambda(p))p$, which is not zero.
For every $h\in\mathfrak h_H$, we have
$\omega(x,h)=\sum_{p\in\mathcal M(H)}c_p\alpha_p(h)=0$.
Also, $\omega(x,q)=0$ for every non-trivial path $q$ in $H$.
Thus $\omega(x,y)=0$ for every $y\in\mathfrak s_H$,
and $\omega$ is degenerate.
\end{proof}

\begin{theorem}\label{thm:recursive}
Starting from $Q$, apply terminal-path or source reductions
until neither is possible, and denote the final quiver
by $Q_\infty$. Then $\gamma(Q)=m(Q)$ if and only if $(Q_\infty)_1=\varnothing$. Moreover,
 $\sq$ is symplectic if and only if 
 $(Q_\infty)_1=\varnothing$ and $|(Q_\infty)_0|$ is even.
\end{theorem}
\begin{proof}

By Theorem~\ref{thm:reduction}, the difference $m-\gamma$
is preserved at every step.
Hence $m(Q)-\gamma(Q)=m(Q_\infty)-\gamma(Q_\infty)$.
If $Q_\infty$ has no arrows, then
$m(Q_\infty)=\gamma(Q_\infty)=0$.
If $Q_\infty$ has an arrow, Lemma~\ref{lem:irreducible}
gives $\gamma(Q_\infty)<m(Q_\infty)$.
This proves the first assertion.

By Theorem~\ref{thm:reduction}, each step also preserves
the existence of a symplectic structure.
Thus $\sq$ is symplectic if and only if
$\mathfrak s_{Q_\infty}$ is symplectic.
If $Q_\infty$ has an arrow, Lemma~\ref{lem:irreducible}
shows that every closed two-form on
$\mathfrak s_{Q_\infty}$ is degenerate.
If $Q_\infty$ has no arrows, then
$\mathfrak s_{Q_\infty}$ is abelian of dimension
$|(Q_\infty)_0|$, and therefore admits a symplectic
structure if and only if this number is even.
This proves the second assertion.

\end{proof}

The criterion for symplecticity is independent of all choices of reductions.
When the final quiver has no arrows and an even number of vertices,
reversing the reductions lifts a symplectic form from its abelian
Lie algebra to $\sq$.
Isolated vertices must be retained, since they contribute
to the parity of the final dimension.

Theorem~\ref{thm:fullclassification} and Theorem~\ref{thm:recursive} also 
show that a reduction preserves the existence of a
forest-matching cover.
For a terminal-path reduction $R=R_P(Q)$ with
$P=a_1\cdots a_\ell$, we can explicitly construct a
forest-matching cover of $Q$ from any such cover $(M,T)$ of $R$.
Replace each path in $V(M)\cup T$ by its corresponding path
in $Q$ under the bijection in
Lemma~\ref{lem:reductionembedding}.
Then add $P$ to $T$ and the edges
$\{a_1\cdots a_i,a_{i+1}\cdots a_\ell\}$ to $M$
for $1\le i<\ell$.
The resulting matching and path set cover all non-trivial
paths of $Q$.
The associated graph is obtained from the forest of the
original cover by adjoining the vertex $t$ and the edge
$\{s,t\}$, where $s=s(P)$ and $t=t(P)$, so it remains a forest.

Given the lists of incoming and outgoing arrows at each vertex,
the symplecticity criterion in
Theorem~\ref{thm:fullclassification} can be checked in
$O(|Q_0|+|Q_1|)$ operations, without enumerating
the non-trivial paths.
The forest, degree, and component conditions have this cost.
For the parity condition, process the vertices in a topological
order and compute
$c_v=\sum_{a\in Q_1,\ t(a)=v}(1+c_{s(a)})$ modulo two.
Since $c_v$ counts non-trivial paths terminating at $v$
modulo two, $m(Q)\equiv\sum_v c_v\pmod2$.
The reduction also gives a construction of symplectic forms.
When $S_Q$ carries a left-invariant symplectic form,
the terminal-path operation has the Hamiltonian realization
described in Section~\ref{sec:solvablegroups}.

\section{Geometry of full path Lie groups}\label{sec:solvablegroups}

For symplectic full path Lie groups, the terminal-path
operation has a geometric realization as Hamiltonian reduction.
We use this reduction to construct Lagrangian subgroups and
then compute the largest possible dimension of an isotropic
normal subgroup.
This distinguishes the existence of Lagrangian subgroups from
the stronger requirement of normality.
We conclude with the cohomology and multiplication ranks of
the complex unit group.

In this section the Lie groups and manifolds are real unless
stated otherwise, and the corresponding Lie algebras are taken
over $\mathbb R$.
By Proposition~\ref{prop:solvablegroups}, $S_Q$ is simply
connected and completely solvable, and it has no lattice
when $Q$ has an arrow.

\subsection{\texorpdfstring{Symplectic forms on $S_Q$ and upper-triangular groups}{Symplectic forms on SQ and upper-triangular groups}}\label{subsec:solvableforms}
We first apply the classification of $\sq$ to left-invariant
symplectic forms on $S_Q$.
We then consider linearly oriented quivers and determine
which real upper-triangular groups with positive diagonal
entries admit such forms.

\begin{theorem}\label{thm:solvableexistence}
The Lie group $S_Q$ admits a left-invariant symplectic form if and
only if the equivalent conditions of
Theorem~\ref{thm:fullclassification} hold and $n(Q)+m(Q)$ is even.
\end{theorem}
\begin{proof}
A left-invariant form is determined by its value at the identity,
and it is symplectic exactly when that value is a symplectic
form on $\sq$. Apply Theorem~\ref{thm:fullclassification}.
\end{proof}

Different from  the symplectic positive path Lie groups in Section~\ref{sec:groups},
symplectic full  path Lie groups occur at arbitrarily large
path length. The linearly oriented quivers already provide
such examples.

For $q\ge1$, let $A_q$ be the quiver consisting of a directed path
on vertices $1,\cdots,q$, and let $p_{ij}$ be its unique path from $i$ to $j$, $i<j$. The map $e_i\mapsto E_{ii}$, $p_{ij}\mapsto E_{ij}$
identifies $\mathfrak s_{A_q}$ with the upper-triangular Lie
algebra $\mathfrak t_q$. It also identifies $S_{A_q}$ with the
group $T_q^+$ of real upper-triangular matrices with positive
diagonal entries.
\begin{theorem}\label{thm:chains}
The Lie group $T_q^+$ admits a left-invariant symplectic form
exactly when $q\equiv0$ or $3\pmod4$.
\end{theorem}
\begin{proof}
Proposition~\ref{prop:permittedblocks} gives
$\gamma(A_q)=m(A_q)$, including $q=1$ when both are zero.
Since $\dim\mathfrak s_{A_q}=q(q+1)/2$,
Theorem~\ref{thm:fullrank} gives the assertion.
\end{proof}

\begin{proposition}\label{prop:chainintegral}
Let $Q=A_q$ if $q(q+1)/2$ is even, and otherwise add
one isolated vertex to $A_q$.
Then $S_Q$ admits a left-invariant symplectic form whose
value at the identity has integer entries and determinant one
in the path basis.
For every $\ell\ge1$, there exists such a quiver with $L(Q)=\ell$.
For $\ell\ge2$, a connected quiver with these properties exists;
no connected quiver of length one has a symplectic $S_Q$.
\end{proposition}
\begin{proof}
If an isolated vertex is added, $\mathcal K_Q$ is obtained
from $\mathcal K_{A_q}$ by adding a zero column, and
$\dim\sq=\dim\mathfrak s_{A_q}+1$.
In either case, $\mathcal K_Q$ has full row rank and
$\sq$ has even dimension.
The integer form of determinant one follows from
Proposition~\ref{prop:integral}, whose proof is independent
of this proposition. Take $q=\ell+1$ for the prescribed length.

For a connected example with $\ell\ge2$, start with $A_{\ell+1}$.
If its full path Lie algebra has odd dimension, add two vertices
$u,w$ and arrows from its starting vertex to $u$ and from $u$ to $w$.
This adds two vertices and three non-trivial paths, preserves
length $\ell$, and satisfies Theorem~\ref{thm:fullclassification}.
Apply Proposition~\ref{prop:integral} again.
If $Q$ is connected and $L(Q)=1$, symplecticity would force
$U(Q)$ to be a tree; then $n(Q)+m(Q)=2n(Q)-1$ is odd.
\end{proof}

\begin{remark}\label{rem:solvablesteps}
For every $d\ge2$, Proposition~\ref{prop:chainintegral}
gives a symplectic Lie group $S_Q$ solvable of step $d$.
Indeed, for $Q$ obtained from $A_{\ell+1}$ as above,
$S_Q$ is solvable of step $\lfloor\log_2\ell\rfloor+2$,
and its derived subgroup is nilpotent of step $\ell$.
\end{remark}

\subsection{\texorpdfstring{Symplectic reductions and isotropic subgroups of $S_Q$}{Symplectic reductions and isotropic subgroups of SQ}}\label{sec:hamiltonian}\label{subsec:fullideals}

We show that, when $S_Q$ is equipped with a left-invariant
symplectic form, the terminal-path reduction $R=R_P(Q)$
can be realized by Hamiltonian reduction.
After the normalization in Lemma~\ref{lem:normalization},
the resulting symplectic manifold is identified with
the full path Lie group $S_R$ equipped with the induced
left-invariant symplectic form.
This correspondence allows successive quiver reductions
to be realized geometrically within the class of full
path Lie groups.

We use symplectic reduction by isotropic normal subgroups
\cite{Medina}, a special case of Marsden--Weinstein
reduction \cite{MW}.
The corresponding reduction of Lie algebras is described in Lemma~\ref{lem:isotropicreduction}, and we use the
momentum convention of Definition~\ref{def:momentum}.

Choose a reducible terminal-path $P=a_1\cdots a_\ell$ and write
$R=R_P(Q)$, $q_i=a_{i+1}\cdots a_\ell$ for $0\le i<\ell$,
as in Definition~\ref{def:terminal}.
Let $I_P=\Span\{q_0,\cdots,q_{\ell-1}\}$ and $U_P=1+I_P$.

\begin{theorem}\label{thm:hamiltonianreduction}
Let $\omega=\beta+d\lambda$ be symplectic on $\sq$, and let
$\Omega$ be its left-invariant extension to $S_Q$.
The left action of $U_P$ is Hamiltonian, with momentum map
$\mu_P(g)(x)=-\lambda(\operatorname{Ad}_{g^{-1}}x)$ for $x\in I_P$.
If $\lambda$ is normalized, the identity component $H_P$ of
$\mu_P^{-1}(\mu_P(1))$ is
$U_P\rtimes\iota(S_R)$, and
$(H_P/U_P,\Omega_{\mathrm{red}})$ is isomorphic to
$(S_R,\Omega_R)$, where the reduced form satisfies
$\pi^*\Omega_{\mathrm{red}}=\Omega|_{H_P}$ for the quotient map
$\pi:H_P\to H_P/U_P$, and $\Omega_R$ is the left-invariant form
determined by $\omega_R=\iota^*\omega$ on $\mathfrak s_R$.
In particular, the reduction decreases dimension by $2\ell$.
\end{theorem}
\begin{proof}
The suffixes $q_i$ all terminate at the sink $t$.
Their mutual products vanish. A path extending one of them
on the left is another suffix, because all indegrees along
$P$ are one. Vertex brackets preserve their span.
Thus $I_P$ is an abelian ideal, and it is isotropic since
$\beta(I_P,\sq)=0$ and $d\lambda(I_P,I_P)=0$.
The group $U_P$ is closed, normal and abelian.

Let $\lambda^L$ be the left-invariant one-form defined by $\lambda$.
For the left action, $x^\#$ is the right-invariant vector field
$x^R$, and $\lambda^L_g(x^R_g)=
\lambda(\operatorname{Ad}_{g^{-1}}x)$.
Its flow consists of left translations, so
$\mathcal L_{x^R}\lambda^L=0$.
The form $\beta^L$ is pulled back from $S_Q/N_Q$;
its contraction with $x^R$ is zero.
Cartan's formula now gives
$\iota_{x^R}\Omega=-d(\lambda^L(x^R))$,
which proves the momentum map formula.
Since $I_P$ is abelian, $\mu_P$ is $U_P$-invariant.
The action is free and proper.
The independence of the fundamental vectors and nondegeneracy
of $\Omega$ imply that $\mu_P$ is a submersion.

The level set $\mu_P^{-1}(\mu_P(1))$ is the stabilizer of
$\lambda|_{I_P}$ under the action of $S_Q$ on $I_P^*$.
Hence $H_P$ has Lie algebra $I_P^{\perp_\omega}$.
For normalized $\lambda$, Lemma~\ref{lem:reductionblocks} gives
$I_P^{\perp_\omega}=I_P\oplus\iota(\mathfrak s_R)$.
The embedding in Lemma~\ref{lem:reductionembedding} is unital
and hence embeds $S_R$ in $S_Q$.
The subgroup $U_P\iota(S_R)$ is connected with Lie algebra
$I_P\oplus\iota(\mathfrak s_R)$, so it equals $H_P$.
Since $U_P\cap\iota(S_R)=\{1\}$ and $U_P$ is normal,
we have $H_P=U_P\rtimes\iota(S_R)$.
The kernel of $\Omega|_{H_P}$ consists of the tangent directions
to the $U_P$-orbits, and $\iota^*\Omega=\Omega_R$.
This proves the quotient and the dimension assertion.

\end{proof}

For any symplectic form $\omega=\beta+d\lambda$ on $\sq$,
Lemma~\ref{lem:maxforest} gives $\lambda(q_0)\ne0$.
We may therefore apply Lemma~\ref{lem:normalization}
to normalize $\lambda$, pulling back $\omega$ by an
inner automorphism and multiplying it by a nonzero scalar.
These operations preserve the existence of the reduction
and of Lagrangian subalgebras.
Thus the theorem applies to every symplectic form after this
specified normalization. In particular, the quotient is taken
from a momentum level, not from all of $S_Q$.

\begin{theorem}\label{thm:lagrangiansubgroup}
Every symplectic $S_Q$ admits a sequence of Hamiltonian
terminal-path or source reductions ending at an even-dimensional
abelian symplectic Lie group. 
Reversing the sequence constructs a left-invariant
symplectic form on $S_Q$ from the form on the final
abelian group, using Theorem~\ref{thm:reduction}.
Every left-invariant symplectic form on $S_Q$ admits a
connected closed Lagrangian subgroup, which can be constructed
along the same sequence.
\end{theorem}
\begin{proof}
Use Theorem~\ref{thm:recursive} and
Theorem~\ref{thm:hamiltonianreduction} at each step.
On the final abelian algebra choose a Lagrangian subspace.
If $L_R$ has already been chosen on the reduced algebra,
its inverse image is $\widetilde L=I_P\oplus\iota(L_R)$.
Lemma~\ref{lem:isotropicreduction} makes $\widetilde L$
a Lagrangian subalgebra for the normalized form
$\widetilde\omega=c\varphi^*\omega$.
Then $\varphi(\widetilde L)$ is Lagrangian for $\omega$.
Continue inductively.
The associated subgroup is connected and closed by the
subgroup correspondence in Section~\ref{sec:prelim}.
\end{proof}

The existence of Lagrangian subalgebras for completely solvable
symplectic Lie algebras follows from \cite{BC}; here they are
constructed along explicit reductions of the quiver. We next determine the largest possible dimension of a
connected closed isotropic normal subgroup and characterize
when a connected closed Lagrangian normal subgroup exists.
Both conclusions are expressed in terms of paths across
downward closed sets.

Let $\Omega$ be a left-invariant symplectic form on $S_Q$, with
value $\omega=\beta+d\lambda$ at the identity.
By Proposition~\ref{prop:underlyingforest}, $U(Q)$ is a forest.
For the downward closed sets of Definition~\ref{def:downward}, let
$a(Q)=\max_D\bigl|\{p\in\pp\mid s(p)\in D,\ t(p)\notin D\}\bigr|$,
where $D$ runs through the downward closed sets.
Write $\mathfrak z_Q=Z(\sq)$.
It is spanned by $\sum_{v\in C}e_v$ over the components
$C$ of $U(Q)$, so $\dim\mathfrak z_Q=c(Q)$.

\begin{theorem}\label{thm:fullisotropicrank}
Write $\rank(\beta|_{\mathfrak z_Q})=2r$.
The maximum dimension of a connected closed isotropic normal
subgroup of $(S_Q,\Omega)$ is $a(Q)+c(Q)-r$.
Consequently, $(S_Q,\Omega)$ has a connected closed Lagrangian
normal subgroup exactly when
$a(Q)+c(Q)-r=\tfrac12(n(Q)+m(Q))$.
Equivalently, the symplectic rank is
$\sigma(\sq,\omega)=a(Q)+c(Q)-r$.
For connected $Q$, this maximum and the existence of such a
Lagrangian normal subgroup are independent of $\Omega$.
\end{theorem}
\begin{proof}
By the subgroup correspondence, it suffices to compute
$\sigma(\sq,\omega)$.
Every isotropic ideal is abelian by
Lemma~\ref{lem:isotropicreduction}.
Every ideal is stable under the diagonal action of $\hq$.
In a forest there is at most one path between given vertices,
so the nonzero path weights are distinct.
Consequently, every abelian ideal $I$ is a direct sum of
$I\cap\hq$ and some one-dimensional spaces $\mathbb Rp$,
where $p$ is a non-trivial path.
If $h\in I\cap\hq$ and $\alpha_p(h)\ne0$,
then $[h,p]=\alpha_p(h)p$ implies $p\in I$,
contrary to $I$ being abelian.
Hence $I\cap\hq\subseteq\mathfrak z_Q$.
Thus $I=Z_0\oplus U$, where $Z_0=I\cap\hq\subseteq\mathfrak z_Q$
and $U$ is spanned by paths.

Since $U$ is an ideal, every nonzero product of a path in $U$
with another path, in either order, belongs to $U$.
Since $U$ is abelian, no two of these paths concatenate.
Let $D$ be their set of starting vertices.
Multiplication on the left shows that $D$ is downward closed.
No terminal vertex of one of these paths can belong to $D$,
since otherwise two paths in $U$ would concatenate.
Therefore $\dim U\le a(Q)$.
Conversely, the span of all paths starting in a downward
closed set $D$ and terminating outside $D$ is an abelian ideal.
Every nonzero product with another path, in either order,
still starts in $D$ and terminates outside $D$,
and no two paths in this span can concatenate.

For every such path ideal $U$, the restriction of $\omega$
to $U$ is zero because $\beta(U,U)=0$ and $[U,U]=0$.
Also $\omega(\mathfrak z_Q,U)=0$.
Thus $Z_0\oplus U$ is isotropic if and only if
$\beta(Z_0,Z_0)=0$.
In a $c(Q)$-dimensional alternating space of rank $2r$, the
largest isotropic subspaces have dimension $c(Q)-r$.
Both bounds can be attained simultaneously, giving the formula.
For connected $Q$ the center is one-dimensional, so $r=0$.
\end{proof}

The maximum dimension of a connected closed isotropic normal
subgroup depends on $\Omega$ only through the restriction
of $\omega$ to the center $\mathfrak z_Q=Z(\sq)$.
The exact part $d\lambda$ of $\omega$ does not affect
this maximum.
For connected $Q$, the formula is therefore determined entirely
by the quiver.

\begin{remark}\label{rem:cutalgorithm}
The number $a(Q)$ can be computed from arrow weights.
Let $p(v)$ and $q(v)$ count the paths ending and starting
at $v$, including the trivial path.
Compute $p(v)$ in a topological order and $q(v)$
in the reverse topological order.
Assign each arrow from $u$ to $v$ the weight $p(u)q(v)$.
Since every path crossing a downward closed set $D$
contains exactly one arrow leaving $D$, $a(Q)$ is the maximum
total weight of arrows leaving such a set.
On each tree, this maximum is computed recursively by
considering whether each vertex belongs to $D$,
with no arrow allowed to enter $D$ from outside.
\end{remark}

\begin{example}\label{ex:triangularrank}

The existence of a Lagrangian subgroup does not imply
the existence of a Lagrangian normal subgroup.
For the Lie groups $T_q^+$, we show more precisely that the
difference between the dimension of a Lagrangian subgroup
and the largest possible dimension of a connected closed
isotropic normal subgroup is unbounded as $q$ varies.

For the directed path $A_q$, every downward closed set
consists of the first $j$ vertices and is crossed by
$j(q-j)$ paths.
Hence $a(A_q)=\lfloor q^2/4\rfloor$.
Since $A_q$ is connected, Theorem~\ref{thm:fullisotropicrank}
gives $\sigma(\mathfrak t_q,\omega)=\lfloor q^2/4\rfloor+1$
for every symplectic form $\omega$, whenever
$q\equiv0,3\pmod4$.
Thus the difference between the dimension of a Lagrangian
subgroup and the maximum dimension of an isotropic normal
subgroup is
$\Delta_q=\frac{q(q+1)}4-\lfloor q^2/4\rfloor-1$.
In particular, $\Delta_{4k}=k-1$ for $k\ge1$ and
$\Delta_{4k+3}=k$ for $k\ge0$.

For $q=7$, the group $T_7^+$ has dimension $28$,
whereas $\sigma(\mathfrak t_7,\omega)=13$.
By Theorem~\ref{thm:lagrangiansubgroup}, every left-invariant
symplectic form on $T_7^+$ admits a connected closed
Lagrangian subgroup of dimension $14$.
However, no such subgroup is normal.
Equivalently, Lagrangian subalgebras exist but Lagrangian
ideals do not.
An explicit symplectic form is $\omega=\beta+d\lambda$, where
$\lambda(E_{17})=\lambda(E_{26})=\lambda(E_{35})=1$,
$\lambda$ vanishes on all other non-trivial paths, and
$\beta=E_{11}^*\wedge E_{22}^*+E_{33}^*\wedge E_{44}^*$.
Its matrix has determinant one by the construction in
Proposition~\ref{prop:integral}.

\end{example}

\subsection{Symplectic forms and cohomology of the complex unit group}\label{subsec:complexunits}
Here the Lie algebras and algebraic forms are over $\mathbb C$.

Consider the algebraic unit group
$G_Q=(\mathbb CQ)^\times$, which need not be simply connected.
As a variety it is
$(\mathbb C^\times)^{n(Q)}\times\mathbb C^{m(Q)}$.
Its left-invariant algebraic symplectic forms satisfy the same
existence criterion as in Theorem~\ref{thm:fullclassification}.
The proposition below shows that their cohomology classes and
multiplication ranks depend only on the vertex part $\beta$.

The same existence criterion also holds for left-invariant
holomorphic symplectic forms on the simply connected complex
Lie group with Lie algebra $\sq$.

Define $t_v:G_Q\to\mathbb C^\times$ by
$t_v(g)=a_v$ for
$g=\sum_{u\in Q_0}a_u e_u+\sum_{p\in\pp}b_p p$,
and let $\xi_v=[dt_v/t_v]$.
For a left-invariant closed algebraic two-form $\Omega$, write
its value at the identity as $\omega=\beta+d\lambda$.
We extend $\beta$ from $\hq$ to $\sq$ by
$\beta(\nq,\sq)=0$ and use the same notation
for its left-invariant extension to $G_Q$.

\begin{proposition}\label{prop:fullcohomology}
For every finite acyclic quiver $Q$,
$H^*_{\mathrm{dR,alg}}(G_Q)\cong\Lambda^*\Span_{\mathbb C}
\{\xi_v\mid v\in Q_0\}$ as graded algebras.
Under this identification $[\Omega]=[\beta]$.
If $\rank\beta=2r$, the map
$H^k_{\mathrm{dR,alg}}(G_Q)\to H^{k+2u}_{\mathrm{dR,alg}}(G_Q)$
defined by $\alpha\mapsto\alpha[\Omega]^u$ has rank
\begin{equation*}
 \sum_{j=0}^{2r}\binom{n(Q)-2r}{k-j}
 \min\left\{\binom{2r}{j},\binom{2r}{j+2u}\right\}
\end{equation*}
for $k,u\ge0$, where $\binom ab=0$ if $b<0$ or $b>a$.
\end{proposition}
\begin{proof}

We first compute the cohomology of the left-invariant forms.
For the element $h_0$ used in Theorem~\ref{thm:closedforms},
each dual basis vector $p^*$ corresponding to a non-trivial
path has a negative integer weight.
Hence every exterior product containing such a vector
has nonzero weight.
On the subcomplex of weight $c\ne0$, the Cartan identity gives
$d(c^{-1}\iota_{h_0})+(c^{-1}\iota_{h_0})d=\id$,
so its cohomology vanishes.
The subcomplex of weight zero is $\Lambda^*\hq^*$,
and its differential is zero.
Since the inclusion
$\Lambda^*\hq^*\hookrightarrow\Lambda^*\sq^*$
preserves exterior products, it induces an isomorphism
$H^*(\sq;\kk)\cong\Lambda^*\hq^*$
of graded algebras over any field $\kk$ of characteristic zero.

Over $\mathbb C$, the left-invariant form corresponding to
$e_v^*$ is $dt_v/t_v$.
The algebraic de Rham cohomology of $\mathbb C^{m(Q)}$
is concentrated in degree zero, while that of $\mathbb C^\times$
is the exterior algebra generated by $[dt/t]$.
Since $G_Q\cong(\mathbb C^\times)^{n(Q)}\times\mathbb C^{m(Q)}$
as varieties, the K\"unneth formula shows that
$H^*_{\mathrm{dR,alg}}(G_Q)$ is the exterior algebra
generated by the classes $\xi_v=[dt_v/t_v]$.
Hence the inclusion of left-invariant algebraic forms induces
the asserted cohomology isomorphism.
Since $d\lambda$ is exact, $[\Omega]=[\beta]$.

Split $\hq$ into a $2r$-dimensional nondegenerate
subspace for $\beta$ and its radical.
On the first summand, the linear Lefschetz decomposition gives
rank $\min\{\binom{2r}{j},\binom{2r}{j+2u}\}$ in degree $j$.
The full operator is the tensor product of this map with the
identity on $\Lambda^*(\rad\beta)^*$.
The latter space has dimension $\binom{n(Q)-2r}{k-j}$
in degree $k-j$. Summing over $j$ gives the formula.

\end{proof}

The same formula computes multiplication by $[\omega]^u$ on
$H^k(\sq;\kk)$ over every field of characteristic zero.
For the real group $S_Q$, the formula applies to the cohomology
of left-invariant forms, rather than ordinary de Rham cohomology,
since $S_Q$ is contractible by
Proposition~\ref{prop:solvablegroups}.

\section{Cohomology and Lefschetz maps of compact nilmanifolds}\label{sec:nilcohomlefschetz}

For compact quotients of symplectic positive path Lie groups,
we compute the cohomology ring and the action of symplectic
classes on it.
The first Lefschetz rank is determined by the quiver, whereas
higher ranks can distinguish symplectic classes on the same
compact nilmanifold.
An explicit decomposition into graded strings determines
all Lefschetz ranks.
When the residual graph $\GR$ has no edges and there is
an active family, we determine the minimum and maximum
ranks of the Lefschetz map on $H^2$, with the minimum
expressed through a weighted perfect matching problem.
In this section the Lie groups and manifolds are real unless
stated otherwise, and the corresponding Lie algebras are taken
over $\mathbb R$.
We use real cohomology unless stated otherwise.

\subsection{Cohomology of compact nilmanifolds}
\label{subsec:nilcohomology}
We first describe a cochain algebra determined by the
paths of $Q$ that computes the cohomology of compact
quotients of $N_Q$.
When $\nq$ is symplectic, we then compute the Betti numbers
and describe the cup product by explicit representatives.

Let $M=\Gamma\backslash N_Q$ be a compact nilmanifold.

For each $p\in\pp$ let $x_p$ be a degree one generator.
The path cochain algebra
$\mathcal C_Q=(\Lambda(x_p\mid p\in\pp),d)$, with
$dx_p=-\sum_{p=uv,\ \ell(u),\ell(v)\ge1}x_u\wedge x_v$,
is the Chevalley--Eilenberg complex of $\nq$.
By Nomizu's theorem \cite{Nomizu}, this cochain algebra
computes $H^*(M)$ as a graded algebra.
For $\Gamma_Q$, the same cochain algebra over $\mathbb Q$
computes $H^*(\Gamma_Q\backslash N_Q;\mathbb Q)$.

For the rest of this section, assume that $\nq$ admits a symplectic form.

The symplecticity criterion restricts the residual components
and active families to a form whose cohomology can be computed
using finite-dimensional $\mathfrak{sl}_2$-modules \cite{FH}.

Use the chosen residual arrows $a_1,\cdots,a_b$ and the
integers $\ell_i,r_{ij},c_\square$ from
Definition~\ref{def:activedata}. Write $S_i$ for the two-dimensional
standard module of the $i$th factor of
$\mathfrak{sl}_2^{\oplus b}$. Let $V$ be the span of the arrows other than $a_i$
in each component $K_{1,\ell_i}$ of $\GR$, the paths
of length two corresponding to the edges of these
components, and the paths $b,ab,bc,abc$ for all active
arrows $b$, where $a,c$ are the unique predecessor
and successor of $b$, respectively.

\begin{proposition}\label{prop:moduledecomposition}
There is an isomorphism
$\nq\cong(\mathfrak a\ltimes V)\oplus
\mathfrak g_\square^{\oplus c_\square}$,
where $\mathfrak a=\Span\{a_1,\cdots,a_b\}$ and $V$ are
abelian, and $\mathfrak g_\square$ is the Lie algebra of
$K_{2,2}$.
The $\mathfrak a$-action on $V^*$ extends to
$\mathfrak{sl}_2^{\oplus b}$,
with
\begin{equation*}
 V^*\cong\bigoplus_i S_i^{\oplus\ell_i}
 \ \oplus\!\bigoplus_{\{i,j\}}(S_i\otimes S_j)^{\oplus r_{ij}}.
\end{equation*}
Under this identification, $a_i$ acts as the raising operator
$E_i$, where
$E_i$ has $\begin{psmallmatrix}0&1\\0&0\end{psmallmatrix}$
in the $i$th summand of $\mathfrak{sl}_2^{\oplus b}$
and zero in all other summands.
\end{proposition}

\begin{proof}
The Lie algebra isomorphism follows by a direct calculation.
On $V^*$, each component $K_{1,\ell_i}$ contributes
$\ell_i$ copies of $S_i$, and each active arrow with
boundary arrows $a_i,a_j$ contributes $S_i\otimes S_j$.
After changing the signs of some dual basis vectors,
$a_i$ acts as $E_i$, giving the stated module
decomposition and extension of the action.
\end{proof}

The preceding module decomposition allows us to compute
the Betti numbers from the coefficients of the graded
character of $\Lambda V^*$.

For $\boldsymbol\epsilon\in\mathbb Z^b$, let
$[z^q\mathbf x^{\boldsymbol\epsilon}]f$ denote the
coefficient of $z^qx_1^{\epsilon_1}\cdots x_b^{\epsilon_b}$
in $f$.
 Let
\begin{equation*}
 \chi_Q(z,\mathbf x)=
 \prod_i\bigl((1+zx_i)(1+zx_i^{-1})\bigr)^{\ell_i}
 \prod_{\{i,j\}}\prod_{\epsilon,\eta\in\{1,-1\}}
 (1+zx_i^\epsilon x_j^\eta)^{r_{ij}}.
\end{equation*}
Write $\kappa_q=\sum_{\boldsymbol\epsilon\in\{0,1\}^b}
[z^q\mathbf x^{\boldsymbol\epsilon}]\chi_Q$ and  
$P_\square(z)=1+4z+14z^2+25z^3+28z^4+25z^5+14z^6+4z^7+z^8$.
The polynomial $P_\square(z)$ is the Poincar\'e polynomial
of $\mathfrak g_\square$.

\begin{theorem}\label{thm:betti}
Since $\nq$ is symplectic, every compact quotient
$\Gamma\backslash N_Q$ has Poincar\'e polynomial
$\operatorname{Poin}_Q(z)=(1+z)^b
\bigl(\sum_q \kappa_qz^q\bigr)P_\square(z)^{c_\square}$.
\end{theorem}
\begin{proof}
The polynomial $\chi_Q$ is the character of
$\Lambda V^*$, with $z$ recording exterior degree.
Every irreducible summand of $\Lambda^qV^*$ has the form
$V_{m_1}\otimes\cdots\otimes V_{m_b}$, where $V_m$ has
highest weight $m$ and dimension $m+1$.
On $V_{m_i}$, both the kernel and the cokernel of $E_i$
are one-dimensional.
Thus each irreducible summand of $\Lambda^qV^*$
contributes $z^q(1+z)^b$ to the Poincar\'e polynomial.
The weights of $V_m$ are $m,m-2,\cdots,-m$.
Exactly one of zero and one occurs, with multiplicity one.
Therefore $\kappa_q$ counts the irreducible summands, with their
multiplicities. Each component of $\GR$ isomorphic to $K_{2,2}$
contributes $P_\square(z)$,
as computed in Appendix~\ref{app:square}.
The stated formula for the Poincar\'e polynomial
follows from the K\"unneth formula.

\end{proof}

The Betti numbers of the graph Lie algebras associated
with $K_{1,\ell}$ are computed in \cite{ACJ}.
The cohomology of graph Lie algebras is studied
further in \cite{Aldi}.
The active families contribute the modules $S_i\otimes S_j$
and the corresponding factors in $\chi_Q$.
To retain the multiplicative information, we next choose explicit
cohomology representatives and compute their products.
Let the lowering operator $F_i$ have
$\begin{psmallmatrix}0&0\\1&0\end{psmallmatrix}$
in the $i$th summand of $\mathfrak{sl}_2^{\oplus b}$
and zero in all other summands.
Extend $E_i$ and $F_i$ to derivations of degree zero
on $\Lambda V^*$.
Set $B_i=\ell_i+2\sum_jr_{ij}$ and
$\mathcal D_i=\{j(m-j+1)\mid1\le j\le m\le B_i\}$.
Define
$P_i^+=\prod_{\delta\in\mathcal D_i}(1-F_iE_i/\delta)$ and
$P_i^-=\prod_{\delta\in\mathcal D_i}(1-E_iF_i/\delta)$.
For $I\subseteq\{1,\cdots,b\}$,  let
$P_I=\prod_{i\notin I}P_i^+\prod_{i\in I}P_i^-$.
Let $\theta_i=a_i^*$.
For $I=\{i_1<\cdots<i_k\}$, write
$\theta_I=\theta_{i_1}\wedge\cdots\wedge\theta_{i_k}$.

\begin{proposition}\label{prop:ringmodel}
There is an isomorphism of graded algebras
$H^*(\mathfrak a\ltimes V)\cong(\mathcal E_Q,\star)$,
where $\mathcal E_Q=\bigoplus_I\theta_I\operatorname{im}P_I$
is the space of representatives defined by the projections $P_I$.
For $f\in\operatorname{im}P_I\cap\Lambda^qV^*$
and $g\in\operatorname{im}P_J$, the product is defined by
\begin{equation*}
 (\theta_If)\star(\theta_Jg)
 =(-1)^{q|J|}\theta_I\wedge\theta_J\wedge
 P_{I\cup J}(f\wedge g).
\end{equation*}
It is zero if $I\cap J\ne\varnothing$.
Consequently,
$H^*(M)\cong(\mathcal E_Q,\star)\otimes
H^*(\mathfrak g_\square)^{\otimes c_\square}$
as graded algebras.
\end{proposition}

\begin{proof}
On an irreducible module of highest weight $m$,
the nonzero eigenvalues of $F_iE_i$ are
$j(m-j+1)$ for $1\le j\le m$.
Every irreducible summand of $\Lambda V^*$ under the
$i$th $\mathfrak{sl}_2$ has highest weight at most $B_i$.
Thus $P_i^+$ and $P_i^-$ project onto the spaces of
highest and lowest weight vectors, respectively.
Operators associated with distinct indices commute.

The Lie algebra cohomology differential is
$d=\sum_i\theta_i\wedge E_i$.
The irreducible decomposition shows that each cohomology
class has a unique representative in $\mathcal E_Q$.
The projection $\theta_If\mapsto\theta_IP_If$ sends
exact cochains to zero and fixes these representatives.
Hence the product is obtained by taking the exterior
product and then applying this projection.
Since $f\wedge\theta_J=(-1)^{q|J|}\theta_J\wedge f$,
this gives the stated multiplication formula.

The formula for $H^*(M)$ follows from the Lie algebra
decomposition in Proposition~\ref{prop:moduledecomposition}
and the K\"unneth formula.
\end{proof}

\subsection{Lefschetz maps}\label{subsec:nillefschetz}
By the Nomizu isomorphism used in
Theorem~\ref{thm:geometricexistence}, we compute Lefschetz
maps on $M$ in $H^*(\nq)$, using a symplectic form on
$\nq$ whose induced form on $M$ represents the given class.

We first determine the Lefschetz map on $H^1(M)$.
We describe its kernel explicitly in terms of the maximal
paths in $Q$ of length at least two and show that its rank
is independent of the symplectic class.
We then decompose cohomology into graded strings to compute
the Lefschetz ranks in all degrees and determine their
dependence on the symplectic class.

Let $\nu(Q)=|\{p\in\mathcal M(Q)\mid\ell(p)\ge2\}|$.

\begin{theorem}\label{thm:firstlefschetz}
Suppose $\dim M=2d>0$, and let $\omega$ be
a symplectic form on $M$.
Then $b_1(M)=|Q_1|$ and
$\rank L_\omega^{d-1}=|Q_1|-\nu(Q)$ on $H^1(M)$.
Equivalently, this rank is the number of tree components
of $\GR$, including isolated vertices.
Let $\eta$ be a closed two-form on $\nq$ whose induced
form on $M$ represents $[\omega]$.
Then $\ker L_\omega^{d-1}
=\{[\iota_z\eta]\mid z\in Z(\nq)\cap[\nq,\nq]\}$,
where $(\iota_z\eta)(x)=\eta(z,x)$ for $x\in\nq$.

\end{theorem}
\begin{proof}

The Nomizu isomorphism identifies $L_\omega$ on $H^*(M)$
with $L_\eta$ on $H^*(\nq)$.
Since $[\eta]^d\ne0$, the form $\eta$ is symplectic.
Let $C=[\nq,\nq]$.
The arrows form a basis of $\nq/C$, so
$H^1(\nq)=\operatorname{ann}C$ has dimension $|Q_1|$.
By Poincar\'e duality, the kernel of $L_\eta^{d-1}$
is the radical of the alternating form
$(\alpha,\beta)\mapsto[\alpha\wedge\beta\wedge\eta^{d-1}]$
on $\operatorname{ann}C$.
This alternating form is a nonzero scalar multiple of
$\eta^{-1}|_{\operatorname{ann}C}$.
Its radical is
$\{\iota_x\eta\mid x\in C\cap C^{\perp_\eta}\}$,
where $C\cap C^{\perp_\eta}=\rad(\eta|_C)$.

If $L(Q)\le2$, then $C\subseteq Z(\nq)$ and $\eta(C,C)=0$.
If $L(Q)=3$, let $C_{\mathrm{act}}$ be the span of all paths
$X_i=ab_i$ and $Y_i=b_ic$ associated with the active
families.
By Lemma~\ref{lem:corenormal},
$\eta|_{C_{\mathrm{act}}}$ is nondegenerate.
Every other path of $C$ is maximal and hence central.
Hence in both cases
$\rad(\eta|_C)=Z(\nq)\cap C$.
By Proposition~\ref{prop:structure}, this space has the
$\nu(Q)$ maximal paths of length at least two as a basis.
The asserted kernel and rank follow.
Finally, maximal paths of length three are in bijection
with active arrows, while the other maximal paths of length
at least two are the edges of $\GR$.
Thus $|Q_1|-\nu(Q)=|V(\GR)|-|E(\GR)|$.
Lemma~\ref{lem:residualstars} shows that this rank equals
the number of tree components of $\GR$, including
isolated vertices.
\end{proof}

\begin{proposition}\label{prop:nilhardlefschetz}
A symplectic compact quotient $M$ has the hard Lefschetz
property if and only if $L(Q)\le1$, in which case $M$ is a torus.
\end{proposition}
\begin{proof}
If $L(Q)\ge2$, then $\nu(Q)>0$, so
Theorem~\ref{thm:firstlefschetz} excludes the property.
For $L(Q)\le1$, the group is abelian and the assertion
follows from linear Lefschetz decomposition.
\end{proof}

Failure of hard Lefschetz for nonabelian compact nilmanifolds
is classical \cite{BG}.
Theorem~\ref{thm:firstlefschetz} identifies the kernel in this
family and computes its dimension as $\nu(Q)$.

To compute the Lefschetz maps in higher degrees, we first
normalize the symplectic form so that the local calculations
can be combined by tensor products.
In the rest of this section, $\omega$ is a symplectic form on
$\nq$ representing a symplectic class on $M$.
The automorphisms used below induce cohomology isomorphisms
and need not preserve the lattice.

By Lemma~\ref{lem:corenormal}, we may assume that
$\omega=\omega_K(J)+\omega_R$, where $\omega_R$
is symplectic on $H_R$.
For the component $K_{1,\ell_i}$ with chosen central arrow $a_i$,
write $u_{i1},\cdots,u_{i\ell_i}$ for the other arrows and
choose signs of the corresponding paths $z_{ij}$ so that
$[a_i,u_{ij}]=z_{ij}$.
Define $B^{(i)}_{jk}=\omega_R(z_{ij},u_{ik})$ and
$k^{(i)}_j=\omega_R(z_{ij},a_i)$.
For an isolated $a_i$, let $B^{(i)}$ be the empty matrix.

\begin{lemma}\label{lem:starnormal}
Each matrix $B^{(i)}$ is symmetric and has corank at most one.
Let $D=\{i\mid\dim\ker B^{(i)}=1\}$.
For $i\in D$, let $u_i$ span the radical of $B^{(i)}$
in $\Span\{u_{i1},\cdots,u_{i\ell_i}\}$ and let $z_i=[a_i,u_i]$.
There is a Lie algebra automorphism $\Phi$ such that,
after replacing $\omega$ by $\Phi^*\omega$, the following
assertions hold:
$\nq=(\mathfrak h\ltimes W)\oplus
\mathfrak g_\square^{\oplus c_\square}$ is a symplectic
direct sum, where
$\mathfrak h=
\bigoplus_{i\in D}\mathfrak h_3(a_i,u_i,z_i)
\oplus\bigoplus_{i\notin D}\mathbb R a_i$,
each $\mathfrak h_3(a_i,u_i,z_i)$ has $[a_i,u_i]=z_i$
as its only nonzero bracket, and
$[W,W]=[u_i,W]=[z_i,W]=0$ for $i\in D$.
Here $W$ is the abelian ideal obtained from $K$ and the
spaces $U_i'\oplus[a_i,U_i']$, where $U_i'$ is a complement
to the radical of $B^{(i)}$ in the span of the $u_{ij}$.
Moreover,
\begin{equation*}
W^*\cong
\bigoplus_i S_i^{\oplus(\ell_i-\mathbf1_D(i))}
\oplus
\bigoplus_{\{i,j\}}(S_i\otimes S_j)^{\oplus r_{ij}}.
\end{equation*}
The restriction of $\omega$ to $\mathfrak h\ltimes W$ is
$\Omega_W+\sum_{i\in D}\theta_i\wedge\zeta_i+\beta$,
where $\Omega_W=\omega|_W$ is nondegenerate,
$\theta_i=a_i^*$ and $\zeta_i=-z_i^*$, so
$d\zeta_i=\theta_i\wedge u_i^*$.
Let
$U=\Span\bigl(\{u_i\mid i\in D\}
\cup\{a_i\mid i\notin D\}\bigr)$.
Then $\beta\in\Lambda^2U^*$ is nondegenerate,
with $\Omega_W$ and $\beta$ extended by zero to the
other summands.
\end{lemma}

\begin{proof}
For one component, write $a=a_i$, $u_j=u_{ij}$, $z_j=z_{ij}$,
$B=B^{(i)}$, $k=k^{(i)}$ and $\ell=\ell_i$.
Closedness on $(a,u_j,u_k)$ gives $B_{jk}=B_{kj}$.
The vectors $z_j$ are central and satisfy
$\omega_R(z_j,z_k)=0$ for all $j,k$.
Moreover, $\omega_R(z_j,r)=0$ for every residual arrow
$r$ outside the component $K_{1,\ell}$ under consideration.
Thus $[k\mid B]$ has row rank $\ell$, so $B$ has
corank at most one.

For each component $K_{2,2}$ of $\GR$, the four central
paths are orthogonal to all paths outside the corresponding
subalgebra.
The matrix of $\omega$ between these four central paths
and the four residual arrows corresponding to the vertices
of this component is therefore nonsingular.
Adding suitable linear combinations of these central paths
first to the residual arrows outside this component and then
to the four residual arrows inside the component
makes the corresponding subalgebra orthogonal to its
complement and makes $\omega$ vanish on pairs of arrows
in this component.
We now consider the remaining summand.

Let $U_0=\Span\{u_1,\cdots,u_\ell\}$ and regard $B$
as a symmetric bilinear form on $U_0$.
If $B$ is nondegenerate, let $U_0'=U_0$.
Otherwise choose $U_0'$ such that
$U_0=U_0'\oplus\rad B$ and $B|_{U_0'}$ is nondegenerate.
In the latter case, choose $0\ne u\in\rad B$ and let
$z=[a,u]$.
Then $\omega_R(u,z)=0$ and $\omega_R(a,z)\ne0$.

Define $W$ as the direct sum of $K$ and the subspaces
$U_0'\oplus[a,U_0']$ obtained from all components
$K_{1,\ell}$.
Then $W$ is an abelian ideal.
The matrix of $\omega$ on $[a,U_0']\times U_0'$
is nonsingular.
The spaces $[a,U_0']$ are mutually orthogonal,
and each is orthogonal to $U_0'$ from every other component.
Since $\omega(K,H_R)=0$ and $\omega|_K$ is nondegenerate,
the restriction $\omega|_W$ is nondegenerate,
regardless of the values of $\omega$ between the
different subspaces $U_0'$.

As in Lemma~\ref{lem:corenormal}, an automorphism replaces
the complementary subalgebra by $W^{\perp_\omega}$,
so that $\omega$ vanishes between $W$ and its complement
in the remaining summand.

Let $\Omega_W$ agree with $\omega|_W$ except that its
values on all pairs of vectors in the spaces $U_0'$ are zero.
Define $T$ by
$\omega(x,y)=\Omega_W(x,(\id_W+T)y)$ for $x,y\in W$.
Then $T$ maps the spaces $U_0'$ into the spaces $[a,U_0']$,
vanishes on the latter and on $K$, and commutes with
the action of the complementary subalgebra.
Thus $T^2=0$, and $\id_W-\frac12T$, extended by the identity
on the complement, is an automorphism transforming
$\omega|_W$ into $\Omega_W$.

The complement is the Lie algebra stated in the lemma.
For each $i\in D$, multiply $u_i$ and $z_i$ by the
same nonzero scalar so that the coefficient of
$\theta_i\wedge\zeta_i$ is one.
Add suitable multiples of $z_i$ to the generators
$a_j,u_j$ other than $a_i$ so that they are orthogonal
to $a_i$, without changing the brackets or the action on $W$.
The remaining form is $\beta$.
Since $\Omega_W$ and
$\sum_{i\in D}\theta_i\wedge\zeta_i$ are nondegenerate
on their respective subspaces, $\omega$ is nondegenerate
if and only if $\beta$ is nondegenerate.
\end{proof}

On the summands corresponding to $K_{1,\ell}$,
$\Omega_W$ is the tensor product of the standard
alternating form on $S_i$ and a symmetric form on
the multiplicity space.
On the summands corresponding to an active family,
$\Omega_W$ is the tensor product of the standard
alternating forms on $S_i,S_j$ and the alternating
form with matrix $J$.
Hence $\Omega_W$ is invariant under
$\mathfrak{sl}_2^{\oplus b}$, and the Lefschetz
decomposition of $\Lambda W^*$ is preserved by
this action.

This compatibility allows the primitive decomposition to be
combined with the coefficient cohomology of the complementary
subalgebra.

The complementary subalgebra is a direct sum of
three-dimensional Heisenberg Lie algebras and
one-dimensional abelian Lie algebras.
The next lemma computes their coefficient
cohomology together with the multiplication operators needed
to assemble the global Lefschetz action.
Let $\mathfrak h_3=\Span\{a,u,z\}$ be the
three-dimensional Heisenberg Lie algebra,
with $[a,u]=z$ and $z$ central.
Let $V_p$ be the irreducible $\mathfrak{sl}_2$-module
of highest weight $p$, regarded as an $\mathfrak h_3$-module
on which $a$ acts by the raising operator and
$u,z$ act trivially.
Let $\theta=a^*$, $\upsilon=u^*$ and $\zeta=-z^*$,
so that $d\zeta=\theta\wedge\upsilon$.
On $H^*(\mathfrak h_3;V_p)$, let $\mathsf U$ and $\mathsf T$
denote multiplication by $[\upsilon]$ and
$[\theta\wedge\zeta]$, respectively.
Let $J_2$ be a graded space with basis $v,w$ in degrees
zero and two, respectively, equipped with the operator
defined by $v\mapsto w$ and $w\mapsto0$.
Write $\kk[j]$ for a one-dimensional space
concentrated in degree $j$.

\begin{lemma}\label{lem:localcoefficients}
The graded vector spaces $H^*(\mathfrak h_3;V_p)$,
together with the operators $\mathsf U$ and $\mathsf T$,
are described by the following table.
All exterior generators have degree one.
\begin{center}
\begin{tabular}{ccl}
\toprule
$p$ & Graded space & Operators \\
\midrule
$0$ & $(\Lambda(\upsilon)\otimes J_2)
       \oplus\kk[1]\oplus\kk[2]$
 & $\mathsf U=\upsilon\wedge$, $\mathsf T$ acts on $J_2$ \\
$1$ & $\Lambda(\upsilon,v,w)$
 & $\mathsf U=\upsilon\wedge$, $\mathsf T=(\upsilon\wedge v)\wedge$ \\
$\ge2$ & $\Lambda(\upsilon,v,w)$
 & $\mathsf U=\upsilon\wedge$, $\mathsf T=0$ \\
\bottomrule
\end{tabular}
\end{center}
For $p=0$, both operators vanish on
$\kk[1]\oplus\kk[2]$.
For the one-dimensional abelian Lie algebra
$\kk a$, with $a$ acting on $V_p$ by the
raising operator, the cohomology
$H^*(\kk a;V_p)$ has dimension one
in each of degrees zero and one.
Multiplication by $[\theta]$ is nonzero
if and only if $p=0$.
\end{lemma}

\begin{proof}
Let $E$ denote the action of $a$ on $V_p$.
Choose a basis $v_0,\cdots,v_p$ of $V_p$ such that
$Ev_j=jv_{j-1}$ for $j\ge1$ and $Ev_0=0$.
Let $B=\Lambda(\upsilon,\zeta)\otimes V_p$,
graded by exterior degree.
The Chevalley--Eilenberg complex is
$\Lambda(\theta)\otimes B$, with differential
$d\xi=\theta\wedge D\xi$ for $\xi\in B$,
where $D=E+\upsilon\partial_\zeta$.
Consequently, $H^k(\mathfrak h_3;V_p)$
is the direct sum of $\ker(D|_{B^k})$
and $\theta\operatorname{coker}(D|_{B^{k-1}})$.

The operator $D$ acts as $E$ on $B^0$ and $B^2$.
On $B^1=\upsilon V_p\oplus\zeta V_p$, it satisfies
$D(\upsilon v)=\upsilon Ev$ and
$D(\zeta v)=\zeta Ev+\upsilon v$.
These formulas give cohomology dimensions
$(1,2,2,1)$ for $p=0$ and $(1,3,3,1)$
for $p>0$.
The ranks of $\mathsf U$ from degrees zero, one, and two
are $(1,0,1)$ for $p=0$ and $(1,2,1)$ for $p>0$.
The operator $\mathsf T$ has rank one from degrees zero
to two and from degree one to three
when $p=0$ or $1$, and vanishes when $p\ge2$.

For $p=0$, the classes
$[1]$, $[\upsilon]$, $[\theta\wedge\zeta]$,
and $[\theta\wedge\upsilon\wedge\zeta]$
span the summand $\Lambda(\upsilon)\otimes J_2$.
The remaining classes $[\theta]$ and
$[\upsilon\wedge\zeta]$ are killed by both operators.

For $p=1$, identify $1,v,w,v\wedge w$ with
$[v_0]$, $[\theta v_1]$,
$[\zeta v_0-\upsilon v_1]$, and
$-[\theta\wedge\zeta\,v_1]$, respectively,
and extend the identification by multiplication
by $\upsilon$.
Direct computation shows that $\mathsf T$ corresponds
to exterior multiplication by $\upsilon\wedge v$.
For $p\ge2$, the graded dimensions and the ranks of
$\mathsf U$ give the same identification of graded vector
spaces with the operator $\mathsf U$, while $\mathsf T=0$.

Finally, the cohomology of the one-dimensional
abelian Lie algebra $\kk a$ with coefficients
in $V_p$ is computed by
$V_p\xrightarrow{E}V_p$.
Multiplication by $[\theta]$ induces the natural map
$\ker E\to\operatorname{coker}E$.
Since $\ker E=\kk v_0$ and
$\operatorname{im}E=\Span\{v_0,\cdots,v_{p-1}\}$
for $p>0$, this map is nonzero exactly when $p=0$.
\end{proof}

We now combine the primitive decomposition of $\Lambda W^*$
with the local coefficient calculations.
For the summand $\mathfrak h\ltimes W$ in
Lemma~\ref{lem:starnormal}, the required data are the
residual coranks, the primitive representation
multiplicities, and the ranks of the restrictions of $\beta$.
For each $\mathfrak g_\square$ factor, we also use the
Lefschetz strings of the restricted symplectic form,
computed as in Appendix~\ref{app:square}.
The rules in Lemma~\ref{lem:stringrules} then determine
every power of the Lefschetz operator.

Use the decomposition and notation of
Lemma~\ref{lem:starnormal}, and set $\dim W=2d_W$.
For $0\le q\le d_W$, the primitive subspace of
$\Lambda^qW^*$ consists of the forms $\xi$ satisfying
$\Omega_W^{d_W-q+1}\wedge\xi=0$.
Let $\chi_W$ be obtained from $\chi_Q$ by replacing
$\ell_i$ with $\ell_i-\mathbf1_D(i)$.
For $0\le q\le d_W$ and $\mathbf p\in\mathbb Z_{\ge0}^b$ let
\begin{equation*}
 c_{q,\mathbf p}=
 \sum_{\boldsymbol\epsilon\in\{0,1\}^b}(-1)^{|\boldsymbol\epsilon|}
 \bigl([z^q\mathbf x^{\mathbf p+2\boldsymbol\epsilon}]\chi_W
 -[z^{q-2}\mathbf x^{\mathbf p+2\boldsymbol\epsilon}]\chi_W\bigr).
\end{equation*}
The number $c_{q,\mathbf p}$ is the multiplicity
of the irreducible $\mathfrak{sl}_2^{\oplus b}$-module
with highest weight $\mathbf p$
in the primitive subspace of $\Lambda^qW^*$.
For each $\mathbf p$, write
$D_0=\{i\in D\mid p_i=0\}$, $D_1=\{i\in D\mid p_i=1\}$,
$D_2=\{i\in D\mid p_i\ge2\}$,
$O_0=\{i\notin D\mid p_i=0\}$ and
$O_1=\{i\notin D\mid p_i>0\}$.
For $1\le i\le b$, let $y_i=u_i$ if $i\in D$
and $y_i=a_i$ otherwise.
For $S\subseteq\{1,\cdots,b\}$, write
$\beta_S=(\beta(y_i,y_j))_{i,j\in S}$,
with the indices arranged in increasing order.

\begin{theorem}\label{thm:allstrings}

All graded Lefschetz strings on $H^*(\nq)$ are obtained
as follows.
For each $0\le q\le d_W$ and
$\mathbf p\in\mathbb Z_{\ge0}^b$ with $c_{q,\mathbf p}>0$,
and for each $A\subseteq D_0$, perform the following steps.
\begin{enumerate}
\item Define
$S=A\cup D_2\cup O_0$, $
e=|A\cup D_1\cup D_2\cup O_0|
  +2|D_1\cup D_2|+|O_1|$, $
2r=2|D_1|+\rank\beta_S$, and $ 
 h=|D_0\setminus A|$.
\item Take the strings in Lemma~\ref{lem:stringrules}
for an exterior algebra on a space of dimension $e$
with an alternating form of rank $2r$.
Tensor them with $(q,d_W-q+1)$ and with $|A|$ copies
of $(0,2)$, using the tensor product rule in that lemma.
\item For each $0\le s\le h$, shift all starting degrees
by $h+s$ and multiply the multiplicities by
$c_{q,\mathbf p}\binom hs$.
\end{enumerate}
Collect the resulting strings over all $q,\mathbf p,A$.
This gives the graded Lefschetz strings on
$H^*(\mathfrak h\ltimes W)$.
Finally, tensor them with the Lefschetz strings of the
restricted symplectic form on each
$\mathfrak g_\square$ factor.
The resulting strings, counted with multiplicity,
are all graded Lefschetz strings on $H^*(\nq)$.

For $k,u\ge0$, each occurrence of a string $(a,\ell)$
contributes one to
$\rank(L_\omega^u:H^k\to H^{k+2u})$
if $k-a$ is nonnegative and even
and $(k-a)/2+u<\ell$, and contributes zero otherwise.

\end{theorem}

\begin{proof}
The Lefschetz decomposition of $\Lambda W^*$ associates
a string $(q,d_W-q+1)$ with each primitive vector
of degree $q$.
The character of the primitive subspace is obtained
by subtracting the character in degree $q-2$
from that in degree $q$.
Taking successive differences between weights
$p_i$ and $p_i+2$ in each coordinate gives
$c_{q,\mathbf p}$, the multiplicity of the irreducible
representation with highest weight $\mathbf p$.

For such a representation, its coefficient cohomology
over $\mathfrak h$ is the tensor product of the spaces
with operators in Lemma~\ref{lem:localcoefficients}.
For each $i\in D_0$, choose either
$\Lambda(\upsilon_i)\otimes J_2$
or one of the summands $\kk[1]$ and $\kk[2]$.
Let $A$ record the first choice.
Among the remaining $h$ indices, choosing $\kk[2]$
at exactly $s$ indices gives a total degree shift
of $h+s$, with multiplicity $\binom hs$.
The indices in $A$ also supply $|A|$ copies of $J_2$.
The remaining exterior generators number $e$.

Let $R=A\cup D_1\cup D_2\cup O_0$.
For $i\in R$, denote by $\upsilon_i$ the exterior
generator corresponding to the $i$th coordinate
of $\beta$.
On this exterior algebra, multiplication by the
remaining part of the symplectic class is represented
by the two-form with matrix $\beta_R$, together with
$\sum_{i\in D_1}\upsilon_i\wedge v_i$.
Replacing the generators $v_i$ by
$v_i$ plus suitable linear combinations of the
$\upsilon_j$ splits off $|D_1|$ nondegenerate
two-dimensional blocks and leaves the matrix
$\beta_S$ on the generators indexed by
$S=R\setminus D_1$.
All other generators lie in the radical.
This two-form therefore has rank
$2|D_1|+\rank\beta_S=2r$.

The formula for strings on an exterior algebra,
followed by tensoring with $(q,d_W-q+1)$
and the $|A|$ strings $(0,2)$, now gives the stated
strings on $H^*(\mathfrak h\ltimes W)$,
with the shifts and multiplicities above.
By Lemma~\ref{lem:starnormal}, the
$\mathfrak g_\square$ factors are symplectic direct
summands.
The K\"unneth isomorphism and the formula for tensor
products of strings therefore give all strings
on $H^*(\nq)$.

Finally, in a string $(a,\ell)$, a basis vector
of degree $k$ exists precisely when
$j=(k-a)/2$ is an integer with $0\le j<\ell$.
Its image under $L_\omega^u$ is nonzero precisely
when $j+u<\ell$, which proves the rank formula.
\end{proof}

\begin{remark}\label{rem:finiteweights}
For fixed $Q$ and $D$, the graded Lefschetz strings on
$H^*(\mathfrak h\ltimes W)$ depend on $\beta$ only through
the ranks of the matrices $\beta_S$.
The contributions of the $\mathfrak g_\square$ factors
are computed in Appendix~\ref{app:square}.

Only finitely many $c_{q,\mathbf p}$ are nonzero.
To compute the strings in Theorem~\ref{thm:allstrings},
these multiplicities can be summed over classes
of highest weights as follows.
For $i\in D$, the classes $p_i=0$, $p_i=1$,
and $p_i\ge2$ correspond, respectively, to the
coefficient operators
$[x_i^0]-[x_i^2]$,
$[x_i^1]-[x_i^3]$, and
$[x_i^2]+[x_i^3]$.
For $i\notin D$, the classes $p_i=0$ and $p_i>0$
correspond, respectively, to
$[x_i^0]-[x_i^2]$ and $[x_i^1]+[x_i^2]$.
Applying the corresponding operator in each variable
to the primitive character
$[z^q]\chi_W-[z^{q-2}]\chi_W$
gives the total multiplicity for each choice of classes.
These identities follow by telescoping the differences
between coefficients at weights $p_i$ and $p_i+2$.
Thus this computation requires only coefficients
for which the exponent of each $x_i$ belongs
to $\{0,1,2,3\}$.

\end{remark}

\subsection{Variation and extremal ranks}\label{subsec:nilvariation}
We apply Theorem~\ref{thm:allstrings} to study the dependence
of Lefschetz ranks on the symplectic class.
An explicit example gives different ranks on $H^3$ and
$H^4$ of the same compact nilmanifold.
Variation of Lefschetz ranks on fixed compact symplectic
manifolds also occurs in \cite{UV}.
We then determine the possible ranks on $H^2$, including
their minimum and maximum, when the residual graph $\GR$
has no edges and at least one active family is present.

\begin{example}\label{ex:lefschetzvariation}

Consider the quiver
\begin{center}
\begin{tikzcd}[column sep=large, row sep=large]
s \arrow[r, "a"]
& u \arrow[r, shift left=1, "b_1"]
    \arrow[r, shift right=1, "b_2"']
    \arrow[d, "f"']
& v \arrow[r, "c"]
& t \\
& w & &
\end{tikzcd}
\end{center}
There are twelve non-trivial paths.
The ideal $K$ has basis
$b_i,ab_i,b_ic,ab_ic$ for $i=1,2$.
Let
$\Omega_K=(ab_1c)^*\wedge b_2^*
-(ab_2c)^*\wedge b_1^*
+(ab_2)^*\wedge(b_1c)^*
-(ab_1)^*\wedge(b_2c)^*$,
and define
$\omega_0=\Omega_K+a^*\wedge(af)^*+f^*\wedge c^*$
and
$\omega_1=\Omega_K+a^*\wedge(af)^*+f^*\wedge c^*
+f^*\wedge(af)^*$.
Both forms are closed and nondegenerate.

The residual graph $\GR$ consists of the edge
$\{a,f\}$ and the isolated vertex $c$.
Choose $a$ as the central vertex of the component
$K_{1,1}$ on $\{a,f\}$.
For $\omega_j$, the corresponding matrix is
$B_j=(\omega_j(af,f))$, where $j=0,1$.
Thus $B_0=(0)$ and $B_1=(-1)$ have corank one
and zero, respectively.
 Theorem~\ref{thm:allstrings} gives the following ranks of
$L^{6-k}:H^k\to H^{12-k}$.
\begin{center}
\begin{tabular}{c|rrrrr}
$k$ & $0$ & $1$ & $2$ & $3$ & $4$ \\
\hline
$\omega_0$ & $1$ & $2$ & $9$ & $30$ & $75$ \\
$\omega_1$ & $1$ & $2$ & $9$ & $32$ & $79$
\end{tabular}
\end{center}
Thus $\omega_0$ and $\omega_1$ induce the same
Lefschetz ranks on $H^1$ and $H^2$, but different ranks on
$H^3$ and $H^4$ of the same compact nilmanifold.

\end{example}

In the case where the residual graph $\GR$ has no edges,
the rank of the Lefschetz map on $H^2$ admits a simpler
description.
It is controlled by the support of the inverse residual form,
which also appears in the Lagrangian ideal criterion.

Assume that $\GR$ has no edges and that there is an active
family. Then $\nq=\mathfrak a\ltimes K$.
By Lemma~\ref{lem:corenormal}, every symplectic form can be
transformed by an automorphism into $\omega_K(J)+\beta$,
where $\beta\in\Lambda^2\mathfrak a^*$ is nondegenerate
and is extended by zero whenever one argument lies in $K$.
Let $\dim\nq=2N$.
Write $B=(\beta(a_i,a_j))_{i,j=1}^b$ and let $\Pi=B^{-1}$.
Let $K_b$ be the complete graph on $\{1,\cdots,b\}$.
Identifying each boundary arrow $a_i$ with its index $i$,
regard $G_A(Q)$ as a subgraph of $K_b$.
For each edge $e=\{i,j\}$ of $G_A(Q)$, let $r_e=r_{ij}$.
All sums over $e$ below are taken over $E(G_A(Q))$.

\begin{proposition}\label{prop:secondlefschetz}
Suppose that $\GR$ has no edges and that at least one active
family exists. For the symplectic form $\omega=\omega_K(J)+\beta$,
the rank on $H^2$ is
\begin{equation*}
 \rank L_\omega^{N-2}
 =\binom b2+\sum_e\binom{r_e}{2}
 +2\sum_{\substack{e=\{i,j\}\\\Pi_{ij}\ne0}}r_e.
\end{equation*}
\end{proposition}
\begin{proof}
By Theorem~\ref{thm:allstrings}, only primitive degrees
zero, one and two contribute.
The trivial representation in degree zero contributes
$1+\binom b2$.
Each active family $e=\{i,j\}$ in degree one contributes
$2r_e$ precisely when $\beta$ is nondegenerate on
$\Span\{a_k\mid k\notin\{i,j\}\}$.
Let $B_{\widehat{i,j}}$ be obtained from $B$
by deleting rows and columns $i,j$.
Since
$\det B_{\widehat{i,j}}=\det B\,\Pi_{ij}^2$,
this condition is equivalent to $\Pi_{ij}\ne0$.

In primitive degree two, only trivial representations
contribute, each contributing one.
For each active family, the invariant subspace of
$\Lambda^2((S_i\otimes S_j)^{\oplus r_e})$
has dimension $\binom{r_e}{2}$,
since $S_i\otimes S_j$ has an invariant symmetric form.
Distinct families give no additional invariants. Indeed, for two distinct edges, some $\mathfrak{sl}_2$ factor
acts through its standard module on one summand
and trivially on the other.
Removing the line spanned by $\omega_K(J)$ gives
$\sum_e\binom{r_e}{2}-1$.
Adding these contributions proves the formula.
\end{proof}

To determine the possible ranks, we must decide which supports
of the inverse form can occur on the active-family graph.

\begin{proposition}\label{prop:secondsupport}
Under the assumptions of Proposition~\ref{prop:secondlefschetz},
for any $F\subseteq E(G_A(Q))$, there exists a
symplectic form $\omega=\omega_K(J)+\beta$ such that  
$\Pi_{ij}\ne0$ for $\{i,j\}\in F$ and
$\Pi_{ij}=0$ for $\{i,j\}\in E(G_A(Q))\setminus F$
if and only if
$K_b\setminus(E(G_A(Q))\setminus F)$
has a perfect matching.
\end{proposition}
\begin{proof}
Let $G=K_b\setminus(E(G_A(Q))\setminus F)$.
If the required form exists, its inverse matrix $\Pi$
is invertible and supported on $G$.
A nonzero term in its Pfaffian gives a perfect
matching of $G$.
Conversely, if $G$ has a perfect matching,
the Pfaffian of an alternating matrix $P$ with
independent entries on the edges of $G$
is a nonzero polynomial.
Its product with $\prod_{\{i,j\}\in F,\ i<j}P_{ij}$
is also nonzero.
Choose real values at which this product is nonzero
and define $\beta$ to have matrix $P^{-1}$.
Then $\omega=\omega_K(J)+\beta$ is symplectic
and $\Pi=P$ has the prescribed entries.
\end{proof}

The support criterion turns the minimum rank into a
minimum-weight perfect matching problem.

\begin{proposition}\label{prop:secondextrema}
Under the assumptions of Proposition~\ref{prop:secondlefschetz},
give the edges of $K_b$ weight $r_e$ on
$G_A(Q)$ and zero elsewhere. 
Let $w_{\min}$ be the minimum among the weights of all
perfect matchings of $K_b$, where the weight of a matching is the sum of the weights of its edges.
Then the minimum rank of $L_\omega^{N-2}$ on $H^2$ is
$\binom b2+\sum_e\binom{r_e}{2}+2w_{\min}$;
its maximum is
$\binom b2+\sum_e\binom{r_e}{2}+2\sum_e r_e$.
\end{proposition}
\begin{proof}
The support graph of $\Pi$ contains a perfect matching,
so its total edge weight is at least $w_{\min}$.
The rank formula therefore gives the stated lower bound.
This bound is attained by taking $\Pi$ supported exactly
on a perfect matching of minimum weight.
The upper bound is attained when $\Pi_{ij}\ne0$ on every
edge of $G_A(Q)$, which is possible by
Proposition~\ref{prop:secondsupport}.
\end{proof}

In this case, if there is a single active family of even
size $r$ and $b=2$, all Lefschetz ranks are independent
of the symplectic form, by Theorem~\ref{thm:allstrings}.
For $r=2$, the Betti numbers are
$(1,4,13,28,46,56,46,28,13,4,1)$, and the ranks of
$L^{5-k}$ for $0\le k\le5$ are $(1,2,6,20,41,56)$.

\section{Comparison and further algebraic consequences}\label{sec:comparison}

We compare the symplecticity criteria for $\nq$ and $\sq$
and derive further algebraic consequences.

Theorem~\ref{thm:fullclassification} gives the following
criteria for the existence of symplectic structures
when $L(Q)\le2$.
\begin{proposition}\label{prop:fulllow}
If $L(Q)\le1$, then $\sq$ is symplectic if and only if $U(Q)$ is a forest and $c(Q)$ is even. If $L(Q)\le2$, then $\sq$ is symplectic if and only if $U(Q)$ is a forest, $n(Q)+m(Q)$ is even, and every internal vertex has degree pair $(1,1)$, $(1,2)$ or $(2,1)$.
\end{proposition}
\begin{proof}

When $L(Q)\le1$, there are no internal vertices and
$m(Q)=|Q_1|$.
If $U(Q)$ is a forest, then
$n(Q)+m(Q)=2|Q_0|-c(Q)$, which is even if and only if
$c(Q)$ is even.
When $L(Q)\le2$, no two internal vertices are adjacent,
as otherwise there would be a path of length three.
Theorem~\ref{thm:fullclassification} therefore reduces
the structural condition to $U(Q)$ being a forest
and every internal vertex having degree pair
$(1,1)$, $(1,2)$ or $(2,1)$.
\end{proof}

Even dimension, Hall's condition, and the forest conditions on
$U(Q)$ and $\Fmax$ do not suffice for symplecticity of $\sq$.
The next example exhibits the additional obstruction at an
internal vertex of degree four.

\begin{example}\label{ex:counterexample}
Let $Q$ have vertices $s,v,t_1,t_2,t_3$ and an arrow $a$ from $s$ to $v$ and arrows $b_i$ from $v$ to $t_i$ for $1\le i\le3$. Write $z_i=ab_i$. Then $\dim\sq=12$, and both $U(Q)$ and $\Fmax$ are forests. 
The bipartite graph $\mathcal B_Q$ has a matching
saturating all seven row vertices, consisting of the edges
$(z_1,s)$, $(z_2,t_2)$, $(z_3,t_3)$, $(b_1,t_1)$,
$(b_2,v)$, $(b_3,a)$, and $(a,b_2)$.
Here $(p,z)$ denotes the edge joining the vertex
corresponding to the row indexed by $p$ in $\mathcal K_Q$
to the vertex corresponding to the column indexed by $z$.
However, the internal vertex $v$ has degree four,
so Theorem~\ref{thm:fullclassification} shows that
$\sq$ is not symplectic.

\end{example}

We next give a sufficient condition on a chosen set of
non-trivial paths that guarantees full row rank of
$\mathcal K_Q(\lambda)$ when $\lambda$ takes the value one
on that set and zero on all other non-trivial paths.

\begin{definition}\label{def:splittinggraph}
For $S\subseteq\pp$, define the splitting graph $D_S$
with vertex set $\pp$, where $p$ and $q$ are adjacent
if and only if $pq\in S$ or $qp\in S$.
\end{definition}
\begin{proposition}\label{prop:splitcertificate}
Let $T$ be the set of isolated vertices of $D_S$. If every non-trivial connected component of $D_S$ has a unique perfect matching, $T\subseteq S$, and $F_T$ is a forest, then $\gamma(Q)=m(Q)$.
\end{proposition}
\begin{proof}

Take $\lambda(p)=1$ for $p\in S$ and $\lambda(p)=0$ otherwise.
Let $D=d\lambda|_{\nq\times\nq}$.
Its support graph is $D_S$, and its nonzero matrix entries
in the path basis are $1$ or $-1$.
On each non-trivial component, the unique perfect matching
gives exactly one nonzero term in the Pfaffian expansion,
as explained after Definition~\ref{def:supportgraph}.
Thus the Pfaffian of the corresponding block is $1$ or $-1$,
so the block is nonsingular and $\rad D=\Span T$.

The linear map $\rad D\to\hq^*$ given by
$x\mapsto d\lambda(x,\cdot)|_{\hq}$ sends each $p\in T$
to $\alpha_p$, since $\lambda(p)=1$.
It is injective by Lemma~\ref{lem:incidence},
because $F_T$ is a forest.
 Therefore $\nq\cap\rad B_\lambda=0$,
so $\mathcal K_Q(\lambda)$ has row rank $m(Q)$.
Consequently, $\gamma(Q)=m(Q)$.
\end{proof}
\begin{remark}\label{rem:splitnotnecessary}
The condition in Proposition~\ref{prop:splitcertificate} is not necessary.
For $Q=Y(1;1,1)$, Proposition~\ref{prop:permittedblocks}
gives $\gamma(Q)=m(Q)$, but no choice of $S$ satisfies it.
Indeed, the two maximal paths are isolated in every $D_S$,
so $T\subseteq S$ forces both into $S$; the three arrows
then form a component $K_{1,2}$, which has no perfect matching.
\end{remark}

Parallel arrows have opposite roles in the two classifications.
They are required for symplecticity of $\nq$ in step three,
but are excluded by symplecticity of $\sq$.

\begin{proposition}\label{prop:opposedthree}
If $L(Q)=3$ and $\nq$ is symplectic, then $\sq$ is not symplectic.
\end{proposition}

\begin{proof}
Since $L(Q)=3$, choose a path $abc$ of length three,
and let $u,v$ be the starting and terminal vertices of $b$.
Then $m_{uv}>0$.
By Theorem~\ref{thm:nilthree}, $m_{uv}$ is even,
so there are at least two arrows from $u$ to $v$.
Hence $U(Q)$ is not a forest, and
Proposition~\ref{prop:underlyingforest} implies that
$\sq$ is not symplectic.
\end{proof}

The two classifications therefore intersect only at path length
at most two, where they give the following criterion.

\begin{theorem}\label{thm:simultaneous}
Both $\nq$ and $\sq$ are symplectic if and only if $L(Q)\le2$, $U(Q)$ is a forest, every internal vertex has degree pair $(1,1)$, $(1,2)$ or $(2,1)$, and both $n(Q)$ and $m(Q)$ are even.
\end{theorem}
\begin{proof}
If both algebras are symplectic, Theorem~\ref{thm:long} gives $L(Q)\le3$ and Proposition~\ref{prop:opposedthree} excludes length three. Proposition~\ref{prop:fulllow} gives the structural conditions. Even dimensions $m(Q)$ and $n(Q)+m(Q)$ give even $n(Q)$ as well.

Conversely, Proposition~\ref{prop:fulllow} implies that
$\sq$ is symplectic.
The degree assumptions give
$d^-(v)d^+(v)\le d^-(v)+d^+(v)$ for every $v\in Q_0$.
Since $m(Q)$ is even, Theorem~\ref{thm:niltwo}
implies that $\nq$ is symplectic.
These arguments also apply when $Q$ has no arrows.
\end{proof}

Simultaneous existence does not make the two forms compatible
under restriction.
If $\nq\ne0$, 
Proposition~\ref{prop:noexact} shows that the restriction
of any closed two-form on $\sq$ to $\nq$ is degenerate.

The constructions also control the coefficients of the form.
Existence over a field of characteristic zero can always be
realized by an integral matrix of determinant one in the
original path basis.

\begin{proposition}\label{prop:integral}

Each of the Lie algebras $\nq$ and $\sq$, whenever symplectic,
admits a symplectic form whose matrix in the path
basis has integer entries and determinant one.
For $\nq$, this matrix can be chosen to have exactly one
nonzero entry in each row, equal to $1$ or $-1$.

\end{proposition}

\begin{proof}
The construction in the proof of \cite[Theorem~3]{PT}
gives a symplectic form on each graph Lie algebra satisfying
its criterion, whose matrix in the basis indexed by the
vertices and edges has exactly one nonzero entry in each row,
equal to $1$ or $-1$.
The same construction works over $\kk$, since its closedness
identities hold over the integers.

Suppose first that $\nq$ is symplectic.
Theorem~\ref{thm:long} gives $L(Q)\le3$.
If $L(Q)\le2$, use $\nq\cong\mathfrak n_{\GR}$ from
Section~\ref{sec:nil}; the preceding construction applies.
If $L(Q)=3$, apply it to
$\nq/K\cong\mathfrak n_{\GR}$ using the isomorphism
in Lemma~\ref{lem:corequotient}.
Choose each $J$ in the construction of $\omega_K$ to be a direct sum
of $\left(\begin{smallmatrix}0&1\\-1&0\end{smallmatrix}\right)$.
The proof of Theorem~\ref{thm:nilthree} then gives a symplectic
form on $\nq$ with the same matrix property.
In both cases, after reordering the path basis, its matrix
is a direct sum of two-dimensional alternating matrices
of determinant one.

Suppose now that $\sq$ is symplectic.
Theorem~\ref{thm:recursive} gives a final quiver
$Q_\infty$ with $(Q_\infty)_1=\varnothing$
and $|(Q_\infty)_0|$ even.
Start with the standard symplectic form on the abelian
Lie algebra $\mathfrak s_{Q_\infty}$ and reverse the reductions
using Theorem~\ref{thm:reduction}.
This construction preserves integer entries.
The changes of basis in Lemma~\ref{lem:reductionblocks}
have integer entries and determinant $\pm1$.
The determinant identity in the same lemma shows that
each step preserves the determinant.
Source reductions have the same properties.
Thus the matrix in the path basis has
integer entries and determinant one.
\end{proof}

\clearpage
\appendix
\section{Cohomology and Lefschetz maps of
\texorpdfstring{$\mathfrak g_\square$}{g square}}\label{app:square}
We describe the cohomology ring of the graph Lie algebra
$\mathfrak g_\square$ associated with $K_{2,2}$ and compute
its graded Lefschetz strings for a given symplectic form.
These calculations supply the contributions of the
$\mathfrak g_\square$ factors in
Section~\ref{sec:nilcohomlefschetz}.

Use the ordered path basis
$(a_1,a_2,b_1,b_2,a_1b_1,a_1b_2,a_2b_1,a_2b_2)$,
where $[a_i,b_j]=a_ib_j$ for $1\le i,j\le2$ and all
other brackets vanish except those determined by skew-symmetry.
Let $e^1,\cdots,e^8$ be the dual basis and write
$e^{i_1\cdots i_k}=e^{i_1}\wedge\cdots\wedge e^{i_k}$.
The Chevalley--Eilenberg differential satisfies
$de^1=\cdots=de^4=0$,
$de^5=-e^{13}$, $de^6=-e^{14}$,
$de^7=-e^{23}$ and $de^8=-e^{24}$.

Let $\tau=[e^{12345678}]$.
For $1\le k\le4$, let $h_{k,i}\in H^k(\mathfrak g_\square)$
be the class of the $i$th representative listed below.
These classes form a basis in each degree.
For $1\le k\le3$, let
$\{\widehat h_{k,i}\}_i$ be the basis of
$H^{8-k}(\mathfrak g_\square)$ determined by
$h_{k,i}\widehat h_{k,j}=\delta_{ij}\tau$.
The unit spans $H^0(\mathfrak g_\square)$, and $\tau$ spans
$H^8(\mathfrak g_\square)$.
Writing $b_k=\dim H^k(\mathfrak g_\square)$, the Betti numbers are
$(b_0,\cdots,b_8)=(1,4,14,25,28,25,14,4,1)$.
The tables below determine all products by Poincar\'e duality.

\subsection{Cohomology representatives}
For each degree $k$, row $i$ gives a representative of $h_{k,i}$.

\begingroup\small
\setlength{\columnsep}{7mm}
\begin{multicols}{2}
\noindent\textbf{Degree 1.}\par\smallskip
\noindent $1:\ e^{1}.$\par
\noindent $2:\ e^{2}.$\par
\noindent $3:\ e^{3}.$\par
\noindent $4:\ e^{4}.$\par
\medskip
\noindent\textbf{Degree 2.}\par\smallskip
\noindent $1:\ e^{12}.$\par
\noindent $2:\ e^{15}.$\par
\noindent $3:\ e^{16}.$\par
\noindent $4:\ e^{25}+e^{17}.$\par
\noindent $5:\ e^{26}+e^{18}.$\par
\noindent $6:\ e^{27}.$\par
\noindent $7:\ e^{28}.$\par
\noindent $8:\ e^{34}.$\par
\noindent $9:\ e^{35}.$\par
\noindent $10:\ e^{37}.$\par
\noindent $11:\ e^{45}+e^{36}.$\par
\noindent $12:\ e^{46}.$\par
\noindent $13:\ e^{47}+e^{38}.$\par
\noindent $14:\ e^{48}.$\par
\medskip
\noindent\textbf{Degree 3.}\par\smallskip
\noindent $1:\ e^{125}.$\par
\noindent $2:\ e^{126}.$\par
\noindent $3:\ e^{127}.$\par
\noindent $4:\ e^{128}.$\par
\noindent $5:\ e^{135}.$\par
\noindent $6:\ e^{136}.$\par
\noindent $7:\ e^{137}.$\par
\noindent $8:\ e^{146}.$\par
\noindent $9:\ e^{147}+e^{138}.$\par
\noindent $10:\ e^{148}.$\par
\noindent $11:\ e^{156}.$\par
\noindent $12:\ e^{237}.$\par
\noindent $13:\ e^{238}.$\par
\noindent $14:\ e^{248}.$\par
\noindent $15:\ e^{256}-e^{167}+e^{158}.$\par
\noindent $16:\ e^{267}-e^{258}-e^{178}.$\par
\noindent $17:\ e^{278}.$\par
\noindent $18:\ e^{345}.$\par
\noindent $19:\ e^{346}.$\par
\noindent $20:\ e^{347}.$\par
\noindent $21:\ e^{348}.$\par
\noindent $22:\ e^{357}.$\par
\noindent $23:\ e^{457}+e^{367}+e^{358}.$\par
\noindent $24:\ e^{467}+e^{458}+e^{368}.$\par
\noindent $25:\ e^{468}.$\par
\medskip
\noindent\textbf{Degree 4.}\par\smallskip
\noindent $1:\ e^{1256}.$\par
\noindent $2:\ e^{1257}.$\par
\noindent $3:\ e^{1258}.$\par
\noindent $4:\ e^{1267}.$\par
\noindent $5:\ e^{1268}.$\par
\noindent $6:\ e^{1278}.$\par
\noindent $7:\ e^{1356}.$\par
\noindent $8:\ e^{1357}.$\par
\noindent $9:\ e^{1367}.$\par
\noindent $10:\ e^{1456}.$\par
\noindent $11:\ e^{1457}+e^{1358}.$\par
\noindent $12:\ e^{1458}.$\par
\noindent $13:\ e^{1467}+e^{1368}.$\par
\noindent $14:\ e^{1468}.$\par
\noindent $15:\ e^{2357}.$\par
\noindent $16:\ e^{2358}.$\par
\noindent $17:\ e^{2367}-e^{1378}.$\par
\noindent $18:\ e^{2368}-e^{1478}.$\par
\noindent $19:\ e^{2378}.$\par
\noindent $20:\ e^{2458}+e^{1478}.$\par
\noindent $21:\ e^{2468}.$\par
\noindent $22:\ e^{2478}.$\par
\noindent $23:\ e^{3456}.$\par
\noindent $24:\ e^{3457}.$\par
\noindent $25:\ e^{3458}.$\par
\noindent $26:\ e^{3467}.$\par
\noindent $27:\ e^{3468}.$\par
\noindent $28:\ e^{3478}.$\par
\medskip
\end{multicols}
\endgroup

\subsection{Products and Poincar\'e duality}
The following blocks determine all products of
positive-degree basis elements whose total degree is at most four.
In the block $(k,l)$, an entry $j:u$ in row $i$ means
$h_{k,i}h_{l,j}=u$, where each bold index $\mathbf r$
in $u$ denotes $h_{k+l,r}$.
Reversed products are determined by graded commutativity.
Every product in these degrees that is neither listed
nor obtained in this way is zero.
\begingroup\small
\setlength{\columnsep}{7mm}
\begin{multicols}{2}
\noindent\textbf{Block $(1,1)$.}\par\smallskip
\noindent $1:$ $2:\mathbf{1}$.\par
\noindent $3:$ $4:\mathbf{8}$.\par
\medskip
\noindent\textbf{Block $(1,2)$.}\par\smallskip
\noindent $1:$ $4:\mathbf{1}$, $5:\mathbf{2}$, $6:\mathbf{3}$, $7:\mathbf{4}$, $9:\mathbf{5}$, $10:\mathbf{7}$, $11:2\,\mathbf{6}$, $12:\mathbf{8}$, $13:\mathbf{9}$, $14:\mathbf{10}$.\par
\noindent $2:$ $2:-\mathbf{1}$, $3:-\mathbf{2}$, $4:-\mathbf{3}$, $5:-\mathbf{4}$, $9:\mathbf{7}$, $10:\mathbf{12}$, $11:\mathbf{9}$, $12:\mathbf{10}$, $13:2\,\mathbf{13}$, $14:\mathbf{14}$.\par
\noindent $3:$ $2:-\mathbf{5}$, $3:-\mathbf{6}$, $4:-2\,\mathbf{7}$, $5:-\mathbf{9}$, $6:-\mathbf{12}$, $7:-\mathbf{13}$, $11:\mathbf{18}$, $12:\mathbf{19}$, $13:\mathbf{20}$, $14:\mathbf{21}$.\par
\noindent $4:$ $2:-\mathbf{6}$, $3:-\mathbf{8}$, $4:-\mathbf{9}$, $5:-2\,\mathbf{10}$, $6:-\mathbf{13}$, $7:-\mathbf{14}$, $9:-\mathbf{18}$, $10:-\mathbf{20}$, $11:-\mathbf{19}$, $13:-\mathbf{21}$.\par
\medskip
\noindent\textbf{Block $(1,3)$.}\par\smallskip
\noindent $1:$ $15:\mathbf{1}$, $16:-\mathbf{3}+\mathbf{4}$, $17:\mathbf{6}$, $22:\mathbf{8}$, $23:\mathbf{9}+\mathbf{11}$, $24:\mathbf{12}+\mathbf{13}$, $25:\mathbf{14}$.\par
\noindent $2:$ $11:-\mathbf{1}$, $15:-\mathbf{3}+\mathbf{4}$, $16:\mathbf{6}$, $22:\mathbf{15}$, $23:2\,\mathbf{16}+\mathbf{17}$, $24:2\,\mathbf{18}+\mathbf{20}$, $25:\mathbf{21}$.\par
\noindent $3:$ $11:-\mathbf{7}$, $15:2\,\mathbf{9}-\mathbf{11}$, $16:\mathbf{16}-\mathbf{17}$, $17:-\mathbf{19}$, $23:\mathbf{24}$, $24:\mathbf{25}+\mathbf{26}$, $25:\mathbf{27}$.\par
\noindent $4:$ $11:-\mathbf{10}$, $15:-2\,\mathbf{12}+\mathbf{13}$, $16:-\mathbf{18}+\mathbf{20}$, $17:-\mathbf{22}$, $22:-\mathbf{24}$, $23:-\mathbf{25}-\mathbf{26}$, $24:-\mathbf{27}$.\par
\medskip
\noindent\textbf{Block $(2,2)$.}\par\smallskip
\noindent $2:$ $5:-\mathbf{1}$, $6:-\mathbf{2}$, $7:-\mathbf{3}$, $10:-\mathbf{8}$, $11:-\mathbf{7}$, $12:-\mathbf{10}$, $13:-\mathbf{11}$, $14:-\mathbf{12}$.\par
\noindent $3:$ $4:\mathbf{1}$, $6:-\mathbf{4}$, $7:-\mathbf{5}$, $9:\mathbf{7}$, $10:-\mathbf{9}$, $11:\mathbf{10}$, $13:-\mathbf{13}$, $14:-\mathbf{14}$.\par
\noindent $4:$ $4:2\,\mathbf{2}$, $5:\mathbf{3}+\mathbf{4}$, $7:-\mathbf{6}$, $9:\mathbf{8}$, $10:-\mathbf{15}$, $11:2\,\mathbf{9}$, $12:-\mathbf{12}+\mathbf{13}$, $13:-2\,\mathbf{16}$, $14:-\mathbf{20}$.\par
\noindent $5:$ $5:2\,\mathbf{5}$, $6:\mathbf{6}$, $9:-\mathbf{9}+\mathbf{11}$, $10:-\mathbf{17}$, $11:2\,\mathbf{12}$, $12:\mathbf{14}$, $13:-2\,\mathbf{18}$, $14:-\mathbf{21}$.\par
\noindent $6:$ $9:\mathbf{15}$, $11:\mathbf{16}+\mathbf{17}$, $12:\mathbf{18}$, $13:-\mathbf{19}$, $14:-\mathbf{22}$.\par
\noindent $7:$ $9:\mathbf{16}$, $10:\mathbf{19}$, $11:\mathbf{18}+\mathbf{20}$, $12:\mathbf{21}$, $13:\mathbf{22}$.\par
\noindent $9:$ $12:-\mathbf{23}$, $13:-\mathbf{24}$, $14:-\mathbf{25}$.\par
\noindent $10:$ $11:\mathbf{24}$, $12:\mathbf{26}$, $14:-\mathbf{28}$.\par
\noindent $11:$ $11:2\,\mathbf{23}$, $13:\mathbf{25}-\mathbf{26}$, $14:-\mathbf{27}$.\par
\noindent $12:$ $13:\mathbf{27}$.\par
\noindent $13:$ $13:2\,\mathbf{28}$.\par
\medskip
\end{multicols}
\endgroup

Define $G=(G_{ij})$ by $h_{4,i}h_{4,j}=G_{ij}\tau$.
The matrix $G$ is symmetric and nonsingular.
Its nonzero entries with $i\le j$ are listed below,
where $(i,j):c$ means $G_{ij}=c$.
\begin{quote}\small $(1,28):1$, $(2,27):-1$, $(3,26):1$, $(4,25):1$, $(5,24):-1$, $(6,23):1$, $(7,22):-1$, $(8,21):1$, $(9,20):-1$, $(10,19):1$, $(11,18):-1$, $(12,17):1$, $(13,16):1$, $(14,15):-1$.\end{quote}

Define $T_{i;jk}$ by
$h_{2,i}h_{3,j}h_{3,k}=T_{i;jk}\tau$.
In row $i$ below, $(j,k):c$ means $T_{i;jk}=c$.
Only entries with $j<k$ are listed, since graded
commutativity gives $T_{i;kj}=-T_{i;jk}$ and $T_{i;jj}=0$.
All other entries with $j<k$ are zero.
\begingroup\small
\setlength{\columnsep}{7mm}
\begin{multicols}{2}
\noindent $1:$ $(22,25):-1$, $(23,24):3$.\par
\noindent $2:$ $(12,25):-1$, $(13,24):1$, $(14,23):-1$, $(16,21):-1$, $(17,19):-1$.\par
\noindent $3:$ $(12,24):1$, $(13,23):-1$, $(14,22):1$, $(16,20):1$, $(17,18):1$.\par
\noindent $4:$ $(7,25):1$, $(9,24):-2$, $(10,23):1$, $(15,21):-2$, $(16,19):-2$.\par
\noindent $5:$ $(7,24):-1$, $(9,23):2$, $(10,22):-1$, $(15,20):2$, $(16,18):2$.\par
\noindent $6:$ $(5,25):-1$, $(6,24):1$, $(8,23):-1$, $(11,21):1$, $(15,19):-1$.\par
\noindent $7:$ $(5,24):1$, $(6,23):-1$, $(8,22):1$, $(11,20):-1$, $(15,18):1$.\par
\noindent $8:$ $(11,17):1$, $(15,16):3$.\par
\noindent $9:$ $(3,25):-1$, $(4,24):1$, $(8,17):-1$, $(10,16):-1$, $(14,15):-1$.\par
\noindent $10:$ $(1,25):1$, $(2,24):-1$, $(8,16):-1$, $(10,15):-1$, $(11,14):-1$.\par
\noindent $11:$ $(3,24):2$, $(4,23):-2$, $(6,17):1$, $(9,16):2$, $(13,15):1$.\par
\noindent $12:$ $(3,23):-1$, $(4,22):1$, $(5,17):-1$, $(7,16):-1$, $(12,15):-1$.\par
\noindent $13:$ $(1,24):-2$, $(2,23):2$, $(6,16):1$, $(9,15):2$, $(11,13):1$.\par
\noindent $14:$ $(1,23):1$, $(2,22):-1$, $(5,16):-1$, $(7,15):-1$, $(11,12):-1$.\par
\end{multicols}
\endgroup

These data determine the full cohomology ring.
Products of total degree eight follow from the defining
relations for $\widehat h_{k,i}$ and the matrix $G$.
A product of total degree five, six or seven is determined
by its products with the basis elements in complementary degree.
Up to reordering, the triples of positive degrees summing
to eight are $(1,1,6)$, $(1,2,5)$, $(1,3,4)$,
$(2,2,4)$ and $(2,3,3)$.
By associativity and graded commutativity, the first four
are determined by the displayed products of total degree
at most four and Poincar\'e duality; the last is given by $T$.
Products of total degree greater than eight vanish.

The tables are obtained by extending the displayed
differential to $\Lambda^*\mathfrak g_\square^*$ as a
derivation of degree one.
The listed representatives are closed, and their classes
are linearly independent modulo exact cochains.
For $0\le k\le8$, the ranks of
$d:\Lambda^k\mathfrak g_\square^*
\to\Lambda^{k+1}\mathfrak g_\square^*$ are
$(0,4,10,21,21,10,4,0,0)$.
These ranks give the stated Betti numbers and hence the
Poincar\'e polynomial $P_\square(z)$ used in
Section~\ref{subsec:nilcohomology}.
Exterior multiplication, followed by reduction modulo
exact cochains, gives the four product blocks.
The coefficients of $e^{12345678}$ in the corresponding
exterior products give $G$ and $T$.

\subsection{Lefschetz maps}
Let $\omega$ be a symplectic form on $\mathfrak g_\square$.
As in the proof of Lemma~\ref{lem:starnormal}, after
replacing $\omega$ by its pullback under a Lie algebra
automorphism obtained by adding central paths to the arrows,
we may assume that $\omega$ vanishes on pairs of arrows.
The form also vanishes on pairs of central paths by
Lemma~\ref{lem:centerorth}.
Thus its matrix in the ordered path basis has the form
$\left(\begin{smallmatrix}0&B\\-B^T&0\end{smallmatrix}\right)$.
Closedness is equivalent to $B$ having the form
\begin{equation*}
 B=\begin{pmatrix}
 t_1&t_2&t_3&t_4\\
 t_3&t_4&t_5&t_6\\
 t_7&t_9&t_8&t_{11}\\
 t_9&t_{10}&t_{11}&t_{12}
 \end{pmatrix}.
\end{equation*}
Nondegeneracy is equivalent to $\det B\ne0$.

For the given form $\omega$, let $M_k$ be the matrix of
$L_\omega:H^k(\mathfrak g_\square)
\to H^{k+2}(\mathfrak g_\square)$ in the bases defined above.
The multiplication tables determine these matrices,
whose entries are linear functions of $t_1,\cdots,t_{12}$.
For $u\ge0$, let $R(k,u)$ be the rank of $L_\omega^u$
in degree $k$.
Thus $R(k,0)=b_k$ and, for $u>0$,
$R(k,u)=\rank(M_{k+2u-2}\cdots M_k)$.
Set $R(k,u)=0$ whenever $k<0$ or $k+2u>8$.
For $a\ge0$ and $\ell\ge1$, the multiplicity of the graded
Lefschetz string $(a,\ell)$ is
\begin{equation*}
 R(a,\ell-1)-R(a-2,\ell)-R(a,\ell)+R(a-2,\ell+1).
\end{equation*}
Indeed, $R(a,\ell-1)-R(a-2,\ell)$ counts the strings
starting in degree $a$ and reaching degree $a+2\ell-2$.
The difference $R(a,\ell)-R(a-2,\ell+1)$ counts those
that also reach degree $a+2\ell$.
Subtracting gives the multiplicity of $(a,\ell)$.
This determines the graded Lefschetz strings contributed
by each $\mathfrak g_\square$ factor in
Theorem~\ref{thm:allstrings}.

For example,
$\omega_t=t e^{15}+e^{17}+e^{25}+e^{18}+e^{26}+e^{28}
+e^{35}+e^{37}+e^{46}$ is closed, and its matrix in the
ordered path basis has determinant $(t-2)^2$.
Thus $\omega_t$ is symplectic for $t\ne2$.
For $t=0$ and $t=1$, the ranks of
$L_{\omega_t}^2:H^2(\mathfrak g_\square)
\to H^6(\mathfrak g_\square)$ are $11$ and $12$,
respectively.
Hence the Lefschetz ranks can depend on the symplectic form.
The tables determine all Lefschetz ranks for each
nonsingular matrix $B$ of the displayed form.
A classification of the rank sequences arising as $B$
varies is not given here.

\noindent\textbf{AI Declaration.}
The proofs of Theorem~\ref{thm:nilidealcriterion} and
Theorem~\ref{thm:fullclassification} were found by ChatGPT 6.
The authors' contribution to these proofs is to verify, organize,
and present them in a rigorous form for a human audience.
The calculations in the appendix were carried out by ChatGPT 6.
The manuscript itself was written by the authors.
ChatGPT 6 was additionally used for language refinement,
proofreading, and checking the presentation of some parts of
the manuscript.
The authors have independently verified the mathematical content
and take full responsibility for the exposition and correctness
of the manuscript.

\bibliographystyle{plain}
\bibliography{symplectic_quivers}

\end{document}